\documentclass{article}
\usepackage{amsmath,amssymb,amsthm}
\usepackage{graphicx}
\usepackage[labelsep=none]{caption}
\usepackage{float}

\newtheorem{theorem}{Theorem}[section]
\newtheorem{lemma}{Lemma}[section]
\newtheorem{remark}{Remark}[section]

\newtheorem{cor}{Corollary}[section]
\newtheorem{prop}{Proposition}[section]
\newtheorem{defn}{Definition}[section]

\title{The Four-Vertex Theorem, The Evolute, and The Decomposition of Polygons}
\author{Wiktor J. Mogilski}

\begin{document}
\date{}
\maketitle
\begin{abstract}
The Four-Vertex Theorem has been of interest ever since a discrete version appeared in 1813 due to Cauchy. Up until now, there have been many different versions of this theorem, both for discrete cases and smooth cases. In 2004 \cite{mus1}, an approach relating the discrete Four-Vertex Theorem to the evolute was published, and here we will give an overview of this paper. We then will define the notion of the decomposition of polygons, and derive some new results about how this notion affects various types of extremality. We will see that from our fresh results we can easily derive discrete Four-Vertex Theorems.
\end{abstract}

\tableofcontents

\section{Introduction}

The geometric notion of curvature has always been present in mathematics and physics. Not only does it play a special role, but it also continues to be a subject of interest to this day. It turns out that the notion of curvature only admits a ``good'' definition only for the class of ($C^2$-)smooth curves. In this thesis, we are particulary interested in the discrete analog of curvature, meaning that we want to define it strictly for polygonal curves, which are definitely not $C^2$.

We will begin by considering the history of the Four-Vertex Theorem in the smooth case. The Four-Vertex Theorem was first proved for convex smooth closed plane curves in 1909 by Syamadas Mukhopadhyaya \cite{muk}, which stated that the curvature function on convex smooth closed plane curves has at least four extreme values. His approach considered the contact of osculating circles with the curve itself. In 1912, Adolph Kneser \cite{knes} proved the more general result, that this theorem holds for all simple closed curves. He proved his result using a projective argument. In 1985, Robert Osserman \cite{oss} published a different and much simpler proof which considered the circumcircle of such curves and the way this circle intersects the curve. It turns out that his proof even generalizes to Jordan curves. To this day, there are still various papers being published on this topic, of particular interest being higher dimensions.

There are many discrete flavors of the Four-Vertex Theorem, and one of the first was considered by Augustin-Louis Cauchy in 1813 \cite{cau}. Cauchy considers two convex polygons which have corresponding sides of same length. He proves that then, either the corresponding angles are the same, or the differences between corresponding angles changes sign at least four times. Another very similar theorem was one of A. D. Aleksandrov \cite{ale}, in which he considered convex polygons with parallel sides with the assumption that no parallel translation places them inside each other. He then showed that the difference between the lengths of corresponding sides changes signs at least four times.

It is interesting that a similar result was proved in 1963 by S. Bilinski \cite{bil}, except the result does not consider two polygons. Bilinski considered an equilateral generic convex polygon which had at least four vertices and proceeded to label the angles in a cyclic fashion. He then showed that the consecutive differences of these angles change signs at least four times. Bilinski's result did require that the polygon also be obtuse, meaning that all angles must be greater than $\pi$. B. Dahlberg \cite{dah} proved a very similar result, relaxing the assumption that the polygon is obtuse. Dahlberg proved the result for all generic convex polygons.

In 1932, R. C. Bose \cite{bos} published a remarkable global result to which a nice generalization was promptly discovered. For his approach, we consider a convex polygon and we the circles passing through any three consecutive vertices of this polygon. The circle is called full if it contains all the other vertices of the curve inside and empty if it contains none of the other vertices inside. We denote by $s_-$ and $s_+$ the numbers of empty and full circles, respectively. Then, we denote by $t_-$and $t_+$ the numbers of empty and full circles passing through three pairwise non-neighboring vertices of our curve, respectively. It follows that $s_- - t_- = s_+ - t_+ = 2.$ This is a very strong result which gives the global discrete Four-Vertex Theorem as a corollary. In his book \cite{pak}, Igor Pak finds an analogue to this result which uses inscribed circles as opposed to circumcircles.

While we will consider the geometric notion of the evolute of a polygonal curve and some results in this direction, in 2000 S. Tabachnikov \cite{tab3} published a very interesting result related to this topic. Tabachnikov considered two convex parallel polygons (polygons with parallel edges), one contained in the other. Defining the notion of the relative evolute, Tabachnikov derived a result which implies the Four-Vertex Theorem in a different flavor: consider extremal edges instead of extremal vertices.

Generalizations of the Four-Vertex Theorem are to this day being considered, with many fairly recent papers by V. I. Arnold \cite{arn}, B. Dahlberg \cite{dah}, B. Wegner \cite{weg1} \cite{weg2}, and S. Tabachnikov \cite{tab3}. While in this thesis we will simply restrict ourselves the planar case, the Four-Vertex Theorem is also of interest in high dimensions. V. I. Arnold \cite{arn}, W. Blaschke \cite{bla1}\cite{bla2} and O. R. Musin \cite{mus1} are well known for their work on Four-Vertex theorems of polyhedra.

In this thesis, we will define three types of extremality: global, local, and radial. We will consider various relationships between these three notions, and we will provide proofs of Four-Vertex Theorems corresponding to each of these notions of extremality. In the next section we will give an overview of O. R. Musin's paper published in 2004 \cite{mus1}, providing all the details. His paper discusses the evolute of a polygonal curve, a rich geometric figure which gives plenty of information about our original curve. We will then show the delicate relationship between the evolute and the original curve. It turns out that there is a very special connection between the winding number of the evolute and the number of extremal vertices of the original polygon.

In the subsequent section, we will introduce the notion of the decomposition of two polygons. While this notion is not particularly new, it turns out that the impact that this notion has on the different types of extremal vertices was never investigated. We will consider this and derive some new results in this direction. To our surprise, we are able to derive our Four-Vertex Theorems almost immediately from our results, providing a new approach to the Four-Vertex Theorems.

Lastly, we will briefly discuss the notion of the gluing two polygons. While we do not have many results in this direction, this notion appears analogous to the notion of the decomposition of two polygons. We conjecture that similar results will hold as for decomposition which should provide us with another approach to the Four-Vertex Theorems.

\section{Different Types of Extremality and Four-Vertex Theorems}
\renewcommand{\thefigure}{\thesection.\arabic{figure}}

We will begin this section with a few definitions and notation. We will denote by $P$ a polygonal curve, which is just simply a piecewise linear curve, with vertices $V_{1},V_{2},...,V_{n}$, $C_{i}=C(V_{i-1}V_{i}V_{i+1})$ the circumcircle formed by the corresponding vertices, $O_{i}=O(V_{i-1}V_{i}V_{i+1})$ the center of $C_{i},$ and $R_i$ the radius of $C_i.$ When we speak of a closed polygonal curve, we will refer to it as a polygon. Also, we will restrict our consideration simply to the planar case. All indices will be taken modulo the number of vertices of the polygonal curve.

\begin{defn}A vertex $V_{i}$ is said to be positive if the left angle with respect to orientation, $\angle V_{i-1}V_{i}V_{i+1}$, is at most $\pi$. Otherwise, it is said to be negative.
\end{defn}

\begin{remark}
Observe that this definition relies on how we travel on the polygonal line. The angle, $\angle V_{i-1}V_{i}V_{i+1}$, will always be on the left side.  We will always assume, without loss of generality, that we are traveling in the positive direction (counterclockwise).
\end{remark}

For simplicity, we will set $\angle V_i=\angle V_{i-1}V_{i}V_{i+1}$ for the rest of this thesis.

The following definitions were coined by Igor Pak \cite{pak}.

\begin{defn}
We say that a polygonal curve is generic if the maximal number of vertices that lie on a circle is three and no three vertices are collinear.
\end{defn}

Observe that all regular polygons are not generic.

\begin{defn}
A polygonal curve $P$ is coherent if for any three consecutive vertices $V_{i-1},V_{i},$ and $V_{i+1}$, the center of the circle $C_{i}$ lies in the infinite cone formed by the vertices $V_{i-1},V_{i},$ and $V_{i+1}$.
\end{defn}

Note that all convex right and obtuse polygons are coherent, as well as all regular polygons.

\subsection{Different Notions of Extremality}
We will now define a few notions of extremality, one global and two local.
\subsubsection{Global Extremality}
\begin{defn}
Let $C_{ijk}$ be a circle passing through any three vertices $V_{i}$, $V_{j}$, $V_{k}$ of a polygonal curve. We say that $C_{ijk}$ is empty if it contains no other vertices of the polygonal curve in its interior, and we say that it is full if it contains all of the other vertices of the polygonal curve in its interior.
\end{defn}

\begin{defn}
We say that the circle $C_{ijk}$ is neighboring if two of the vertices are adjacent to the third, disjoint if no two vertices are adjacent, and intermediate if only one pair of the vertices is adjacent.
\end{defn}

\begin{remark}
For simplicity, we will denote a neighboring circle passing through vertices $V_{i-1},$ $V_i$ and $V_{i+1}$ by $C_i.$
\end{remark}

We will denote by $s_+,$ $t_+,$ and $u_+,$ the number of full circles which are neighboring, disjoint, and intermediate, respectively, and by $s_-$, $t_-$, and $u_-$, the number of empty circles which are neighboring, disjoint, and intermediate, respectively.

\begin{defn}
We call a neighboring full or empty circle $C_i$ an extremal circle. We refer to the corresponding vertex $V_i$ as a globally extremal vertex.
\end{defn}

We will take this opportunity to discuss the triangulation of polygons, since many of our results will be using triangulation arguments. Consider all of the empty circles passing through any three distinct points of a polygon $P$ ($s_-$, $t_-$ and $u_-$). It was proven by B. Delaunay in his 1934 paper \cite{del} that the triangles formed by each of the three points corresponding to an empty circle actually form a triangulation of the polygon $P$. This triangulation is called a \emph{Delaunay triangulation}, and is usually denoted by $DT(P)$. As we can see, this notion is very closely related to globally extremal vertices.

\begin{remark}
We take this opportunity to mention that, if we assume convexity on our polygon, then if we similarly consider all of the full circles passing through any given three points, the triangles given by these three points also form a triangulation. This is commonly known as the Anti-Delaunay triangulation.
\end{remark}

\subsubsection{Discrete Curvature and Two Types of Local Extremality}

Before we define locally extremal vertices, we must define the notion of \emph{discrete curvature}. Assume that a vertex $V_i$ is positive. We say that the curvature of the vertex $V_{i}$ is greater than the curvature at $V_{i+1}$ ($V_{i}\succ V_{i+1})$ if the vertex $V_{i+1}$ is positive and $V_{i+2}$ lies outside the circle $C_{i}$ or if the vertex $V_{i+1}$ is negative and $V_{i+2}$ lies inside the circle $C_{i}$. (See Figure 2.1)

\begin{figure}[H]
\centerline{\includegraphics[scale=0.85]{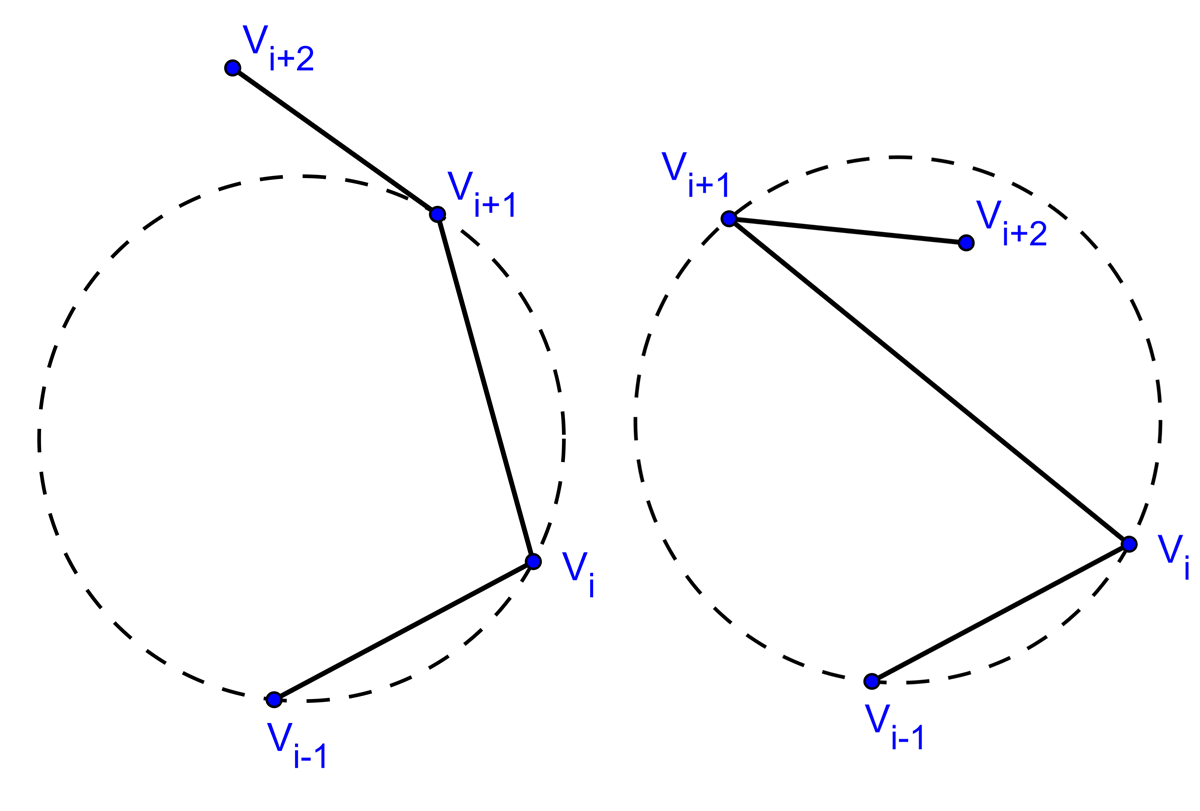}}
\caption[An Example of Discrete Curvature (1)]{}
\end{figure}

By switching the word ``inside" with the word ``outside" in the above definition (and vice-versa), we obtain that $V_{i}\prec V_{i+1}$, or that the curvature at $V_{i}$ is less than the curvature at $V_{i+1}$. (See Figure 2.2)

\begin{figure}[H]
\centerline{\includegraphics[scale=1.1]{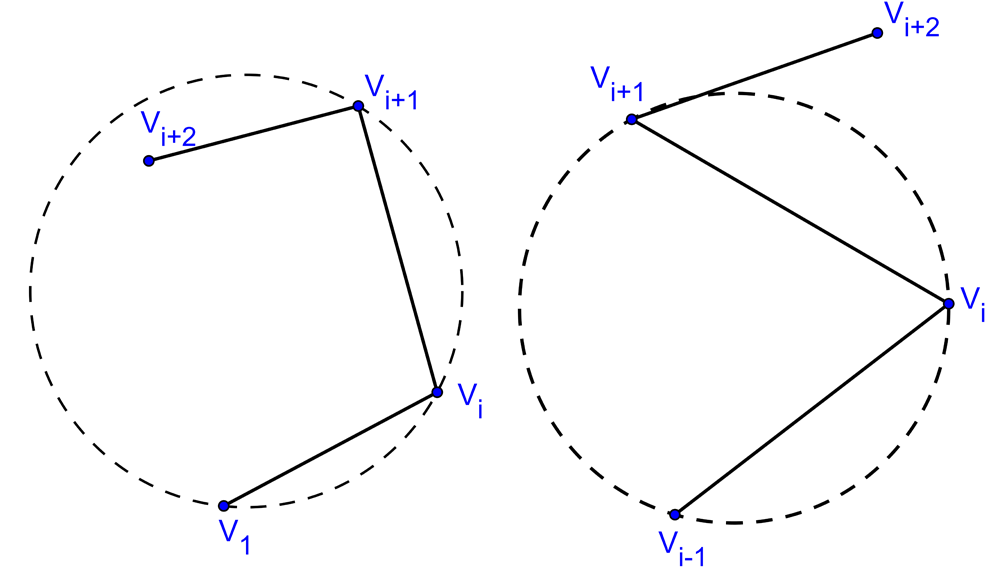}}
\caption[An Example of Discrete Curvature (2)]{}
\end{figure}

In the case that the vertex $V_{i}$ is negative, simply switch the word ``greater" with the word ``less", and the word ``outside" by the word ``inside". Figure 2.3 shows a case where $V_{i} \prec V_{i+1}$.

\begin{figure}[H]
\centerline{\includegraphics[scale=0.52]{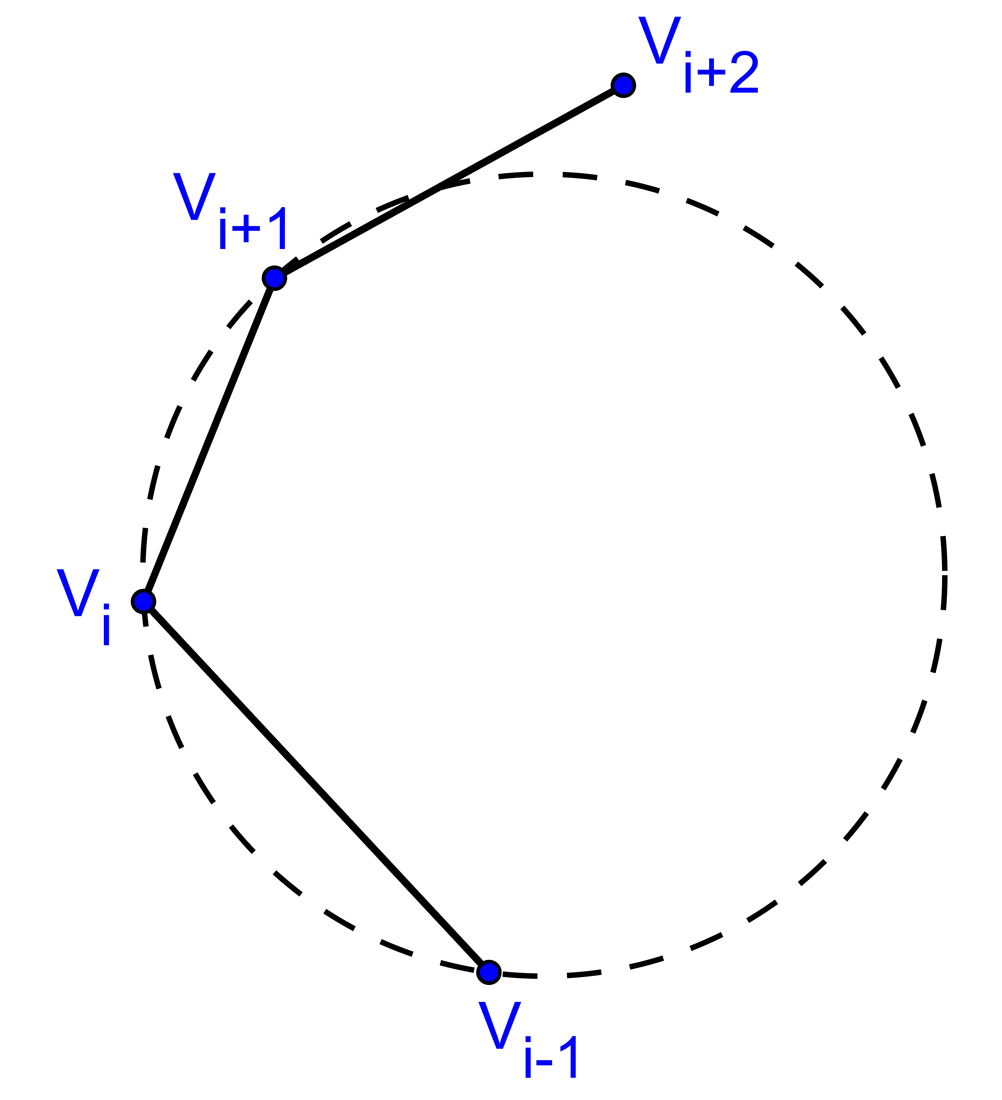}}
\caption[An Example of Discrete Curvature (3)]{}
\end{figure}

Now that we have defined the notion of discrete curvature, we can now speak about locally extremal vertices.

\begin{defn}A vertex $V_{i}$ of a polygonal line $P$ is locally extremal if

\centerline{$V_{i-1}\prec V_{i} \succ V_{i+1}$ or $V_{i-1}\succ V_{i} \prec V_{i+1}$.}

\end{defn}

\begin{remark}
 Note that in our definition of locally extremal vertices, need to consider five vertices. Hence our definition will not be interesting for a triangle. Also, observing the definition of locally extremal vertices closely, we can see that if we assume convexity on our polygon, we really are considering the position of the vertices $V_{i-2}$ and $V_{i+2}$ with respect to the circle $C_{i}.$ Our vertex $V_{i}$ will be locally extremal if they both lie inside or both lie outside the circle $C_{i}.$
\end{remark}

For a smooth curve our extremal vertices are critical points of the curvature function of the curve, positive extremal vertices are the maxima and negative extremal vertices are the minima. Also, we must note that a critical point of the curvature is simply called a vertex.
\vskip2mm
Now that we have one type of local extremality, we can define one more. For this definition, we simply consider the radii $R_{i-1},$ $R_i,$ and $R_{i+1}$ of the circles $C_{i-1},$ $C_i,$ and $C_{i+1},$ respectively.

\begin{defn}
We say that a vertex $V_{i}$ is radially extremal if $$R_{i-1}<R_i>R_{i+1}$$ or $$R_{i-1}>R_i<R_{i+1}.$$
\end{defn}

\subsection{Important Facts and Relationships Between Various Extremal Vertices}

We take this opportunity to discuss maximality and minimality of vertices. In our global sense, we mentioned empty and full neighboring circles. If a circle $C_i$ is empty, then we say that the corresponding vertex $V_i$ is maximal, and if $C_i$ is full, then we say $V_i$ is minimal. For locally extremal vertices, we call a vertex maximal if $V_{i-1}\prec V_{i} \succ V_{i+1}$ and minimal if $V_{i-1}\succ V_{i} \prec V_{i+1}.$ For radially extremal vertices, a vertex is maximal if $R_{i-1}<R_i>R_{i+1}$ and minimal if $R_{i-1}>R_i<R_{i+1}.$

From now on, we will denote the number of globally maximal-extremal vertices of a polygonal curve $P$ by $s_{-}(P)$ and globally minimal-extremal vertices by $s_{+}(P)$. Similarly, for locally extremal vertices we will attribute the notation $l_{-}(P)$ and $l_{+}(P)$ and for radially extremal vertices $r_{-}(P)$ and $r_{+}(P)$, respectively.

Before we show some relationships between our various types of extremality, we need to prove some important simple results.

\begin{prop}
\label{prop:maxmin}
Let $P$ be a generic convex polygon. Then: $$l_{+}(P)=l_{-}(P)$$ and $$r_{+}(P)=r_{-}(P).$$
\end{prop}
\begin{proof}
 For locally extremal vertices this follows directly from the fact that we need to be able to close the polygon in a convex way. For radially extremal vertices, it immediately follows from the fact that we are considering a finite set of real numbers.
\end{proof}

\begin{remark}
Note that it was very important for us to to add the assumption that our polygon is generic in our last proposition, since this eliminates the possibility of having two extremal vertices adjacent to each other.
\end{remark}

In the above proposition, observe that we did not mention globally extremal vertices. This is because the equality $s_{+}(P)=s_{-}(P)$ does not hold. In fact, we cannot form any relationship between globally maximal-extremal and globally minimal-extremal vertices. Consider the following figure:

\begin{figure}[H]
\centerline{\includegraphics{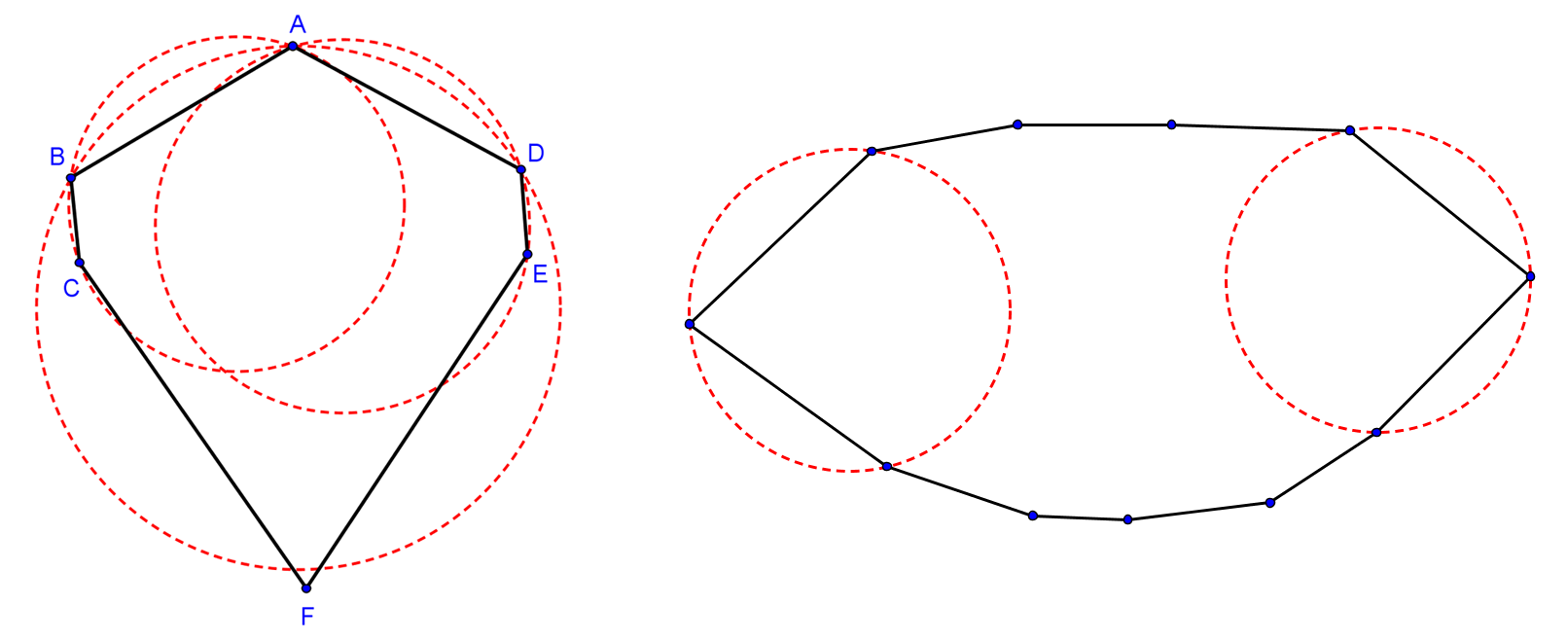}}
\caption[Examples When $s_{+}(P)\neq s_{-}(P)$]{}
\end{figure}

The polygon on the left has more globally maximal-extremal vertices than globally minimal-extremal vertices, while the polygon on the right has more globally minimal-extremal vertices than globally maximal-extremal vertices. This means that if we would like to prove results pertaining to globally extremal vertices, we will need to prove them separately for maximal-extremal and minimal-extremal vertices.

\begin{prop}
\label{prop:existenceprop}
Let $P$ be a generic convex polygon. Then $P$ has one of each type of maximal-extremal vertices.
\end{prop}

\begin{proof}
For radially extremal vertices, this is trivial, since we have a finite set of real numbers, so there must be a maximum. For locally extremal vertices this follows since we need to be able to close our polygon. For globally extremal vertices, we apply a Delaunay triangulation to $P$ and we definitely obtain a maximal-extremal circle.
\end{proof}

Now we will observe the close relationship between globally, locally, and radially extremal vertices. We will show that, if we assume convexity, there is a nice relationship between globally extremal vertices and locally extremal vertices. Moreover, with one more assumption, we will show that the notions of local and radial extremality are equivalent.

\begin{prop}
\label{prop:globaltolocal}
Let $P$ be a generic convex polygon. If $V_i$ is a globally extremal vertex, then $V_{i}$ is a locally extremal vertex.
\end{prop}

\begin{proof}
Without loss of generality, we will prove the result simply for extremal empty circles, since the proof is exactly the same for extremal full circles. Assume that we have an extremal empty circle passing through the vertices $V_{i-1},$ $V_{i}$ and $V_{i+1}$. As mentioned in Remark 2.3, it is sufficient to check the positions of the vertices $V_{i-2}$ and $V_{i+2}$ with respect to the circle $C_i$. Since we assumed convexity, all vertices are positive, and since the circle is empty, $V_{i+2}$ and $V_{i-2}$ must lie outside the circle. It follows that $V_{i-1}\prec V_{i}\succ V_{i+1}.$
\end{proof}

\begin{prop}
\label{prop:localtoradial}
Let $P$ be a generic convex coherent polygon. Then:
\begin{enumerate}
\item $V_{i-1}\succ V_{i} \iff R_{i-1}<R_{i}$
\item $V_{i-1}\prec V_{i} \iff R_{i-1}>R_{i}$
\end{enumerate}

where $R_{i-1}$ is the radius of the circle $C_{i-1}$ and $R_{i}$ is the radius of the circle $C_{i}.$
\end{prop}

\begin{proof}
The idea of the proof is as follows:
Consider the vertex $V_{i+1}.$ As this vertex moves away outside (respectively inside) of the circle $C_{i-1}$, the angle $\angle V_{i}$ increases (respectively decreases). For the (1), observe Figure 2.4:

\begin{figure}[H]
\centerline{\includegraphics[scale=1.3]{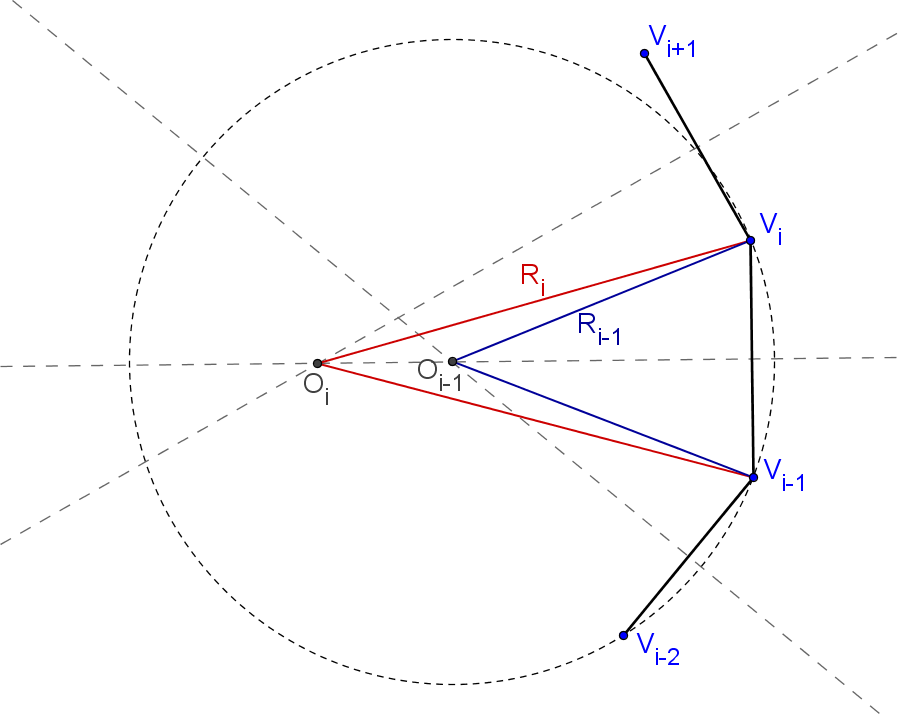}}
\caption[Proving Proposition \ref{prop:localtoradial} (1)]{}
\end{figure}

Observe that the red and blue triangles both are isosceles and share the side $V_{i}V_{i+1}.$ Now by an application of the Law of Sines, it directly follows that $R_{i-1}<R_{i}$.
\vskip2mm
For (2), we have the situation of Figure 2.5, for which we apply exactly the same argument and achieve that $R_{i-1}>R_i$.

\begin{figure}[H]
\centerline{\includegraphics[scale=1.4]{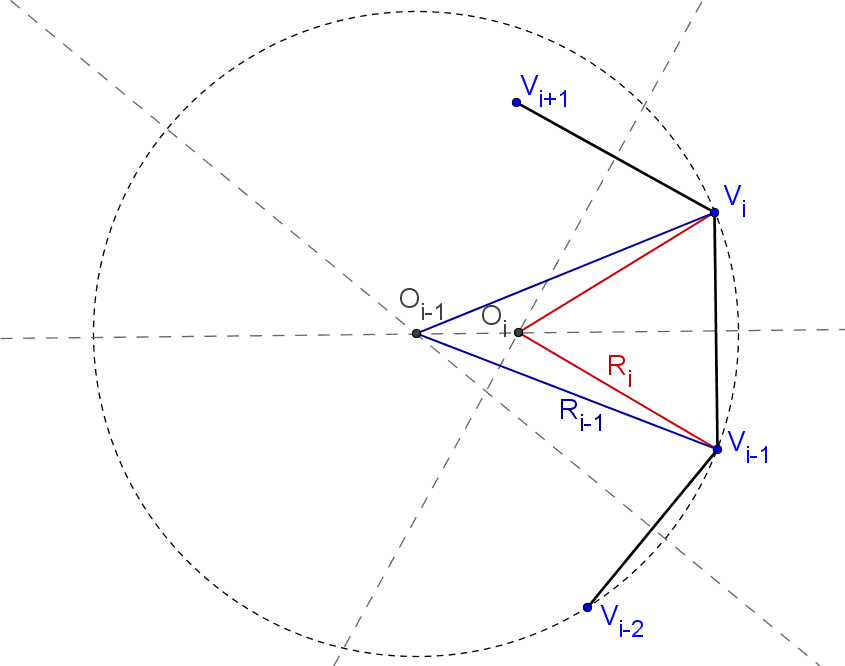}}
\caption[Proving Proposition \ref{prop:localtoradial} (2)]{}
\end{figure}
\end{proof}

\begin{remark}
Observe that we have strict inequalities. This is because of our original assumption that our polygonal curve is generic. If we were to assume that two radii are equal, then this condition will be violated. In general, we may allow such a situation, and in this case both of the vertices corresponding to the radii would be considered extremal. Also, observe that Proposition \ref{prop:localtoradial} agrees with the case of a smooth curve. If we have that $V_{i-1}\prec V_{i} \succ V_{i+1}$, then we have $R_{i-1}>R_{i}<R_{i+1}$, and similarly for the other case of locally extremal. Recall that for a smooth convex plane curve, the curvature is simply $\frac{1}{R}$, where $R$ is the radius of the osculating circle at the point of the curve.
\end{remark}

We notice that we needed to add the extra assumption of coherent to Proposition \ref{prop:localtoradial} for the proof to work. So an interesting question would be, is there any relationship between local and radial extremality if we relax this assumption? Section 2.4.1 gives us an example where we have less radially extremal vertices than locally extremal vertices. It would be nice if this was always the case, but the following example shows otherwise:

\begin{figure}[H]
\centerline{\includegraphics[scale=0.8]{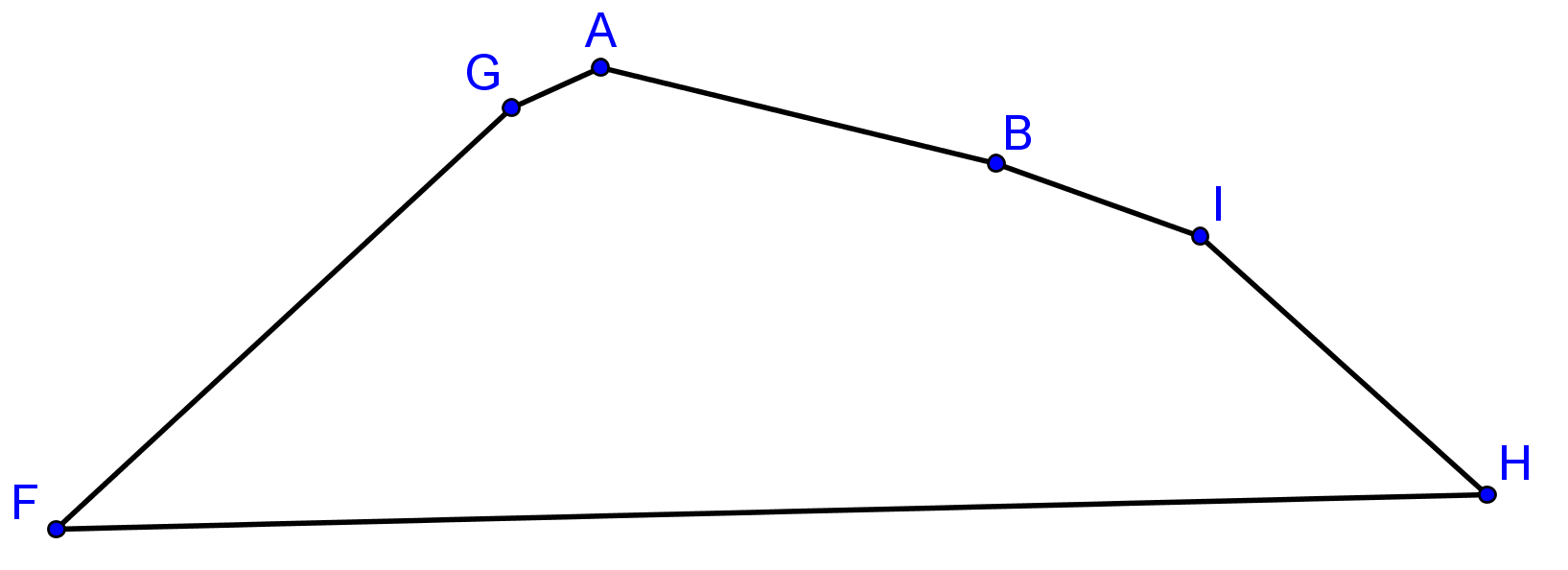}}
\caption[An Example When Local Extremality is Different from Radial Extremality]{}
\end{figure}

It is a routine verification that vertices $A$, $I$ and $F$ are radially maximal-extremal, and it turns out that only $A$ and $I$ are locally maximal-extremal. This shows that these two notions are entirely different once we lose the assumption of coherent.

\subsection{The Four-Vertex Theorems}
Now that we have our three notions of extremality, we turn our attention to three beautiful theorems. It turns out that for special polygonal curves, we have at least four of each type of extremal vertices. We now state these theorems.

\begin{theorem}[The Global Four-Vertex Theorem]
\label{thm:globalfourvertex}
Every generic convex polygon with four or more vertices has at least four globally extremal vertices.
\end{theorem}

\begin{theorem}[The Local Four-Vertex Theorem]
\label{thm:localfourvertex}
Every generic convex polygon with four or more vertices has at least four locally extremal vertices.
\end{theorem}

\begin{theorem}[The Radial Four-Vertex Theorem]
\label{thm:radialfourvertex}
Every generic coherent convex polygon with four or more vertices has at least four radially extremal vertices.
\end{theorem}

Now we turn our attention to proving Theorem \ref{thm:globalfourvertex}, Theorem \ref{thm:localfourvertex}, and Theorem \ref{thm:radialfourvertex}. It turns out that these results follow immediately from a result which arose in the paper of R. C. Bose in 1932 \cite{bos}. We now state this result.

\begin{theorem}
\label{thm:bose}
Let $P$ be a generic convex polygon with $n$ vertices. Then:
$$s_+-t_+=s_- - t_- =2$$ and $$s_+ + t_+ + u_+=s_- + t_- +u_-=n-2.$$
\end{theorem}

Now that we have this result, we can combine it with our results in Section 2.2 to prove our three theorems.

\begin{proof}[Proof of Theorem 2.1]
This follows from the fact that:
$$s_+\geq s_+ - t_+=2$$ and $$s_-\geq s_- - t_-=2.$$
\end{proof}

\begin{proof}[Proof of Theorem 2.2]
From Theorem \ref{thm:globalfourvertex} it follows that we have at least four globally extremal vertices. By Proposition \ref{prop:globaltolocal}, it then follows that we have at least four locally extremal vertices.
\end{proof}

\begin{proof}[Proof of Theorem 2.3]
By Theorem \ref{thm:localfourvertex}, we have at least four locally extremal vertices. With the assumption of coherent and application of Proposition \ref{prop:localtoradial}, we have at least four radially extremal vertices.
\end{proof}

\subsubsection{Alternate Proofs of the Four-Vertex Theorems}

Essentially, the above results followed from the much stronger Theorem \ref{thm:bose}. To give some geometric flavor, we now provide alternate proofs of our three Four-Vertex Theorems. To begin, we first will prove these facts for the most simple of polygons, a quadrilateral.

\begin{prop}
\label{prop:quad}
Let $P$ be a generic convex quadrilateral. Then $P$ has four globally extremal and locally extremal vertices. If $P$ is coherent, then it has four radially extremal vertices.
\end{prop}
\begin{proof}
For globally extremal vertices, we apply a Delaunay triangulation to $P$, which immediately yields two globally maximal-extremal vertices and two minimal-extremal vertices. Proposition \ref{prop:globaltolocal} and \ref{prop:localtoradial} immediately yield the remaining part of our statement.
\end{proof}

While the following proposition is quite obvious, it will be a vital proposition that will be used frequently in Section 4.

\begin{prop}
\label{prop:circleprop}
Let $A$, $B$, $C$ and $X'$ be four points in the plane in a generic arrangement, $C_B$ be the corresponding circle passing through $A$, $B$ and $C$, and let $C_A$ be the circle passing through the points $X'$, $A$ and $B$. We denote by $\widetilde{C}_A$ and $\widetilde{C}_B$ the open discs bounded by $C_A$ and $C_B$, respectively. Denote by $H^{+}_{AB}$ the half-plane formed by the infinite line $AB$ containing the point $C$ and by $H^{-}_{AB}$ the half-plane formed by the infinite line $AB$ not containing the point $C$.

If $X'$ lies in $\widetilde{C}_B\bigcap H^{+}_{AB}$, then $C$ lies in $H^{+}_{AB}\setminus \widetilde{C}_A$. If $X'$ lies in $H^{+}_{AB}\setminus \widetilde{C}_B$, then $C$ lies in $\widetilde{C}_A$.

Analogously, if $X'$ lies in $\widetilde{C}_B\bigcap H^{-}_{AB}$, then $C$ lies in $\widetilde{C}_A$. If $X'$ lies in $H^{-}_{AB}\setminus \widetilde{C}_B$, then $C$ lies in $H^{+}_{AB}\setminus \widetilde{C}_A$.
\end{prop}
\begin{proof}
It is sufficient just to check our situation around the origin. Let $a,r,s>0$. Set $A=(-a,0)$, $B=(a,0)$, $O=(0,r)$ and $O'=(0,s)$. So, $O$ is the center of the circle $C_B$ and $O'$ is the center of the circle $C_A$ in our statement. So we have the following equations for circles $C_A$ and $C_B$:
$$
\begin{array}{cccccc}
  C_{A}: & x^{2} & + & (y-s)^{2} & = & a^2+s^2 \\
  C_{B}: & x^{2} & + & (y-r)^{2} & = & a^2+r^2
\end{array}$$

For the first part of the statement, set $s<r$. Now observe the following system of equations:
$$
\begin{array}{ccccc}
  x^{2} & + & (y-s)^{2} & < & a^2+s^2 \\
  x^{2} & + & (y-r)^{2} & = & a^2+r^2
\end{array}$$

This system only has solutions when $y>0$. If we pick $C$ lying on $C_B$ not satisfying this system of equations, it then follows that $C$ lies outside the circle $C_A$. The following figure illustrates our situation:

\begin{figure}[H]
\centerline{\includegraphics{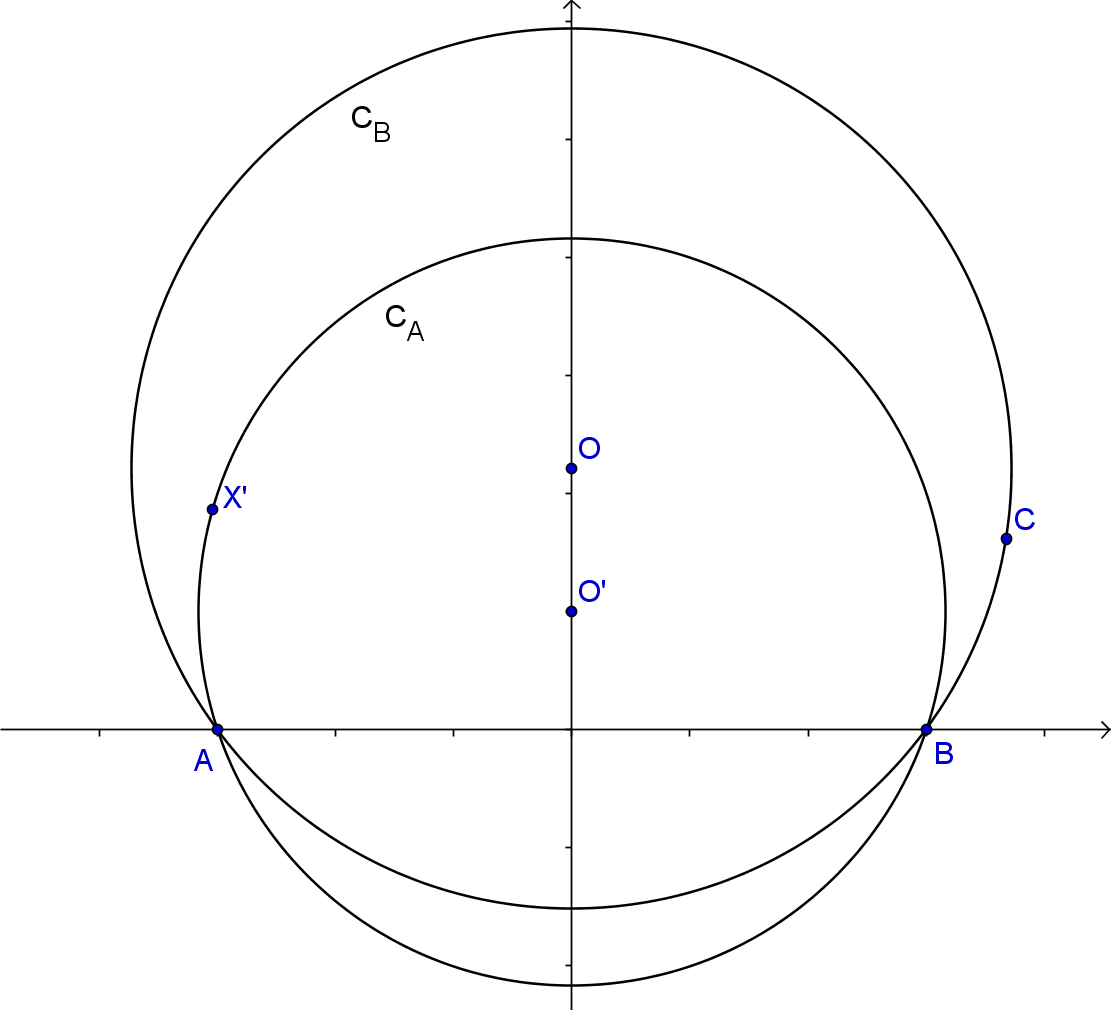}}
\caption[Proving Proposition \ref{prop:circleprop} (1)]{}
\end{figure}

For the proof of the second part of our statement, set $s>r$. Now observe the following system of equations:
$$
\begin{array}{ccccc}
  x^{2} & + & (y-s)^{2} & = & a^2+s^2 \\
  x^{2} & + & (y-r)^{2} & < & a^2+r^2
\end{array}$$
It is clear that this system only has solutions when $y<0$. So, if we pick $X'$ lying on $C_A$ not satisfying the above system, then $X'$ lies outside of the circle $C_B$. It follows that $C$ must lie inside of the circle $C_A$. The following figure illustrates our situation:

\begin{figure}[H]
\centerline{\includegraphics[scale=0.8]{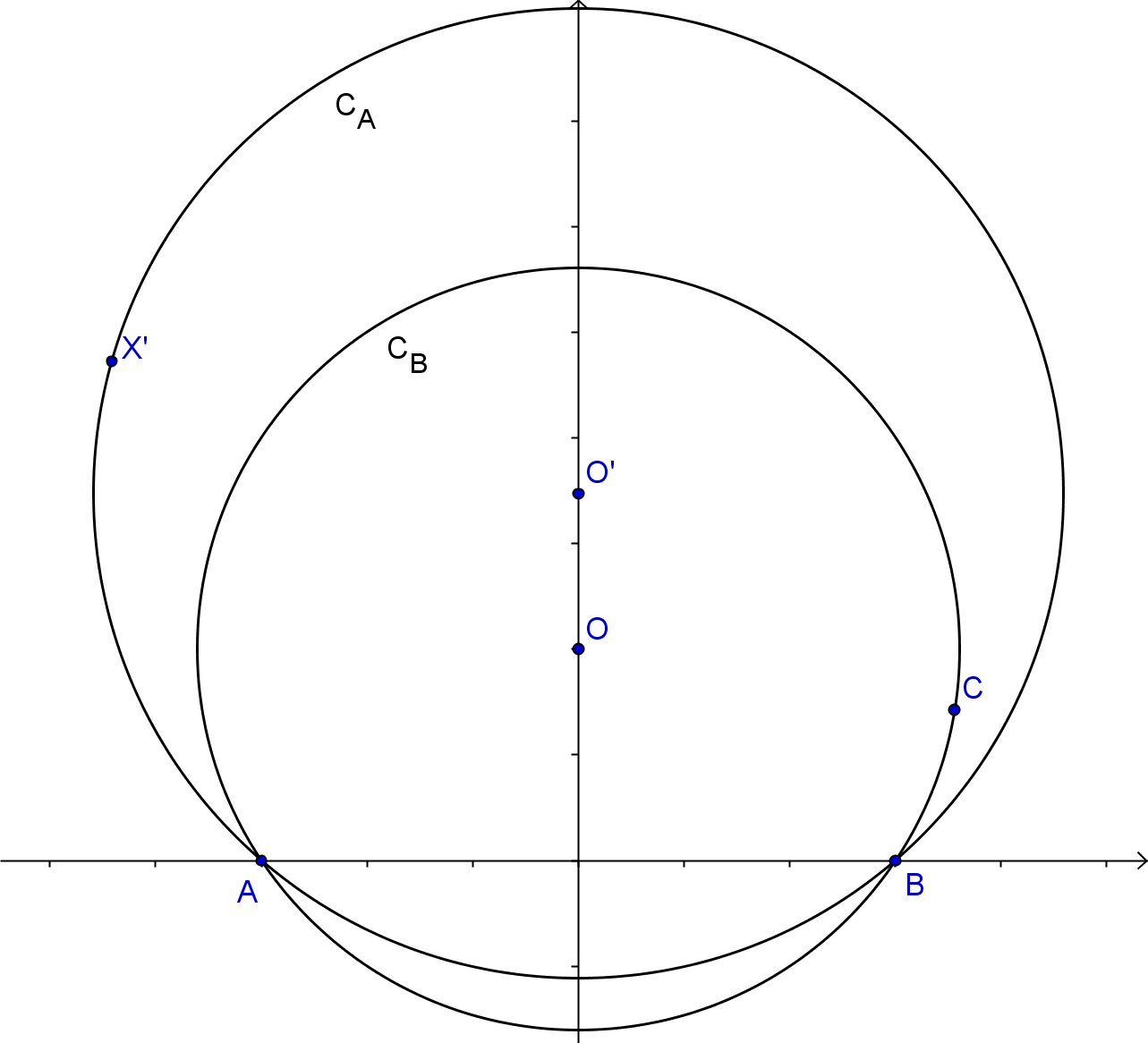}}
\caption[Proving Proposition \ref{prop:circleprop} (2)]{}
\end{figure}

For the analogous statement, the proof follows the exact same scheme as above.
\end{proof}

\begin{lemma}
\label{lemma:globalinductionmax}
Let $P$ be a convex generic polygon with at least five vertices, and let $V_i$ be a globally maximal-extremal vertex. Let $P'$ be the polygon obtained by removing $V_i$ and connecting the vertices $V_{i-1}$ and $V_{i+1}$ by an edge. Then either $s_{-}(P)=s_{-}(P')$ or $s_{-}(P)=s_{-}(P')+1$.
\end{lemma}

\begin{proof}
We essentially must show that, if we have a globally maximal-extremal vertex at any point other than $V_{i-1}$ and $V_{i+1}$ in $P'$, then it stays globally maximal-extremal in $P$. Let $C_i$ be the circle passing through vertices $V_{i-1}$, $V_i$ and $V_{i+1}$. If we pick a globally maximal-extremal vertex $X$ of $P'$, then we have two possibilities: either it is a neighbor of $V_{i-1}$ or $V_{i+1}$, or not.
\vskip2mm
\noindent\textit{Case 1:} $X$ is a neighbor of $V_{i-1}$ or $V_{i+1}$.
\vskip1mm
Without loss of generality, assume that $X$ is the neighbor of $V_{i-1}$. Let $Y$ be the second neighbor of $X$ in $P'$, and denote by $C_X$ the circle passing through the vertices $Y$, $X$ and $V_{i-1}$. Now, since $V_i$ is globally maximal-extremal in $P$, it follows that the circle $C_i$ contains no vertices of $P$, as well as no vertices of $P'$. So, it follows that $X$ and $Y$ lie outside $C_i$.

We must now show that $V_i$ lies outside of $C_X$. Since $X$ is globally maximal-extremal, it follows that $V_{i+1}$ lies outside of the circle $C_X$. So now, denote the point of intersection of the circle $C_X$ with the circle $C_i$ by $Y'$. The following figure illustrates our situation:

\begin{figure}[H]
\centerline{\includegraphics{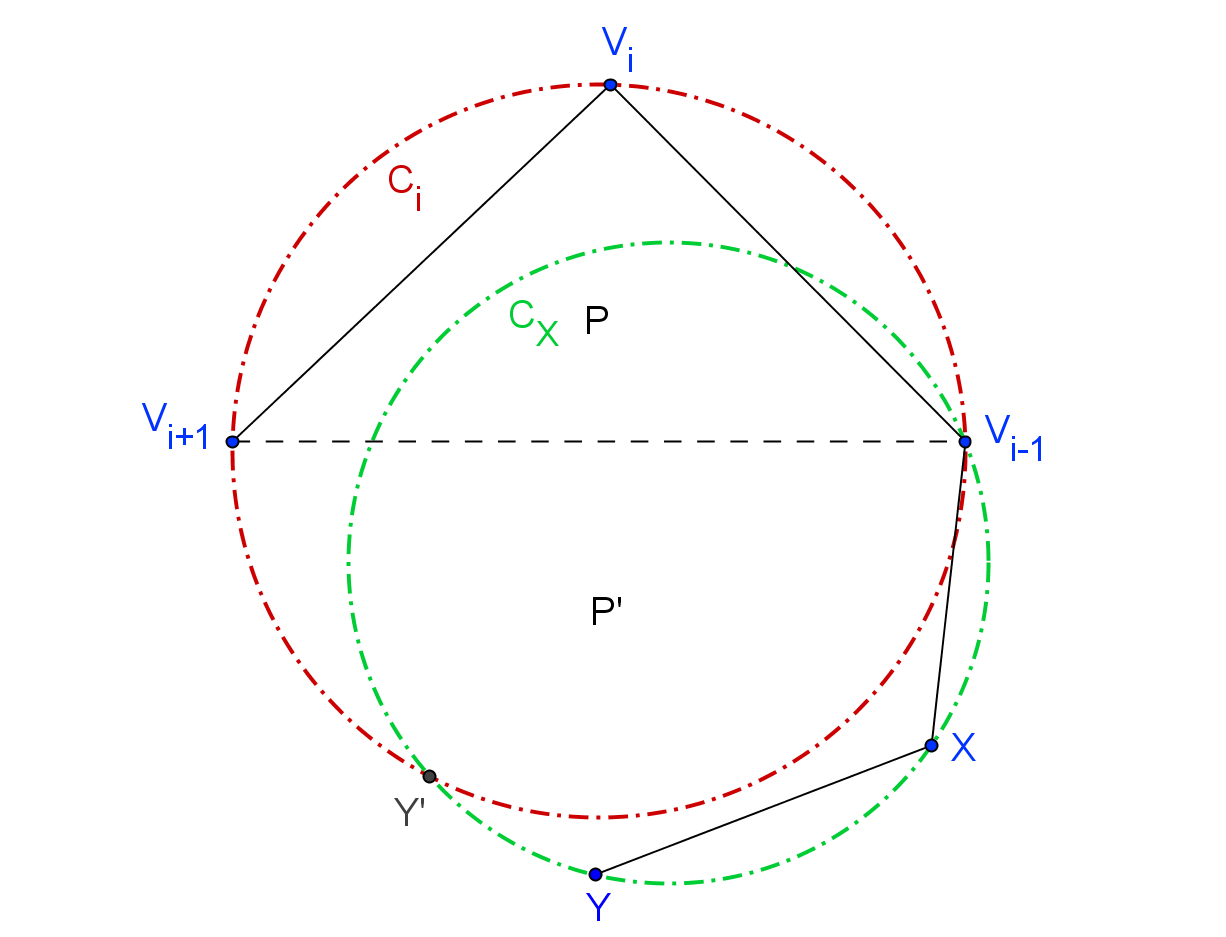}}
\caption[Proving Lemma \ref{lemma:globalinductionmax} (1)]{}
\end{figure}

To finish the proof, we simply apply Proposition \ref{prop:circleprop} to the points $X$, $V_{i-1}$, $Y'$ and $V_i$, observing that $X$ lies in a different half-plane than $V_i$ with respect to the infinite line $Y'V_{i-1}$. This yields that $V_i$ lies outside of the circle $C_X$, proving Case 1.
\vskip2mm
\noindent\textit{Case 2:} $X$ is not a neighbor of $V_{i-1}$ and is not a neighbor of $V_{i+1}$.
\vskip1mm
Denote by $X$ a globally maximal-extremal vertex of $P'$ with neighbors $Y$ and $Z$, and denote the circle passing through $Y$, $X$ and $Z$ by $C_X$. Now, since $V_i$ is globally maximal-extremal in $P$, it follows that the circle $C_i$ contains no vertices of $P$, as well as no vertices of $P'$. So, it follows that $X$, $Y$ and $Z$ lie outside $C_i$.

Now, we have that $X$ is globally maximal in $P'$, so we know that vertices $V_{i-1}$ and $V_{i+1}$ must lie outside the circle $C_X$. Our goal is now to show that $C_X$ does not contain the vertex $V_i$. If $C_X$ does not intersect $C_i$, then we are done. So, we assume that $C_X$ intersects $C_i$ and denote the intersection points of the circle $C_X$ with the circle $C_i$ by $Y'$ and $Z'$. The following figure illustrates our situation:

\begin{figure}[H]
\centerline{\includegraphics{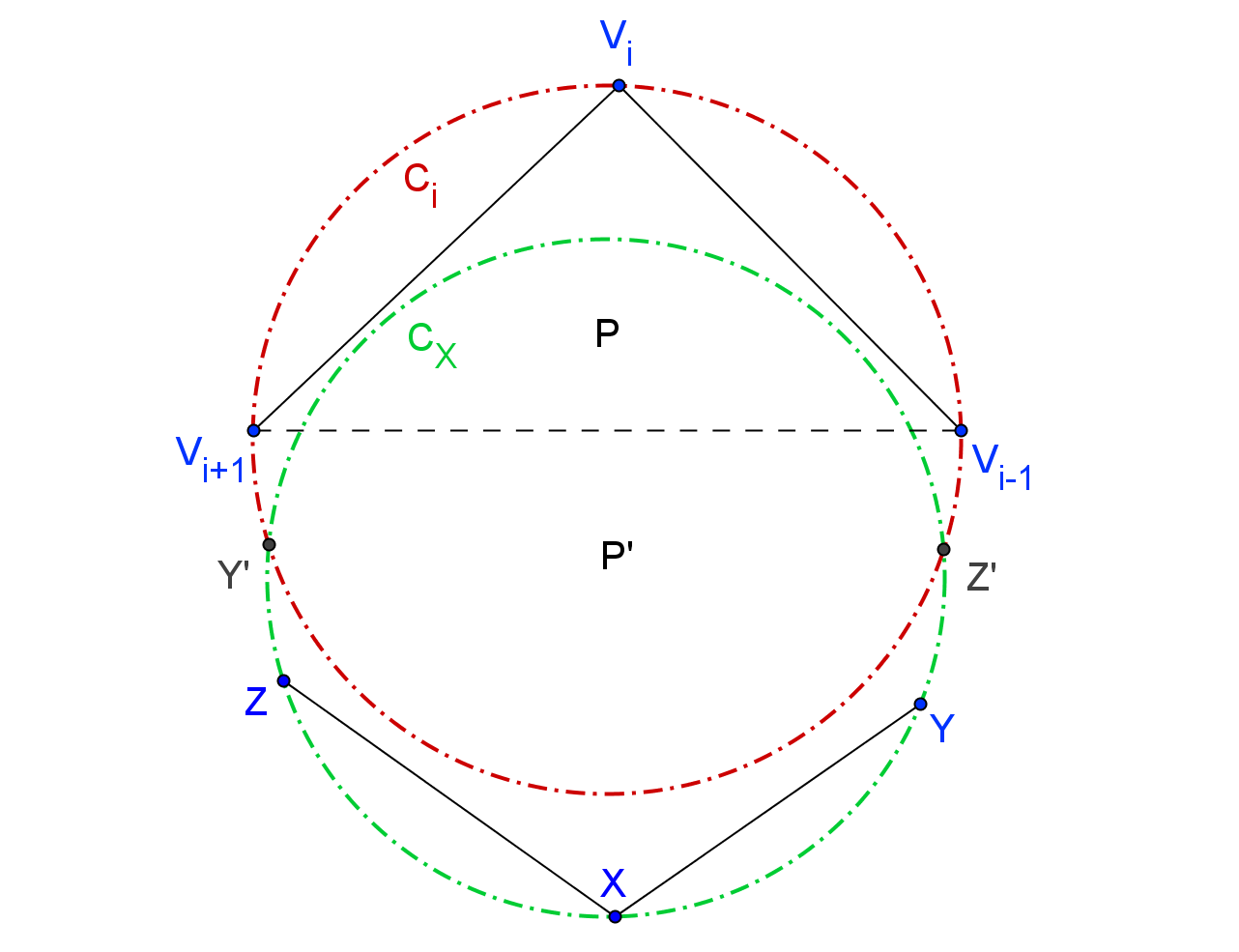}}
\caption[Proving Lemma \ref{lemma:globalinductionmax} (2)]{}
\end{figure}

We know that $X$ lies on a different half-plane than $V_i$ with respect to the infinite line $Y'Z'$, as well as outside the circle $C_i$. An application of Proposition \ref{prop:circleprop} to points $X$, $Y'$, $Z$' and $V_i$ yields that $V_i$ lies outside the circle $C_X$.
\end{proof}

The previous section states that we do not necessarily have the same number of globally maximal-extremal vertices as globally minimal-extremal vertices. So, for an alternate proof of the Global Four-Vertex Theorem, we must prove a similar lemma as Lemma \ref{lemma:globalinductionmax} for globally minimal-extremal vertices.

\begin{lemma}
\label{lemma:globalinductionmin}
Let $P$ be a convex generic polygon with at least five vertices, and let $V_i$ be a globally minimal-extremal vertex. Let $P'$ be the polygon obtained by removing $V_i$ and connecting the vertices $V_{i-1}$ and $V_{i+1}$ by an edge. Then either $s_{+}(P)=s_{+}(P')$ or $s_{+}(P)=s_{+}(P')+1$.
\end{lemma}
\begin{proof}
We essentially must show that, if we have a globally maximal-extremal vertex at any point other than $V_{i-1}$ and $V_{i+1}$ in $P'$, then it stays globally maximal-extremal in $P$. Let $C_i$ be the circle passing through vertices $V_{i-1}$, $V_i$ and $V_{i+1}$. If we pick a globally maximal-extremal vertex $X$ of $P'$, then we have two possibilities: either it is a neighbor of $V_{i-1}$ or $V_{i+1}$, or not.
\vskip2mm
\noindent\textit{Case 1:} $X$ is a neighbor of $V_{i-1}$ or $V_{i+1}$.
\vskip1mm
Without loss of generality, assume that $X$ is the neighbor of $V_{i-1}$. Let $Y$ be the second neighbor of $X$ in $P'$, and denote by $C_X$ the circle passing through the vertices $Y$, $X$ and $V_{i-1}$. Now, since $V_i$ is globally minimal-extremal in $P$, it follows that the circle $C_i$ contains all vertices of $P$, as well as all of $P'$. So, it follows that $X$ and $Y$ lie outside $C_i$.

We must now show that $V_i$ lies inside of $C_X$. Since $X$ is globally minimal-extremal, it follows that $V_{i+1}$ lies inside of the circle $C_X$. So now, denote the point of intersection of the circle $C_X$ with the circle $C_i$ by $Y'$. The following figure illustrates our situation:

\begin{figure}[H]
\centerline{\includegraphics{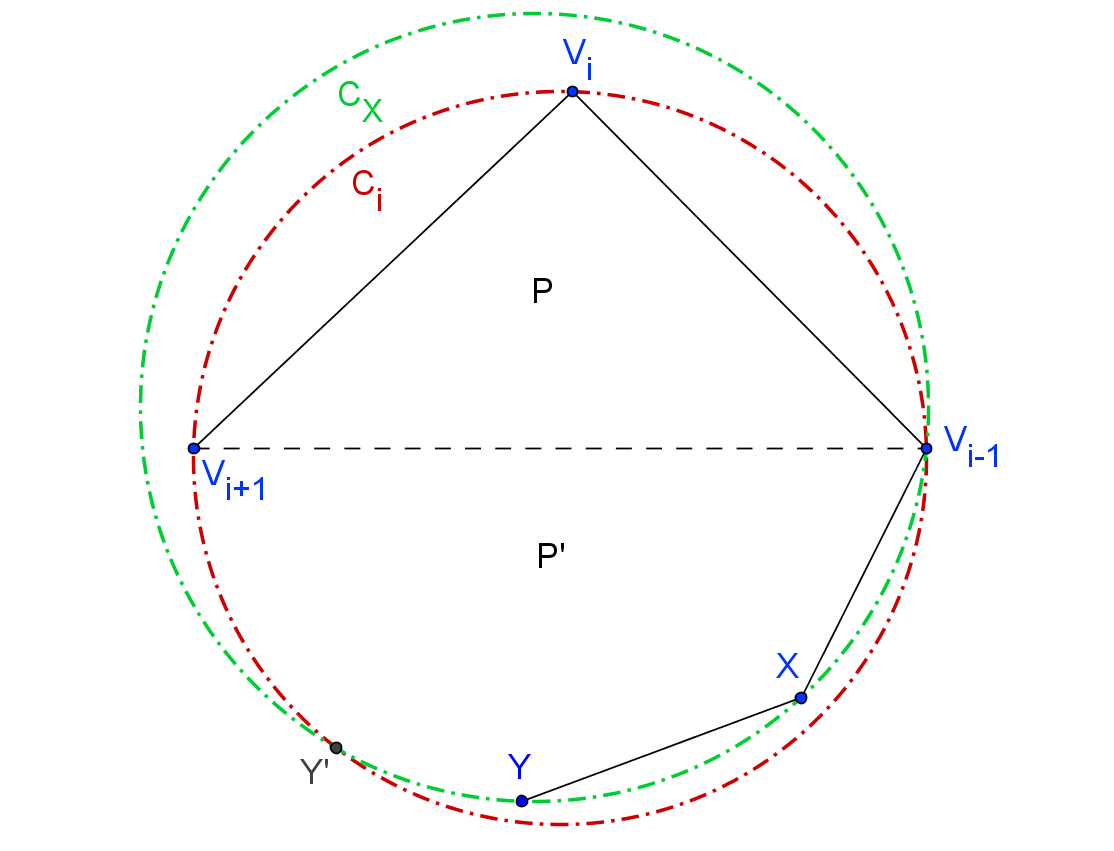}}
\caption[Proving Lemma \ref{lemma:globalinductionmin} (1)]{}
\end{figure}

To finish the proof, we simply apply Proposition \ref{prop:circleprop} to the points $X$, $V_{i-1}$, $Y'$ and $V_i$, observing that $X$ lies in a different half-plane than $V_i$ with respect to the line $Y'V_{i-1}$. This yields that $V_i$ lies inside of the circle $C_X$, proving Case 1.
\vskip2mm
\noindent\textit{Case 2:} $X$ is not a neighbor of $V_{i-1}$ and is not a neighbor of $V_{i+1}$.
\vskip1mm

Denote by $X$ a globally minimal-extremal vertex of $P'$ with neighbors $Y$ and $Z$, and denote the circle passing through $Y$, $X$ and $Z$ by $C_X$. Now, since $V_i$ is globally minimal-extremal in $P$, it follows that the circle $C_i$ contains all vertices of $P$, as well as $P'$. So, it follows that $X$, $Y$ and $Z$ are lie inside $C_i$.

Now, we have that $X$ is globally minimal in $P'$, so we know that vertices $V_{i-1}$ and $V_{i+1}$ must lie inside the circle $C_X$. Our goal is now to show that $C_X$ contains the vertex $V_i$. If $C_X$ does not intersect $C_i$, then we are done. So, we assume that $C_X$ intersects the circle $C_i$ and denote the intersection points of the circle $C_X$ with the circle $C_i$ by $Y'$ and $Z'$. The following figure illustrates our situation:

\begin{figure}[H]
\centerline{\includegraphics[scale=0.9]{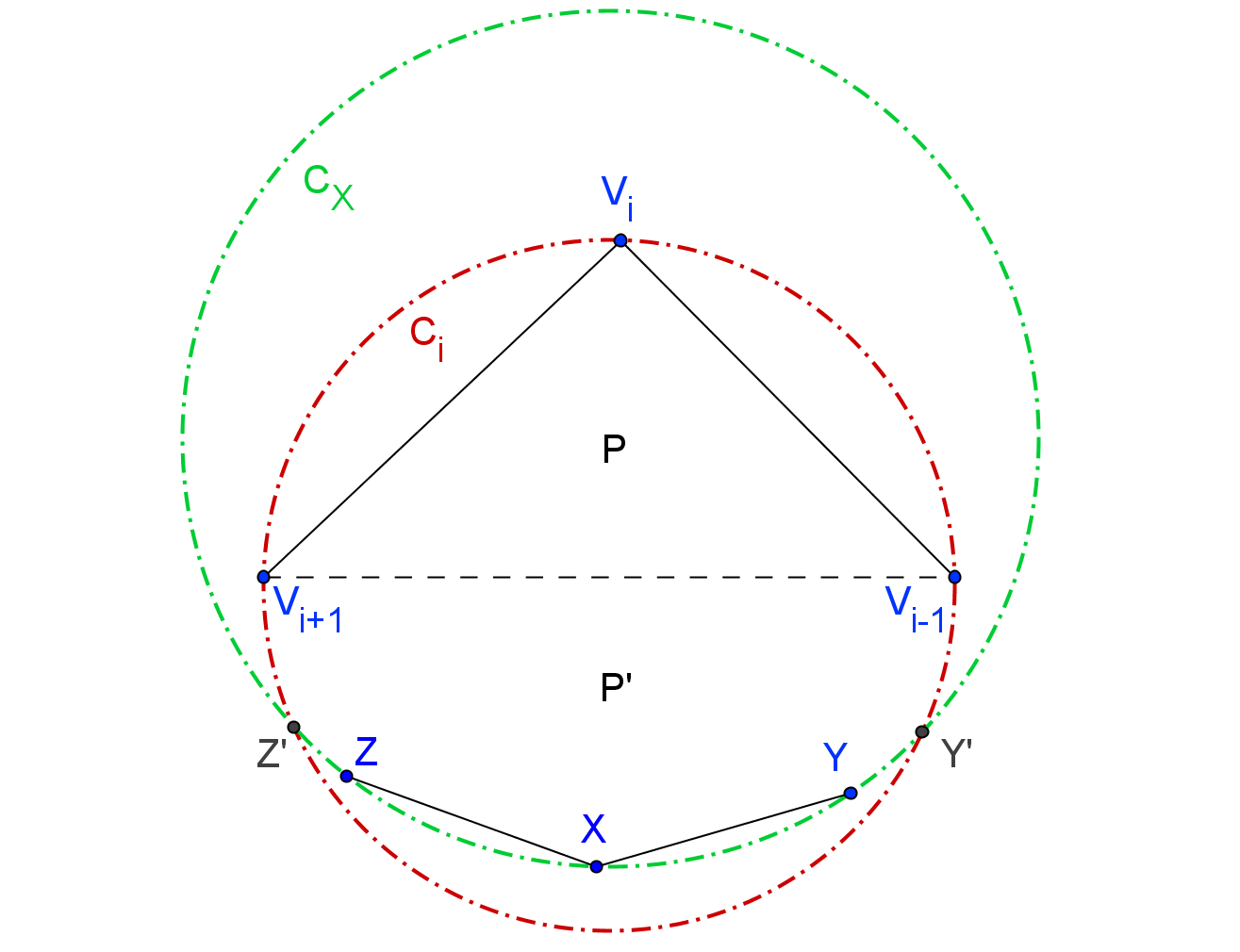}}
\caption[Proving Lemma \ref{lemma:globalinductionmin} (2)]{}
\end{figure}

We know that $X$ lies on a different half-plane than $V_i$ with respect to the line $Y'Z'$, as well as inside the circle $C_i$. An application of Proposition \ref{prop:circleprop} to points $X$, $Y'$, $Z$' and $V_i$ yields that $V_i$ lies inside the circle $C_X$.
\end{proof}

These two lemmas now give us an alternate proof of Theorem \ref{thm:globalfourvertex}.

\begin{proof}[Alternate Proof of Theorem \ref{thm:globalfourvertex}]
We will prove the theorem first for globally maximal vertices by induction. Proposition \ref{prop:quad} takes care of the case where $n=4$. Let $P$ be a closed generic convex polygon with vertices $V_{1},V_{2},...,V_{n}.$ Proposition \ref{prop:existenceprop} guarantees us at least one globally maximal-extremal vertex, without loss of generality, say $V_{i}.$ We now remove this vertex to obtain $P'$, a polygon with $n-1$ vertices. By induction, it follows that $P'$ has at least two globally extremal vertices. By Lemma \ref{lemma:globalinductionmax} we have that either $s_{-}(P)=s_{-}(P')$ or $s_{-}(P)=s_{-}(P')+1$, so it follows that $P$ has at least two globally maximal vertices.

Mimicking the exact same technique as above and using Lemma \ref{lemma:globalinductionmin}, we obtain that we have at least two globally minimal vertices. So, our assertion is proved.
\end{proof}

With the help of one lemma, we will provide an alternate proof of Theorem \ref{thm:localfourvertex}.

\begin{lemma}
\label{lemma:inductionmax}
Let $P$ be a convex generic polygon with at least five vertices, and let $V_i$ be a locally maximal-extremal vertex. Denote by $V_{i-1}$ and $V_{i+1}$ the neighboring vertices of $V_i$, and by $X$ a neighboring vertex of either $V_{i-1}$ or $V_{i+1}$. Consider the polygon $P'$ formed by removing the vertex $V_i$ and joining the vertices $V_{i-1}$ and $V_{i+1}$ by an edge. If $X$ is locally maximal-extremal in $P'$, then it is locally maximal-extremal in $P$.
\end{lemma}
\begin{proof}
Without loss of generality, we assume that our vertex $X$ is the neighbor of vertex $V_{i-1}$, and $Y$ is the neighbor of $X$. Denote the circle passing through $V_{i-1}$, $V_{i}$ and $V{i+1}$ by $C_i$, the circle passing through $X$, $V_{i-1}$ and $V_{i}$ by $C_{i-1}$, and the circle passing through $V_{i-1}$, $X$ and $Y$ by $C_X$. The following figure illustrates our situation:

\begin{figure}[H]
\centerline{\includegraphics[scale=0.7]{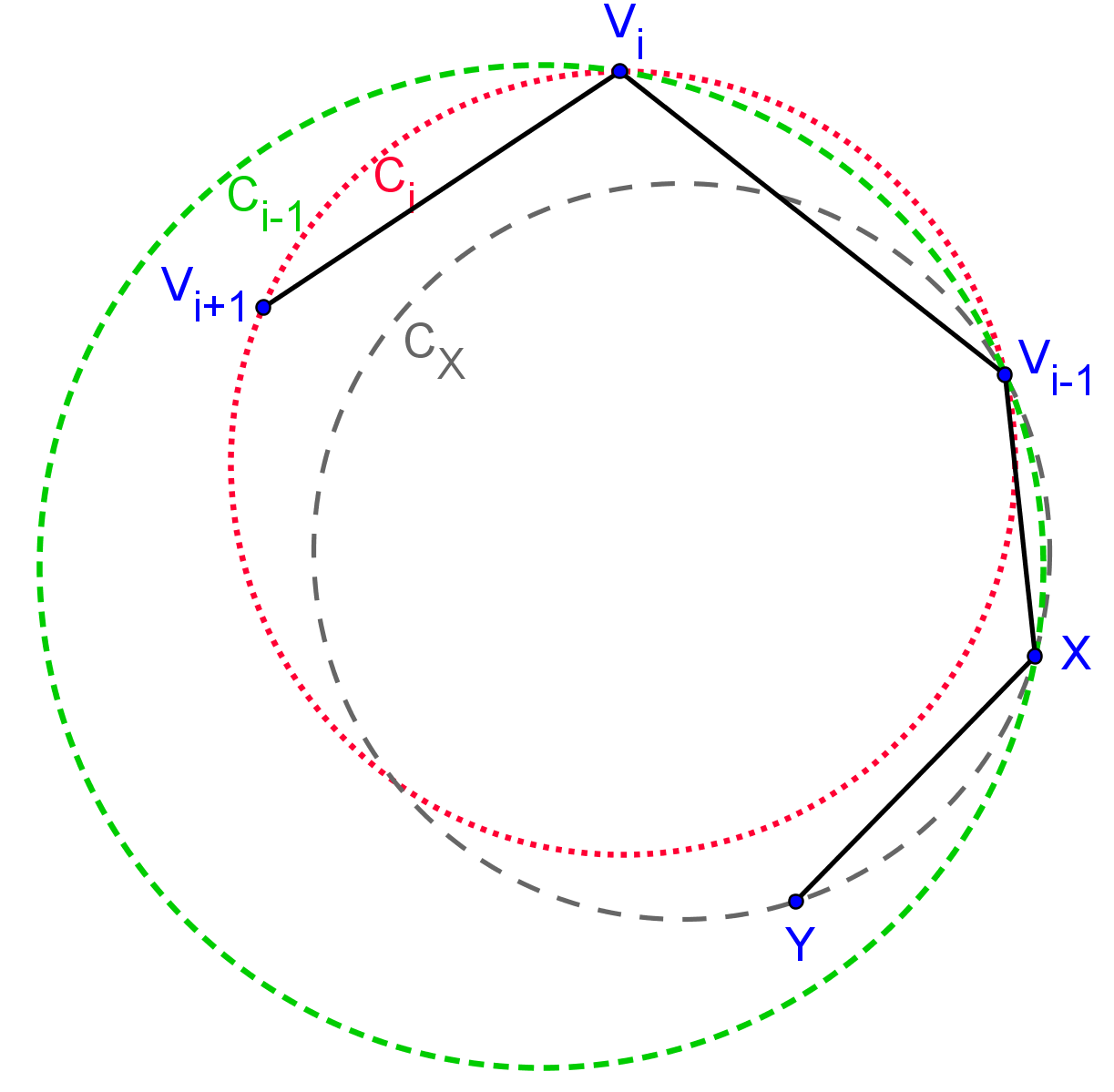}}
\caption[Proving Lemma \ref{lemma:inductionmax}]{}
\end{figure}

Since $V_i$ is maximal-extremal in $P$, it follows that $X$ must lie outside of the circle $C_{i}$. By Proposition \ref{prop:circleprop}, it follows that the circle $C_{i-1}$ contains $V_{i+1}$. Since $X$ is maximal-extremal in $P'$, it follows that $V_{i+1}$ lies outside of the circle $C_X$. It follows by Proposition \ref{prop:circleprop} that $Y$ must lie inside the circle $C_{i-1}$, because otherwise $C_X$ would contain $V_{i+1}$. Since $Y$ lies inside the circle $C_{i-1}$, Proposition \ref{prop:circleprop} tells us that $V_i$ must lie outside the circle $C_X$. Therefore $X$ is maximal-extremal in $P$.
\end{proof}

\begin{proof}[Alternate Proof of Theorem \ref{thm:localfourvertex}]
We will prove the theorem by induction. Proposition \ref{prop:quad} takes care of the case where $n=4$. Let $P$ be a closed generic convex polygon with vertices $V_{1},V_{2},...,V_{n}.$ Proposition \ref{prop:existenceprop} guarantees us at least one locally maximal-extremal vertex, without loss of generality, say $V_{i}.$ If we remove this vertex, we obtain a convex polygon $P'$ with $n-1$ vertices. By Lemma \ref{lemma:inductionmax}, it follows that by adding $V_{i}$ to $P'$, we will lose at most one locally maximal-extremal vertex of $P'.$ Applying our inductive assumption to $P'$ and reattaching the vertex $V_i$ to obtain $P$, we see that $P$ has at least two locally maximal-extremal vertices. Proposition \ref{prop:maxmin} yields two locally minimal-extremal vertices.
\end{proof}

\begin{proof}[Alternate Proof of Theorem \ref{thm:radialfourvertex}]
This proof is thanks to O. R. Musin \cite{mus2}.

Since we have only a finite number of radii, we already have two extremal vertices, since one of them must have the largest size and one must have the smallest size.
We will prove the theorem indirectly:
\vskip2mm
Without loss of generality, we can assume that these vertices are $V_{1}$ and $V_{k}$.
Let $V_{1}$ have minimal radius and $V_{k}$ have maximal radius. (By radius, we mean the radius $R_{i}$ of the of the circle through $V_{i-1}V_{i}V_{i+1}$.)
\vskip2mm
Then we have:
\begin{enumerate}
\item $R_{1}<R_{2}<...<R_{k}$
\item $R_{1}<R_{n}<R_{n-1}<...<R_{k+1}<R_{k}$
\end{enumerate}
\vskip2mm
We now construct the perpendicular bisector for each line segment of the polygonal curve.
Two consecutive bisectors for line segments $l_{i-1}$ and $l_{i}$ will intersect at the point $O_{i}$, which is the center of the circle through $V_{i-1}V_{i}V_{i+1}$. (See Figure 2.15)

\begin{figure}[H]
\centerline{\includegraphics[width=120mm,height=95mm]{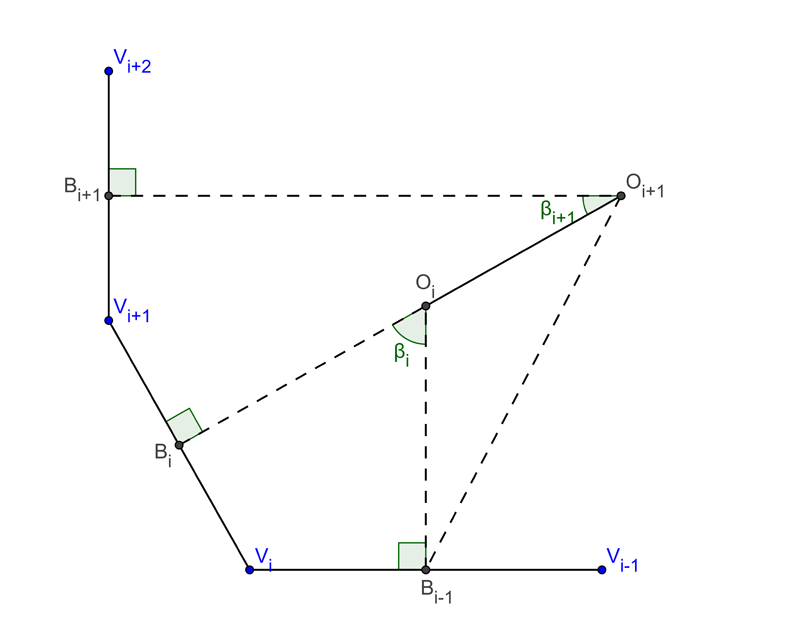}}
\caption[Proving Theorem \ref{thm:radialfourvertex}]{}
\end{figure}

By the definition of $R_{i}$, we have $V_{i}O_{i}=R_{i}.$
Let $B_{i}$ be the midpoint of the line segment $V_{i}V_{i+1}$.
Then $\beta_{i}=\angle B_{i-1}O_{i}B_{i}=\pi -\angle V_{i-1}V_{i}V_{i+1}$.
\vskip2mm
We will prove that if $R_{i}<R_{i+1}$, then $\angle B_{i-1}O_{i+1}B_{i+1}<\beta_{i}+\beta_{i+1}.$
Let us notice that $B_{i}O_{i}<B_{i}O_{i+1}$. This is because the right triangles $B_{i}O_{i}V_{i}$ and $B_{i}O_{i+1}V_{i+1}$ have $B_{i}V_{i}$ and $B_{i}V_{i+1}$ of equal length.
Thus the angle $\beta_{i}$ is supplementary to $\angle B_{i-1}O_{i}O_{i+1},$ and necessarily $\angle B_{i-1}O_{i+1}B_{i}<\beta_{i}.$ From this, we have that $\angle B_{i-1}O_{i+1}B_{i+1}<\beta_{i}+\beta_{i+1}.$, which is our desired conclusion.
\vskip2mm
Applying this last result by step-by-step considering the angles $B_{1}O_{3}B_{3},B_{1}O_{4}B_{4},...$, we obtain that:
$$\angle B_{1}O_{k}B_{k}<\beta_{2}+\beta_{3}+...+\beta_{k}.$$
Similarly, using (2), we obtain that
$$2\pi-\angle B_{1}O_{k}B_{k}<\beta_{1}+\beta_{n}+\beta_{n-1}+...+\beta_{k-1}.$$
Combining the last two inequalities, we obtain
$$2\pi<\beta_{1}+\beta_{2}+...+\beta_{n}=\pi-\angle V_{1}+\pi-\angle V_{2}+...+\pi-\angle V_{n}=$$
$$=n\pi-(\angle V_{1}+\angle V_{2}+...+\angle V_{n})=n\pi-\pi(n-2)=2\pi$$
Here we used the fact that the sum of the angles of an n-gon equals to $\pi(n-2)$.
Thus, $2\pi<2\pi$, a contradiction.
\end{proof}

\subsection{Some Counterexamples}
\subsubsection{A Counterexample to the Radial Four-Vertex Theorem}

In this section we investigate the reason why we added the extra assumption that the polygonal curve is coherent. Observe the following figure:

\begin{figure}[H]
\centerline{\includegraphics[width=120mm,height=38mm]{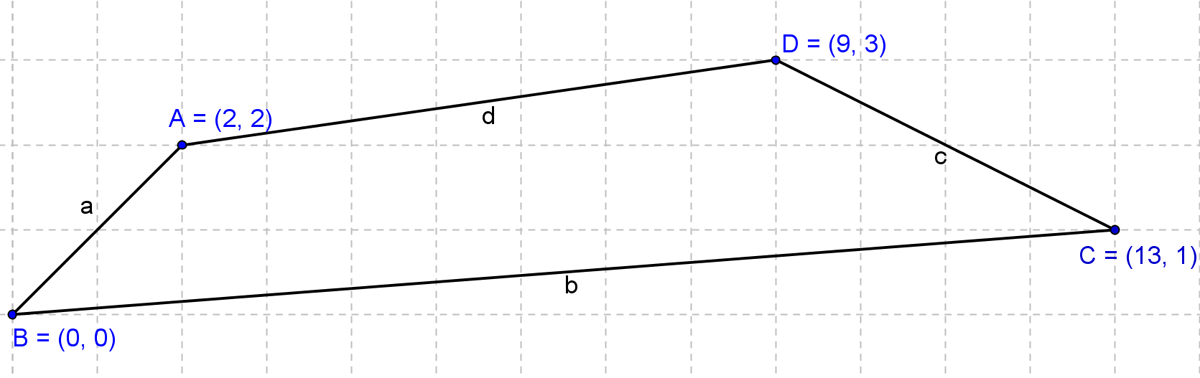}}
\caption[A Counterexample to the Radial Four-Vertex Theorem]{}
\end{figure}

We first observe that our polygon is not coherent. We can easily verify this by computing the equations of the circles $C_{A}, C_{B}, C_{C}, C_{D}$:
$$
\begin{array}{cccccc}
  C_{A}: & (x-6.5)^{2} & + & (y+4.5)^{2} & = & 62.5 \\
  C_{B}: & (x-6.92)^{2} & + & (y+4.92)^{2} & = & 72.01 \\
  C_{C}: & (x-7)^{2} & + & (y+6)^{2} & = & 85 \\
  C_{D}: & (x-6.78)^{2} & + & (y+6.44)^{2} & = & 94.14
\end{array}$$

We observe that $R_{A}<R_{B}<R_{C}<R_{D}.$ So, we conclude that we have exactly two radially extremal vertices, vertex $A$ and vertex $D$. This example shows why this extra assumption was necessary.

\subsubsection{Local Does Not Imply Global}

Here we have a simple counterexample as to why a locally extremal vertex is not necessarily globally extremal. Consider the following figure:

\begin{figure}[H]
\centerline{\includegraphics[scale=1]{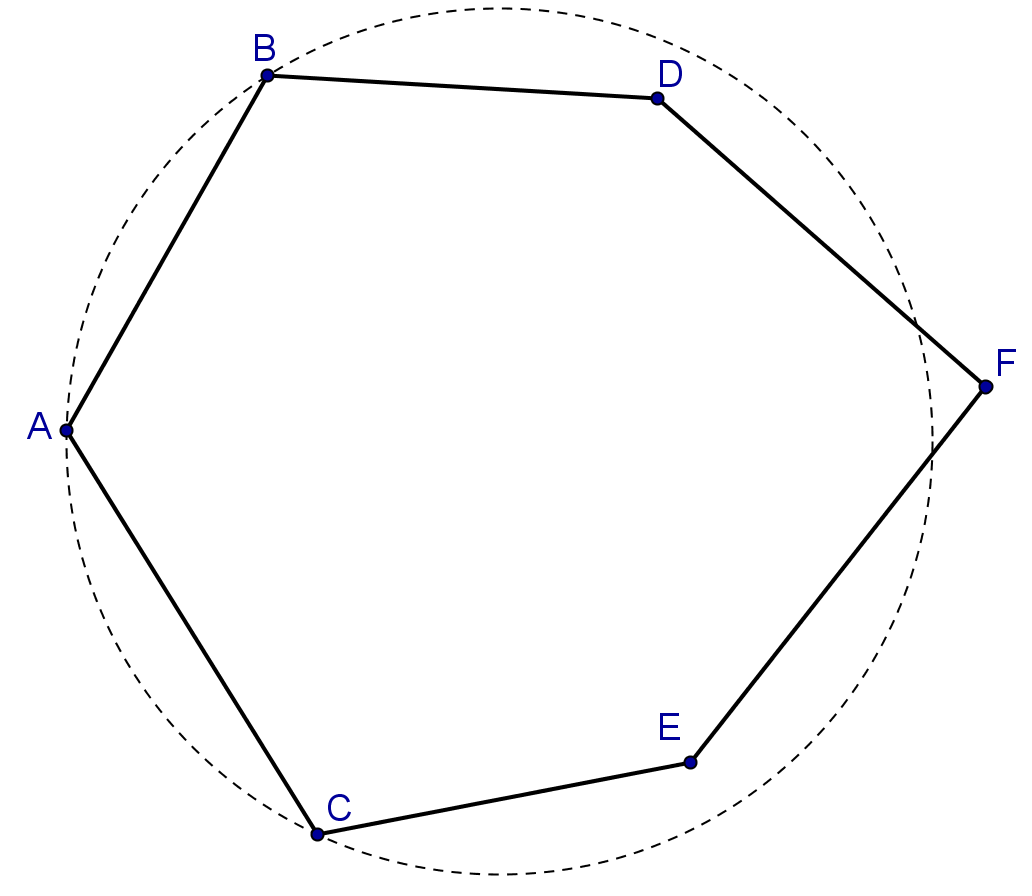}}
\caption[An Example Where Local Extremality Does Not Imply Global Extremality]{}
\end{figure}

The vertex under consideration is vertex $A.$ Since vertices $B$ and $E$ lie inside the circle $C_A,$ it follows that the vertex $A$ is locally extremal. But, the vertex $F$ lies outside of the circle $C_A,$ so we see that vertex $A$ is not globally extremal.

\section{The Local Four-Vertex Theorem and the Evolute}
\renewcommand{\thefigure}{\thesection.\arabic{figure}}
\setcounter{figure}{0}

In this section, we give an overview and fill in all details of O. R. Musin's 2004 paper \cite{mus1}. We will define the notion of the evolute of a curve and see how closely this notion is tied to the Local Four-Vertex Theorem. In this section we will simply restrict ourselves to locally extremal vertices, and will denote the number of positive locally extremal vertices by $N_{+}$ and the number of negative locally extremal vertices by $N_{-}$.
\vskip2mm
\begin{remark}
Observe that if we assume convexity, we have $N_{-}=0$.
\end{remark}
\subsection{The Evolute}

Here we will consider the evolute of a polygonal curve.

\begin{defn}
The figure formed by the centers $O_{1},O_{2},...,O_{n}$ is called the evolute $E(P)$ of $P$. For a smooth curve, the evolute is the curve formed by the centers of the osculating circles at each point of the curve (centers of curvature).
\end{defn}

\begin{remark}
It is important to recall from elementary geometry that the center of a circle circumscribing three points is the intersection of the three perpendicular bisectors of the triangle formed by the three points. This fact will be heavily used for proving the results of this section.
\end{remark}

\vskip2mm

We consider some examples of evolutes to build some geometric intuition. A simple example would be the one of a square, which has a point as its evolute. This is because the four vertices of a square lie on a circle. Figure 3.1 and 3.2 give examples of evolutes of polygons. (Note that the polygon is blue and the evolute is green.)

\begin{figure}[H]
\centerline{\includegraphics[width=80mm,height=63mm]{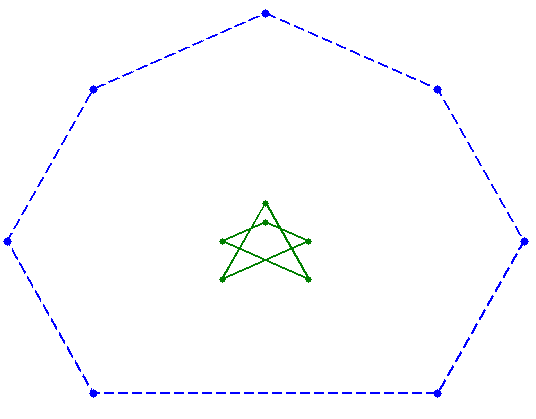}}
\caption[The Evolute of a Pentagon]{}
\end{figure}

Figure 3.1 shows the evolute of the 7-gon with vertices given by the following matrix:

$$\left(
  \begin{array}{ccccccc}
    2 & 3 & 2 & 0 & -2 & -3 & -2 \\
    0 & 2 & 4 & 5 & 4 & 2 & 0 \\
  \end{array}
\right)$$

\begin{figure}[H]
\centerline{\includegraphics[width=80mm,height=61mm]{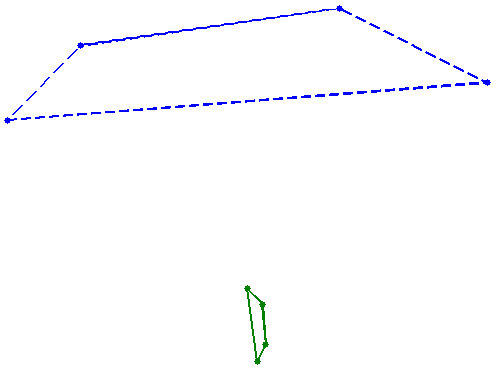}}
\caption[The Evolute of Figure 2.16]{}
\end{figure}

Figure 3.2 shows the evolute of Figure 2.16.

\begin{figure}[H]
\centerline{\includegraphics[width=80mm,height=58mm]{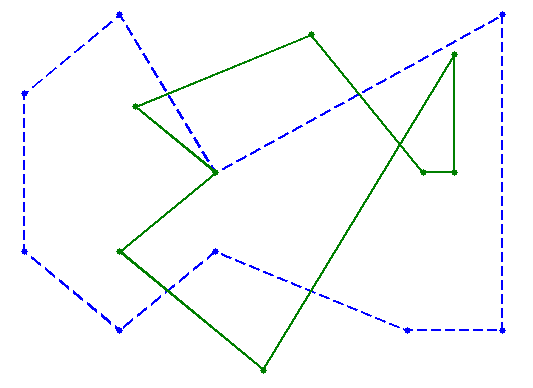}}
\caption[The Evolute of a Nonagon]{}
\end{figure}

Figure 3.3 shows the evolute of a 9-gon with vertices given by the following matrix:

$$\left(
  \begin{array}{ccccccccc}
    0 & 1 & 3 & 4 & 4 & 1 & 0 & -1 & -1 \\
    1 & 2 & 1 & 1 & 5 & 3 & 5 & 4 & 2 \\
  \end{array}
\right)$$

\vskip4mm
The following figures give some examples of evolutes of smooth curves:

\begin{figure}[H]
\centerline{\includegraphics[width=80mm,height=69mm]{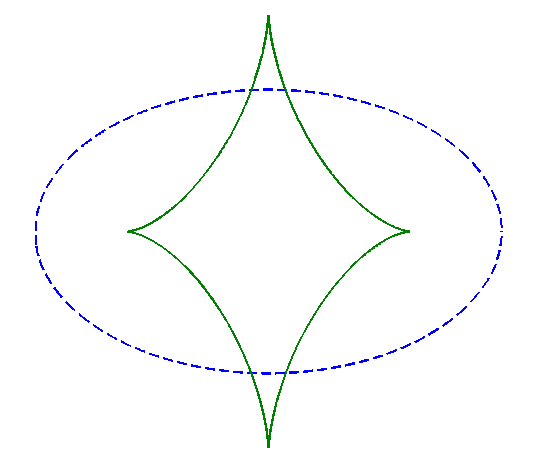}}
\caption[The Evolute of an Ellipse]{}
\end{figure}

Figure 3.4 shows the evolute of the ellipse:

$$\begin{array}{ccc}
  x(t) & = & cos(t) \\
  y(t) & = & 0.63sin(t)
\end{array},$$

where $0\leq t\leq 2\pi.$

\begin{figure}[H]
\centerline{\includegraphics[width=80mm,height=56mm]{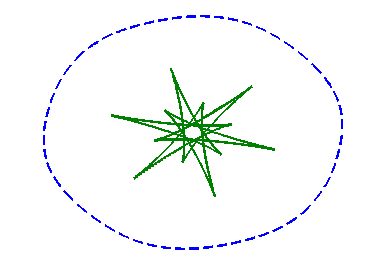}}
\caption[An Evolute With Many Cusps]{}
\end{figure}

Figure 3.5 shows the evolute the curve:
$$\begin{array}{ccc}
  x(t) & = & r(t)cos(t) \\
  y(t) & = & r(t)sin(t)
\end{array},$$

where
$$
\begin{array}{ccccc}
  r(t) & = & 1 & + & \frac{1}{100}sin(6t)
\end{array}
$$

and $0\leq t\leq 2\pi.$
\vskip1mm
This curve has a very interesting geometry since it almost resembles a circle, but has much less extremal vertices.

\begin{figure}[H]
\centerline{\includegraphics[width=80mm,height=68mm]{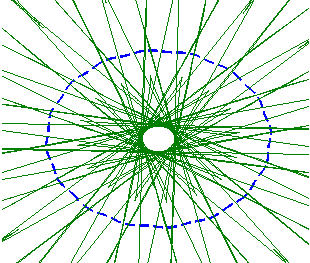}}
\caption[An Evolute With Even More Cusps]{}
\end{figure}

Figure 3.6 shows the evolute the curve:
$$\begin{array}{ccc}
  x(t) & = & r(t)cos(t) \\
  y(t) & = & r(t)sin(t)
\end{array},$$

where
$$
\begin{array}{ccccc}
  r(t) & = & 1 & + & \frac{1}{100}sin(16t)
\end{array}
$$

and $0\leq t\leq 2\pi.$
\vskip1mm
Observe the slight difference in the equations of Figure 3.5 and Figure 3.6.

\vskip2mm
We now investigate important behavior of the evolute near extremal vertices, which is particularly evident in the last three figures. We observe that the evolute appears to have a sharp ``point" near an extremal vertex. This leads us to the next definition.
\vskip2mm
\begin{defn}
A vertex of the evolute is said to be a cusp if
$$\angle V_{i} - \angle O_{i} = \pm\pi.$$
\end{defn}
\vskip2mm
The following lemma illustrates the importance of this notion.
\vskip2mm
\begin{lemma}
\label{lemma:evoluteangles}
A vertex $V_{i}$ is locally extremal if and only if $O_{i}$ is a cusp. Also,
A vertex $V_{i}$ is not locally extremal if and only if $$\angle V_{i} - \angle O_{i} = 0.$$
\end{lemma}
\begin{proof}
The proof of this lemma is a case by case analysis. For the first part, we must run through both cases of extremality of the vertex $V_{i}$ and then run through all the cases of positivity and negativity of the vertices $V_{i-1}$ and $V_{i+1}$, in the mean time considering the positions of the vertices $V_{i-2}$ and $V_{i+2}$ with respect to the circle $C_{i}$. This will give us a total of sixteen cases for the first part. For the second part of the lemma, we must consider all of the cases of non-extremality of the vertex $V_{i}.$ Similarly, this will give us again sixteen cases. So to fully prove this lemma, we must consider a total of 32 cases. We will simply consider a total of six cases, since the the technique of the proof in every case is exactly the same.
\vskip2mm
We now prove two cases for the first part of our lemma, which will have two subcases. Since $V_{i}$ is extremal, either $V_{i+1}\prec V_{i} \succ V_{i-1}$ or $V_{i+1}\succ V_{i} \prec V_{i-1}.$

\vskip2mm
Case 1: $V_{i+1}\prec V_{i} \succ V_{i-1}$
\vskip2mm
For our first subcase, assume that $V_{i-1}$, $V_{i}$, and $V_{i+1}$ are positive and $V_{i-2}$ and $V_{i+2}$ lie outside of the circle $C_{i}.$ Now, label the perpendicular bisector of the line segment formed by the vertices $V_{i-1}$ and $V_{i}$ as $b_{i}$. We now focus our attention on the following figure, which illustrates our situation:

\begin{figure}[H]
\centerline{\includegraphics[scale=0.8]{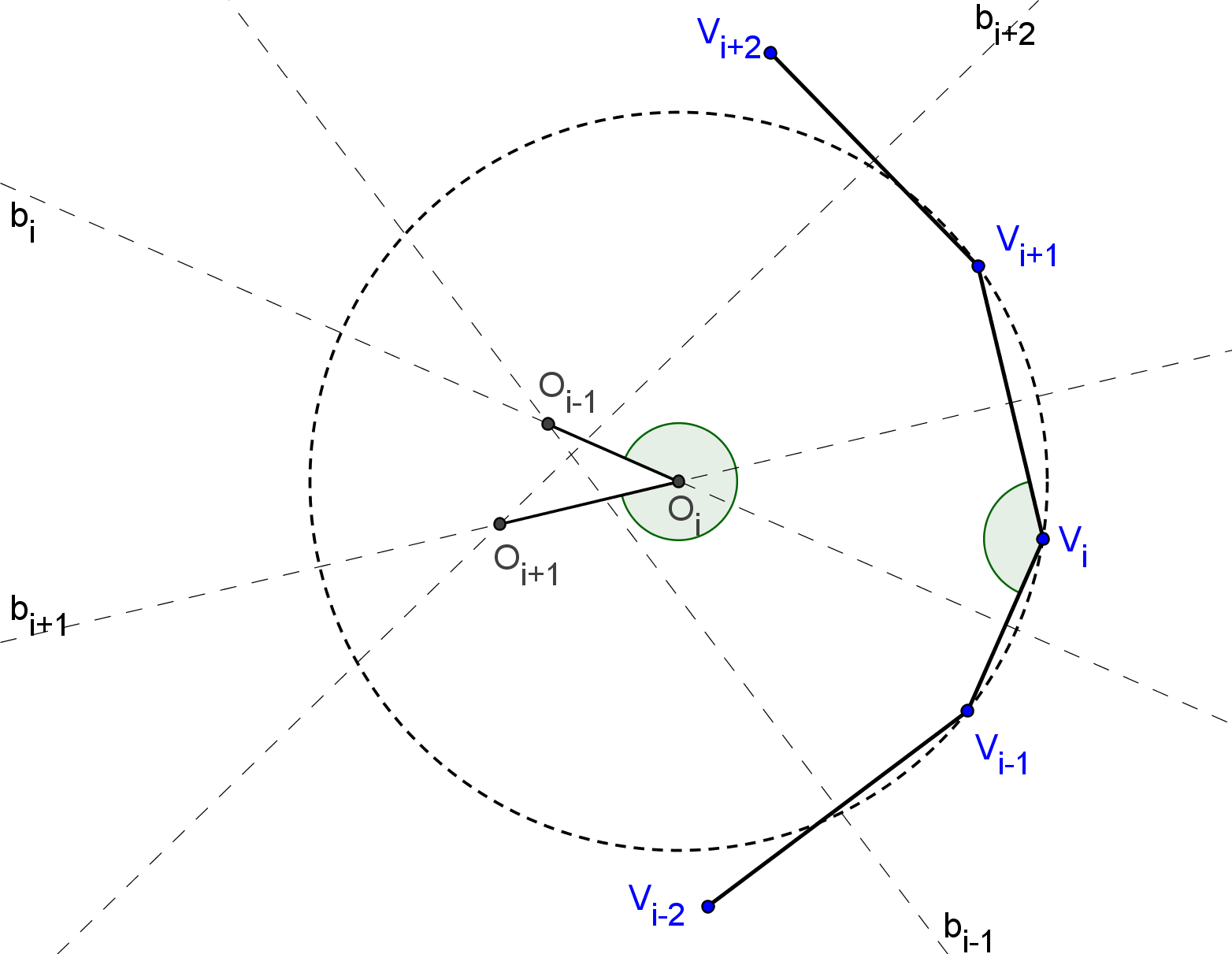}}
\caption[Proof of Lemma \ref{lemma:evoluteangles} (1)]{}
\end{figure}

We observe that the green angles are our angles of consideration. We see that since $b_{i}$ and $b_{i+1}$ are perpendicular to the segments of our polygonal line, the angle formed by the intersection of these two line segments is supplementary to $\angle V_{i}.$
\vskip1mm
By our convention mentioned earlier, we always consider the left angle with respect to the orientation. We notice that we are traveling counterclockwise on our polygonal line and observe that the position of the vertices $V_{i-2}$ and $V_{i+2}$ affects how the perpendicular bisectors intersect.
\vskip1mm
So now, we know that the right angle at the vertex $O_{i}$ is supplementary to $\angle V_{i}$, and since we are considering the left angle at the vertex $O_{i}$, we have the following:
$$\angle O_{i}=2\pi -(\pi -\angle V_{i})=\pi +\angle V_{i}$$
From this, our assertion follows.
\vskip2mm
Now, lets consider another subcase. Assume that $V_{i-1}$, $V_{i}$, and $V_{i+1}$ are negative and $V_{i-2}$ and $V_{i+2}$ lie inside of the circle $C_{i}.$ Now, constructing everything as before, we have the following figure:

\begin{figure}[H]
\centerline{\includegraphics[scale=1.2]{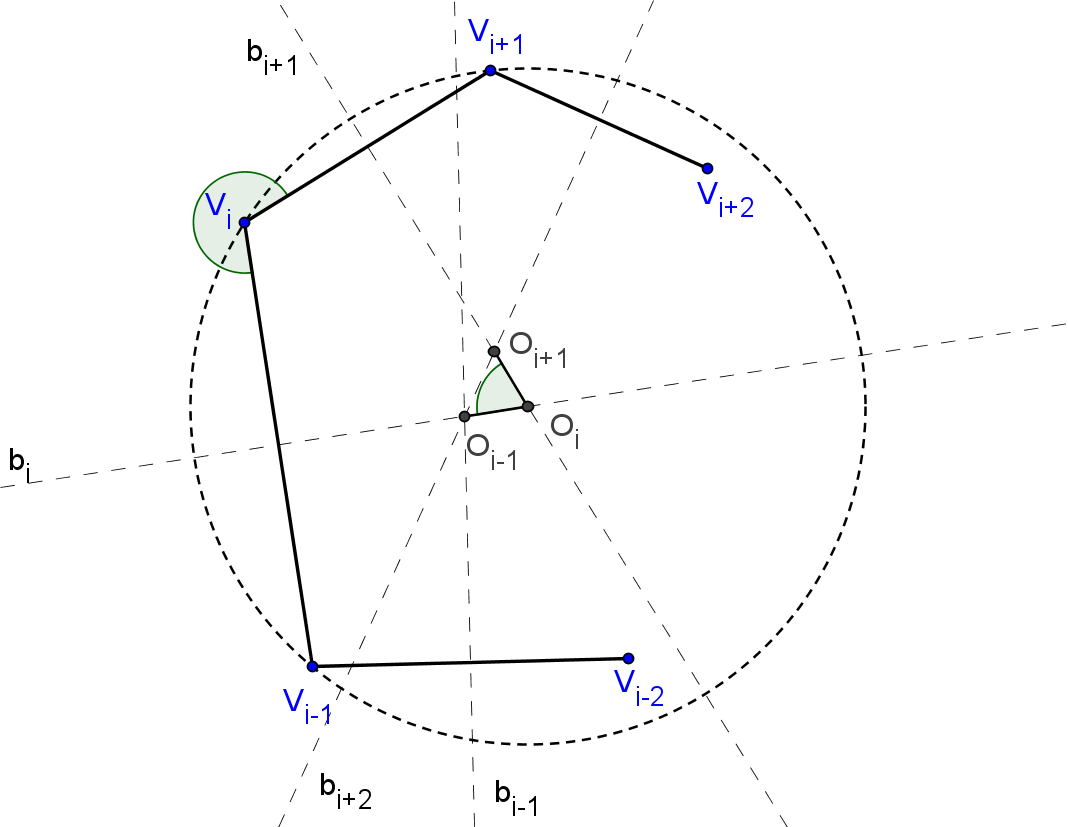}}
\caption[Proof of Lemma \ref{lemma:evoluteangles} (2)]{}
\end{figure}

Observe now that the orientation of our polygonal line is clockwise and the vertices $V_{i-2}$ and $V_{i+2}$ lie inside the circle $C_{i}$. So, our situation has changed. Our left angle $\angle V_{i}$ is no longer supplementary to the one formed by the intersection of the corresponding bisectors, which is the one on the evolute under consideration. But, the angle to the right of $V_{i}$ is. So we have the following:
$$\angle O_{i}=\pi -(2\pi-\angle V_{i})=\angle V_{i}-\pi$$
From this, our assertion follows.

\vskip3mm
Case 2: $V_{i+1}\succ V_{i} \prec V_{i-1}.$
\vskip2mm
For our first subcase, assume that $V_{i-1}$, $V_{i}$, and $V_{i+1}$ are positive and $V_{i-2}$ and $V_{i+2}$ lie inside of the circle $C_{i}.$ Working similarly as we did the other case, we have the following situation:

\begin{figure}[H]
\centerline{\includegraphics[scale=1.1]{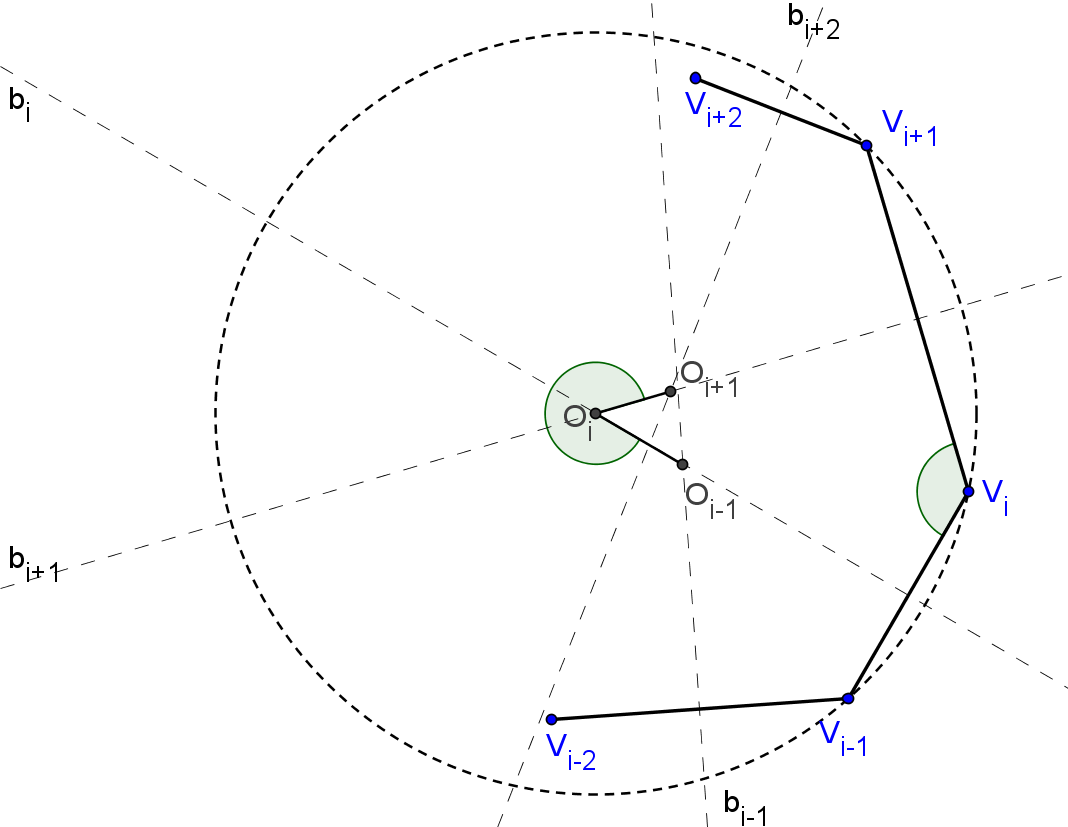}}
\caption[Proof of Lemma \ref{lemma:evoluteangles} (3)]{}
\end{figure}

Observe that the orientation of our polygonal line in this case is counterclockwise. We also observe that in this case, $\angle V_{i}$ is supplementary to the angle formed by the corresponding perpendicular bisectors. Hence, we have the following:
$$\angle O_{i}=2\pi -(\pi-\angle V_{i})=\angle V_{i}+\pi,$$
which gives us the desired conclusion.
\vskip2mm
Now we consider another subcase. Assume that all vertices are negative and that vertices $V_{i-2}$ and $V_{i+2}$ lie outside the circle $C_{i}.$ Working similarly as the other cases, we have the following situation:

\begin{figure}[H]
\centerline{\includegraphics[scale=1.4]{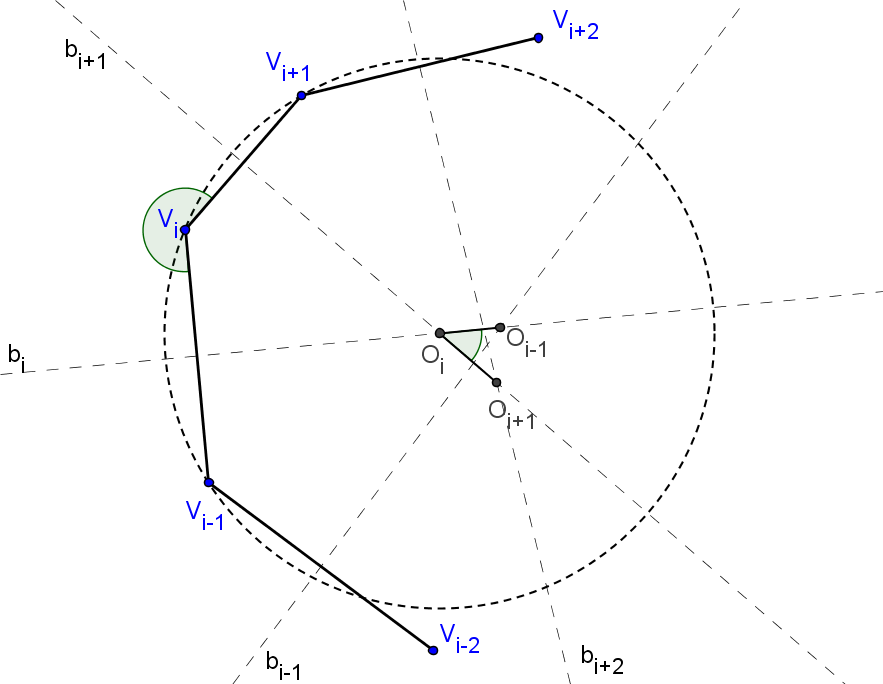}}
\caption[Proof of Lemma \ref{lemma:evoluteangles} (4)]{}
\end{figure}

The orientation of our polygonal line is clockwise in this case, hence we have:
$$\angle O_{i}=\pi-(2\pi-\angle V_{i})=\angle V_{i}-\pi,$$
which gives us our assertion.
\vskip2mm

Now we focus our attention on the second portion of our lemma and look at two cases where our vertex is not extremal.
\vskip2mm
Case 1:
\vskip2mm
Assume that all vertices are positive, and that the vertex $V_{i-2}$ lies inside the circle $C_{i},$ while $V_{i+2}$ lies outside. The following figure illustrates out situation:

\begin{figure}[H]
\centerline{\includegraphics[scale=1.35]{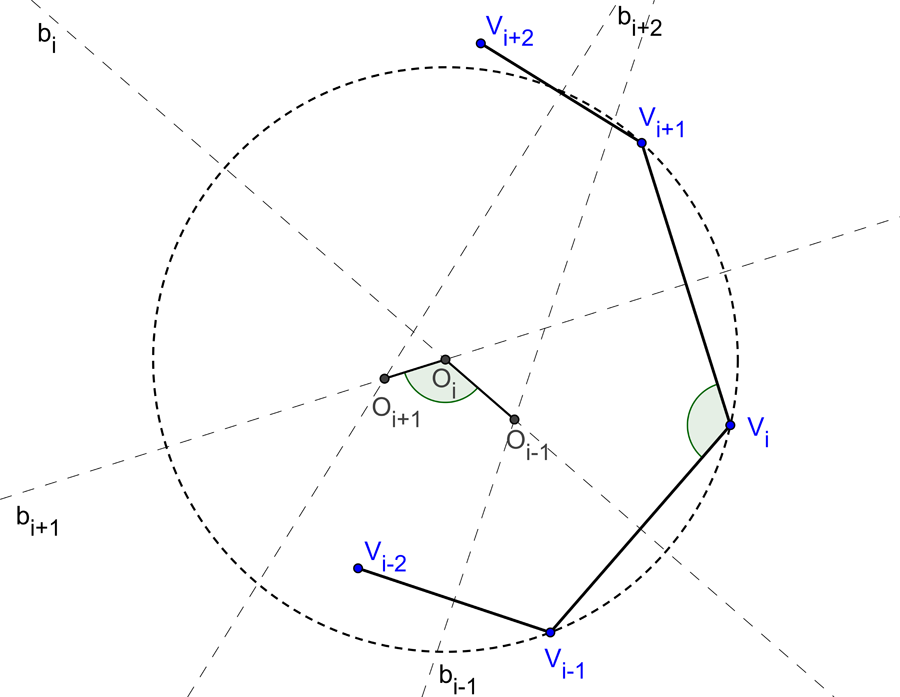}}
\caption[Proof of Lemma \ref{lemma:evoluteangles} (5)]{}
\end{figure}

We clearly see that since $\angle V_{i}$ is supplementary to the angle formed by the intersection of the corresponding perpendicular bisectors, which is supplementary to $\angle O_{i}$. So, it follows that:
$$\angle V_{i}=\angle O_{i},$$
giving us the desired conclusion.
\vskip2mm
Case 2:
\vskip2mm
Assume that all vertices are negative, and that the vertex $V_{i-2}$ lies inside the circle $C_{i},$ while $V_{i+2}$ lies outside. The following figure illustrates out situation:

\begin{figure}[H]
\centerline{\includegraphics[scale=1.35]{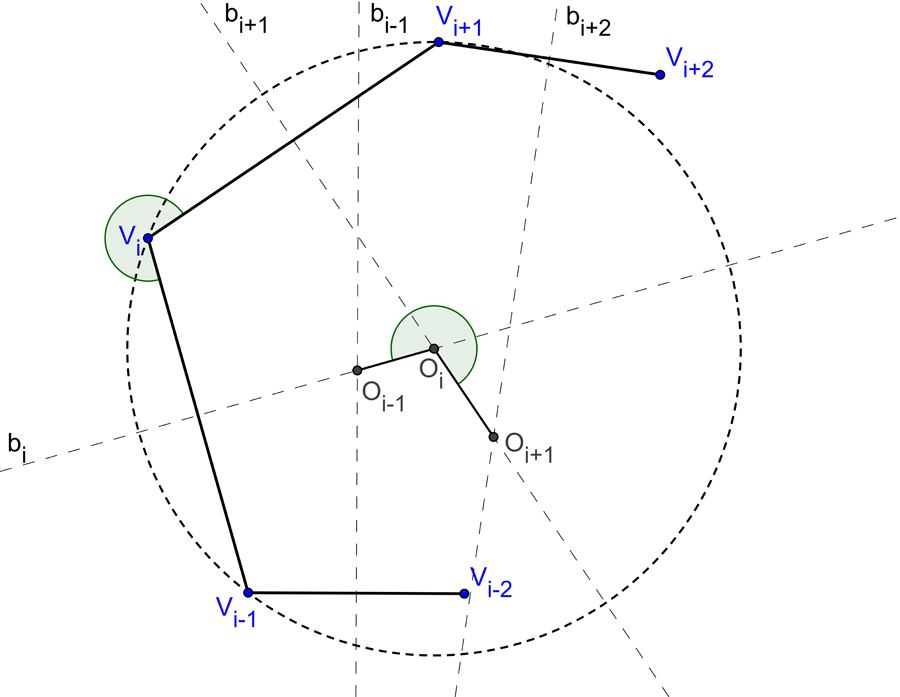}}
\caption[Proof of Lemma \ref{lemma:evoluteangles} (6)]{}
\end{figure}

We observe that we have the following:
$$\angle O_{i}=2\pi-(2\pi-\angle V_{i})=\angle V_{i},$$
which proves our assertion.
\vskip2mm
The proof of the reamaining 26 cases follows exactly the same pattern as these cases, so we will not consider them.
\end{proof}

\subsection{The Winding Number}

\begin{defn}
Let $P$ be a polygonal line with vertices $V_1,V_2,...,V_n.$ Then the discrete winding number $wind(P)$ of $P$ is defined to be:
$$wind(P)=\frac{1}{2\pi}\sum_{i}^n(\pi-\angle V_{i}).$$
\end{defn}

The following figure shows the angles which are considered when computing the discrete winding number:

\begin{figure}[H]
\centerline{\includegraphics[width=100mm,height=81mm]{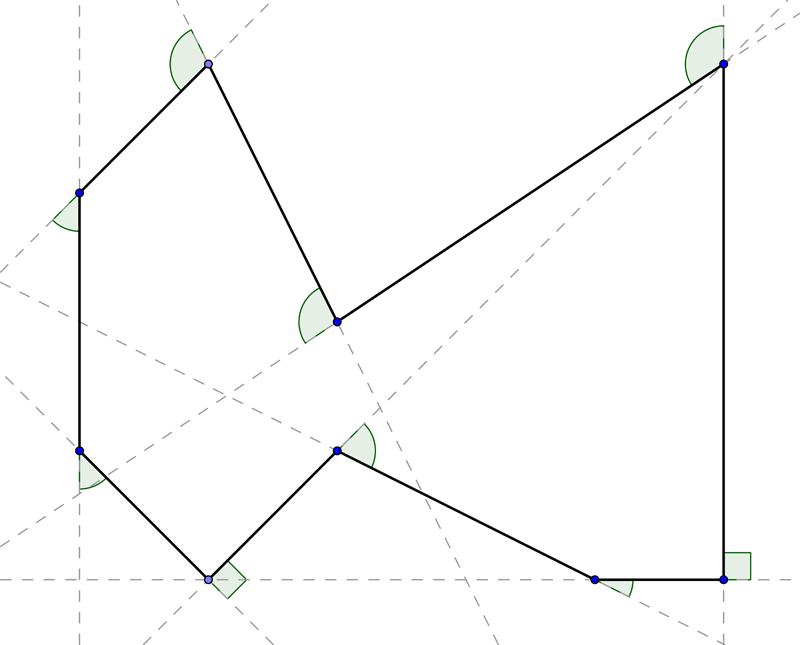}}
\caption[Computing the Winding Number]{}
\end{figure}

In this example, we add all of the angles at positive vertices, and then subtract the angles at the negative vertices.
The discrete winding number (as well as the smooth winding number) is always an integer. In our example, and more generally for all simple polygons, the discrete winding number is always equal to 1. This can be easily seen with a triangulation argument.

\begin{defn}
Let $P$ be a smooth curve. Then:
$$wind(P)=\frac{1}{2\pi}\int k(s)ds,$$
where $k(s)$ is the curvature function of the curve.
\end{defn}
\vskip2mm

It is a well known fact that every smooth curve can be uniformly approximated by a polygonal line. So, we need to see that by finer and finer approximations and passage of the limit, these definitions agree. So, we have to show that:
$$\lim_{n\rightarrow\infty}\sum_{i=1}^{n}(\pi-\angle V_{i})=\int k(s)ds.$$
But, it turns out that this is a result from differential geometry known as Fox-Milnor's Theorem \cite{mil}, so indeed our definitions agree.
\vskip2mm
Now, lets divert our attention to the evolute. For the discrete case, since the evolute of a closed polygon is a closed polygon, the discrete definition of the winding number we gave does hold. But, for the smooth case, the definition we gave does not hold. We notice that when the curvature vanishes on our original curve, the radius of curvature goes to infinity, leaving an asymptote in our evolute. Also we notice that at cusps of the evolute, the tangent vector is not defined. The following figure illustrates a basic curve for which we encounter such problems:

\begin{figure}[H]
\centerline{\includegraphics[width=100mm,height=75mm]{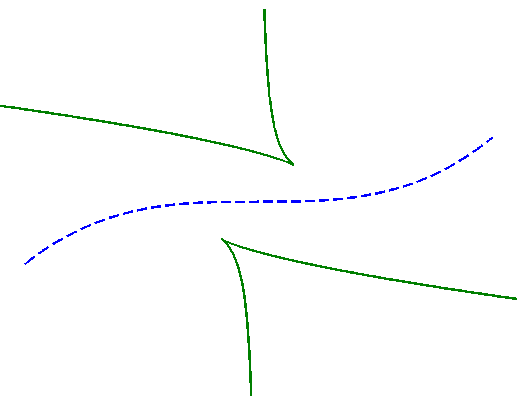}}
\caption[An Example of an Evolute with an Asymptote]{}
\end{figure}

In this figure we let $-1\leq t\leq 1$ and

$$\begin{array}{ccc}
  x(t) & = & t \\
  y(t) & = & t^{3}
\end{array}$$

This example illustrates that whenever we have an inflection point on our curve, our evolute has an asymptote. This leads us to the following rigorous definition for the winding number of the evolute:

\begin{defn}
$$wind(E(P))=\frac{1}{2\pi}\lim_{n\rightarrow\infty}\sum_{i=1}^{k_{n}}(\pi-\angle P_{i}),$$
where $Q_{n}=P_{1}P_{2}...P_{k_{n}}$ are polygonal approximations of $E(P).$
\end{defn}

\subsection{The Local and Smooth Four-Vertex Theorem and the Evolute}

In this section, we will prove a result that relates the number of locally extremal vertices to the winding numbers of a curve and the corresponding evolute.

\begin{theorem}
\label{thm:evolutethm}
$$N_{+} - N_{-} = 2wind(P) - 2wind(E(P)),$$ where $P$ is a smooth curve or a polygonal curve.
\end{theorem}

\begin{proof}
We first consider the case of a closed polygonal curve with vertices $V_{1},V_{2},...,V_{n}$.
By Lemma \ref{lemma:evoluteangles}, we see that if a vertex $V_{i}$ is extremal, then
$$\angle V_{i} - \angle O_{i}=\pm\pi$$
and if the vertex is not extremal, then
$$\angle V_{i} - \angle O_{i}=0$$
So,
$$2wind(P)-2wind(E(P))=\frac{1}{\pi}\sum_{i=1}^{n}(\pi-\angle V_{i})-\frac{1}{\pi}\sum_{i=1}^{n}(\pi-\angle O_{i})=$$
$$=\frac{1}{\pi}\sum_{i=1}^{n}(\angle O_{i}-\angle V_{i})=\frac{1}{\pi}(\pi N_{+}-\pi N_{-})=N_{+}-N_{-}$$
\vskip2mm
For the case of a smooth planar curve, this theorem follows immediately by passage of the limit. Since we proved it for polygonal curves, we have that this theorem holds for every polygonal curve approximating the smooth curve. Now if we pass to the limit, Fox-Milnor's Theorem \cite{mil} gives us the desired integral, and by definition of the winding number of the evolute of a smooth curve, we have our desired equality.
\end{proof}
\vskip1mm
Now, we will see what this theorem tells us about the winding number of the evolute of a curve.

\begin{cor}
\label{cor:windevolute}
Let $P$ be a simple polygon or a simple closed planar curve. Then
$$wind(E(P))=\frac{2-(N_+ - N_{-})}{2}.$$
If $P$ is convex, then
$$wind(E(P))=\frac{2-N_+}{2}.$$
\end{cor}
\begin{proof}
For a simple polygon, the first part of our statement simply follows from the fact that $wind(P)=1$ for simple polygons. For a smooth simple closed planar curve $P$, we know that all such curves are homotopic to a circle, which has winding number equal to 1. Hence, $wind(P)=1.$ Combining these facts with Theorem \ref{thm:evolutethm}, we get our assertion. The second part of the statement simply follows from our convexity assumption, which gives us that $N_{-}=0.$
\end{proof}

We now recall the Four-Vertex Theorem for the smooth case.

\begin{theorem}[The Smooth Four-Vertex Theorem]
Every smooth simple convex closed planar curve has at least four vertices.
\end{theorem}

We now observe how the Local and Smooth Four-Vertex Theorems are related to the winding number of the evolute. We see that in fact, Corollary \ref{thm:evolutethm} tells us that the $wind(E(P))$ depends on the number of locally extremal vertices. In fact, we have that $wind(E(P))<0$ if and only if $P$ has at least four locally extremal vertices. So, by proving the Local Four-Vertex Theorem, we essentially are proving that $wind(E(P))<0$.

\section{The Decomposition of Polygons}
\renewcommand{\thefigure}{\thesection.\arabic{figure}}
\setcounter{figure}{0}
In this section we will investigate a very natural concept: the decomposition of a polygon into two smaller polygons. We will try to see what kind of impact this notion has on our various types of extremal vertices. For the remainder of this section, we need to introduce one more definition.

\begin{defn}
We call an edge of a polygon a Delaunay edge if there exists an empty circle passing through the corresponding vertices of that edge. If there exists a full circle passing through this edge, then we call the edge an Anti-Delaunay edge.
\end{defn}

So what exactly does it mean to decompose a polygon? Here, the notion of decomposing a polygon will simply be the cutting of a polygon $P$ by passing a line segment through any two vertices so that the line segment lies in the interior of the polygon. We will call this line segment a \textit{diagonal}. Also, we will denote the two new polygons formed by a decomposition by $P_1$ and $P_2$ and also require that they each have at least four vertices. By this last condition, it automatically follows that $P$ must have at least six vertices to successfully perform a decomposition.

Now that we have introduced this new concept, we search for some relationships between the number of various extremal vertices of the original polygon $P$ and the two smaller polygons ($P_1$ and $P_2$) obtained after a decomposition. At first glance, we would believe that the number of extremal vertices of $P$ is less than the sum of the extremal vertices of the two smaller polygons for every type of extremality. But, alas, this is not true. Consider the following figure:

\begin{figure}[H]
\centerline{\includegraphics[scale=0.8]{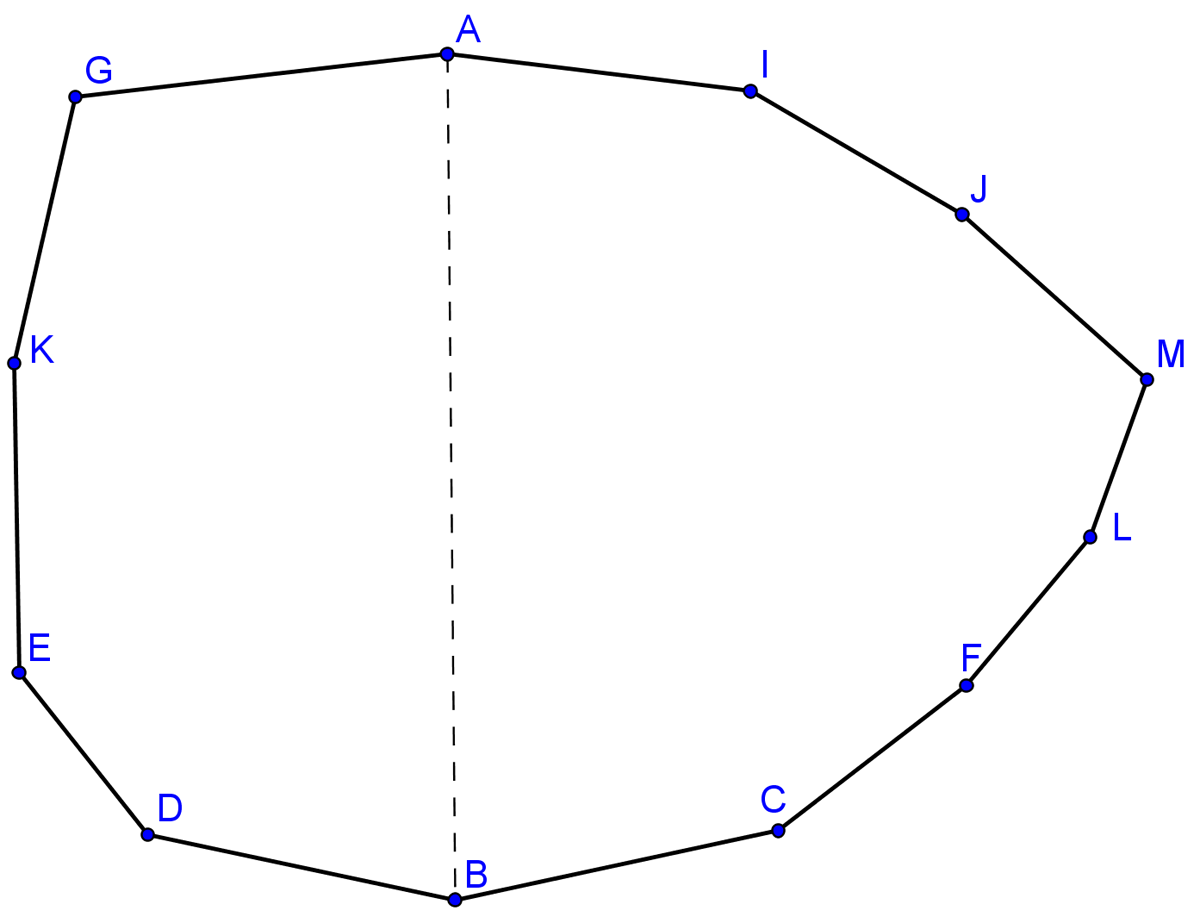}}
\caption[A Counterexample to $s_{-}(P)\leq s_{-}(P_1)+s_{-}(P_2)$]{}
\end{figure}

This is a very delicate counterexample. Let $P_1$ be the polygon on the left after the cutting, and $P_2$ be the polygon on the right. For the polygon $P,$ we have maximal globally extremal vertices $G,D,C,M$ and $I.$ So, $s_{-}(P)=5.$ For $P_1,$ only $G$ and $D$ are maximal, and for $P_2$ only $B$ and $M$ are maximal. So it follows that $s_{-}(P_1)=s_{-}(P_2)=2.$ So, we see this does not hold for globally extremal vertices.
\vskip1mm
It is interesting to note that for this example we have exactly the same situation for locally extremal vertices as for global, but everything holds for radially extremal vertices. So the next natural question to ask is whether there exists an example for which the sum of extremal vertices of $P_1$ and $P_2$ is greater than the sum of extremal vertices of $P$ for local extremality, but not for global. This question was answered by Arseniy Akopian, who found the following example:

\begin{figure}[H]
\centerline{\includegraphics[scale=0.7]{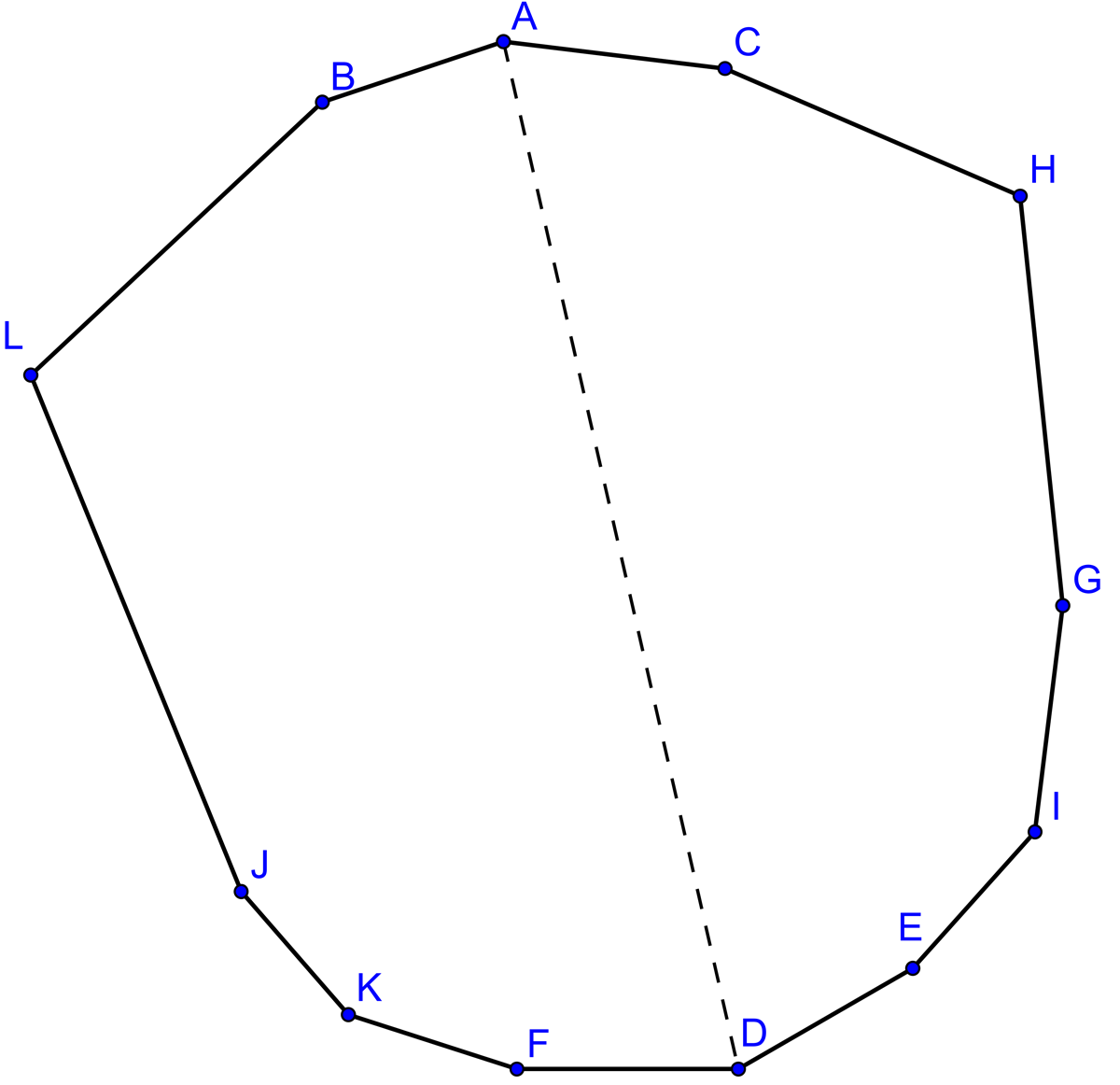}}
\caption[An Example Where an Inequality Holds for Local Extremality, But Not for Global Extremality]{}
\end{figure}

For Akopian's example this is not true for globally extremal vertices, but surprisingly holds for locally extremal and radially extremal vertices, even with equality. This tells us that, if we would like to derive some concrete inequalities, we need to consider each of the cases of extremality \emph{separately}.

We now will prove some results about the decomposition of polygons and the impact it has on different types of extremal vertices.

\subsection{Decomposition and Globally Extremal Vertices}

In this section, we will simply restrict our situation to globally extremal vertices. We are particularly interested on how many we gain or lose when performing a decomposition, and see if we can put a bound on this.

\begin{lemma}
\label{lemma:decompglobal}
Let $P$ be a generic convex polygon with 6 vertices and let $P_1$ and $P_2$ be the resulting polygons of a decomposition. Then
$$s_{-}(P)\geq s_{-}(P_1)+s_{-}(P_2)-2.$$
\end{lemma}

\begin{proof}
Due to the fact that we have a 6-gon and no two maximal vertices can be next to each other, we see that $P$ can have at most 3 globally maximal-extremal vertices. Also, we see that if we decompose $P$ into two smaller polygons $P_1$ and $P_2$, they each will always have four vertices. Let $A$ and $D$ be the vertices of the cutting diagonal, with $B$ and $C$ neighbors of $A$ and $E$ and $F$ neighbors of $D$, respectively. So, $P_1=ABED$ and $P_2=ACFD$. The following figure illustrates our situation:

\begin{figure}[H]
\centerline{\includegraphics[scale=0.9]{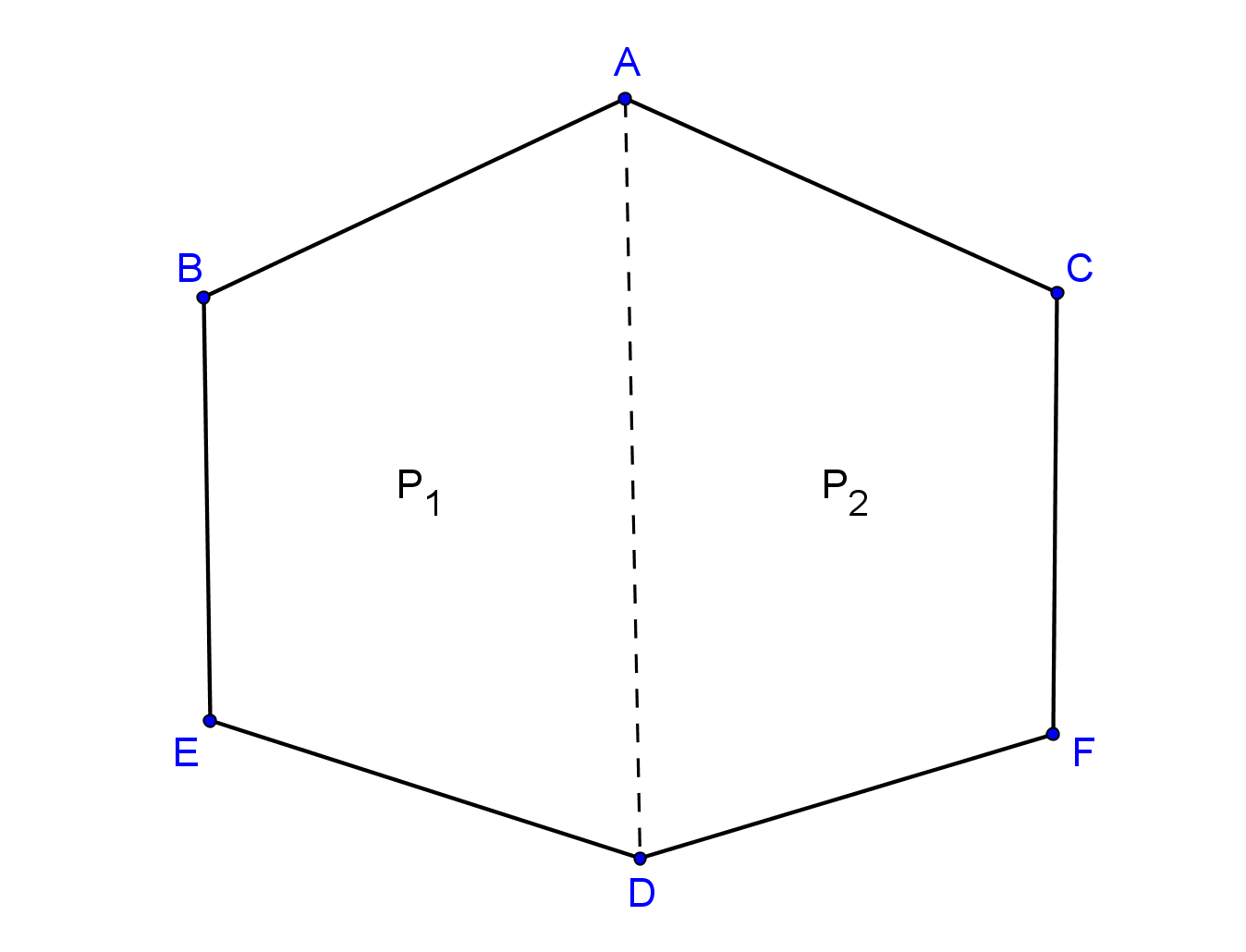}}
\caption[Proving Lemma \ref{lemma:decompglobal}]{}
\end{figure}

We now simply run through all possible cases. We first assume that $P$ has three globally maximal-extremal vertices. It follows that one must lie on the cutting diagonal and the other two must be neighbors of the adjacent vertex on the diagonal. Without loss of generality, say the vertex on the diagonal is $A$ and the other two vertices are $E$ and $F$. Then, $E$ is maximal-extremal for $P_1$ and $F$ is maximal-extremal for $P_2$. Our inequality holds, since we have that $A$ could be extremal for both $P_1$ and $P_2$, for only one of them, or for neither of them.

Now we assume that $P$ has exactly two globally maximal-extremal vertices. Then we have three possibilities, either both are neighbors of a vertex of our cutting diagonal, both lie on our cutting diagonal, or one lies on a diagonal vertex and one is the neighbor of the opposite diagonal vertex. A similar routine checking for each of the cases as the one above yields our inequality.

For the case where $P$ has exactly one globally maximal vertex, we have that this vertex either must lie on our cutting diagonal or on one of the neighboring vertices. For the former, assume that $A$ is extremal. Then, since none of our neighboring vertices can be extremal, it follows that $A$ or $D$ is extremal for $P_1$ and similarly for $P_2$. So, our inequality follows. For the latter case, we see if the neighboring vertices are extremal in $P$, then they are extremal in one of the smaller polygons, without loss of generality say $P_1$. Since we can achieve at most one more maximal vertex in $P_1$ and at most two in $P_2$, our inequality follows.
\end{proof}

\begin{remark}
Notice that if we had assumed the Global Four-Vertex Theorem, the proof of the above lemma would be much simpler, since we would have less cases to check. For a reason that we will disclose later, we assumed that we have no such theorem.
\end{remark}

\begin{theorem}
\label{thm:globaldecomp}
Let $P$ be a generic convex polygon with at least 6 vertices and let $P_1$ and $P_2$ be the resulting polygons of a decomposition. Then
$$s_{-}(P)\geq s_{-}(P_1)+s_{-}(P_2)-3.$$
\end{theorem}

\begin{remark}
We will first approach the proof with a much stronger inequality, $s_{-}(P)\geq s_{-}(P_1)+s_{-}(P_2)-2$, since this would be much more desirable. Then, at the end will see one particular case where our approach fails, forcing us to weaken our inequality.
\end{remark}

\begin{proof}
Before we begin, we turn our attention to a lemma from Section 2, Lemma \ref{lemma:globalinductionmax}. In particular, we focus our attention on what was proven. We essentially showed that, if we make such a cutting of a maximal-extremal vertex, apply induction to $P'$, and then proceed to add the maximal-extremal vertex to $P'$ to retrieve $P$, then we will lose \emph{at most} one maximal vertex of $P'.$ So, in essence, we showed that either
\vskip2mm
$s_{-}(P)=s_{-}(P')$ or $s_{-}(P)=s_{-}(P')+1.\hskip20mm (1)$
\vskip2mm
We will use this fact heavily, since we will apply a similar induction for the proof of this theorem. We will find an extremal vertex, remove it, and assume that this new polygon satisfies our stronger inequality. We then will re-attach our vertex and retrieve our original polygon, and show that our inequality holds.
\vskip2mm
Lemma \ref{lemma:decompglobal} takes care of the case where $n=6.$ So we let $P$ be a convex polygon with at least 7 vertices. Proposition \ref{prop:existenceprop} guarantees us an extremal vertex, so without loss of generality we assume that $V_i$ is this vertex. We let $V_{i-1}$ and $V_{i+1}$ be the neighboring vertices of $V_i$. We now remove our vertex $V_i$ from our polygon $P$ and join vertices $V_{i-1}$ and $V_{i+1}$ by a line segment to obtain a new polygon $P'$. There are two ways we can now decompose our polygon $P'$: either we use a diagonal passing through either the vertex $V_{i-1}$ or $V_{i+1}$, or we use a diagonal passing through any other vertices of $P'$.
\vskip2mm
\noindent\textit{Case 1:} We use any other vertices of $P'$.
\vskip1mm
So now, our diagonal does not pass through either $V_{i-1}$ or $V_{i+1}$. We obtain two polygons, $P_1'$ and $P_2'$. We observe that if we were to add the vertex $V_i$ back to the polygon $P'$, then we would either be adding it to $P_1'$ or $P_2'$. Without loss of generality, assume that we are adding it to $P_2'$ and denote the new polygon obtained after adding $V_i$ by $P_2$. The following figure illustrates our situation:

\begin{figure}[H]
\centerline{\includegraphics{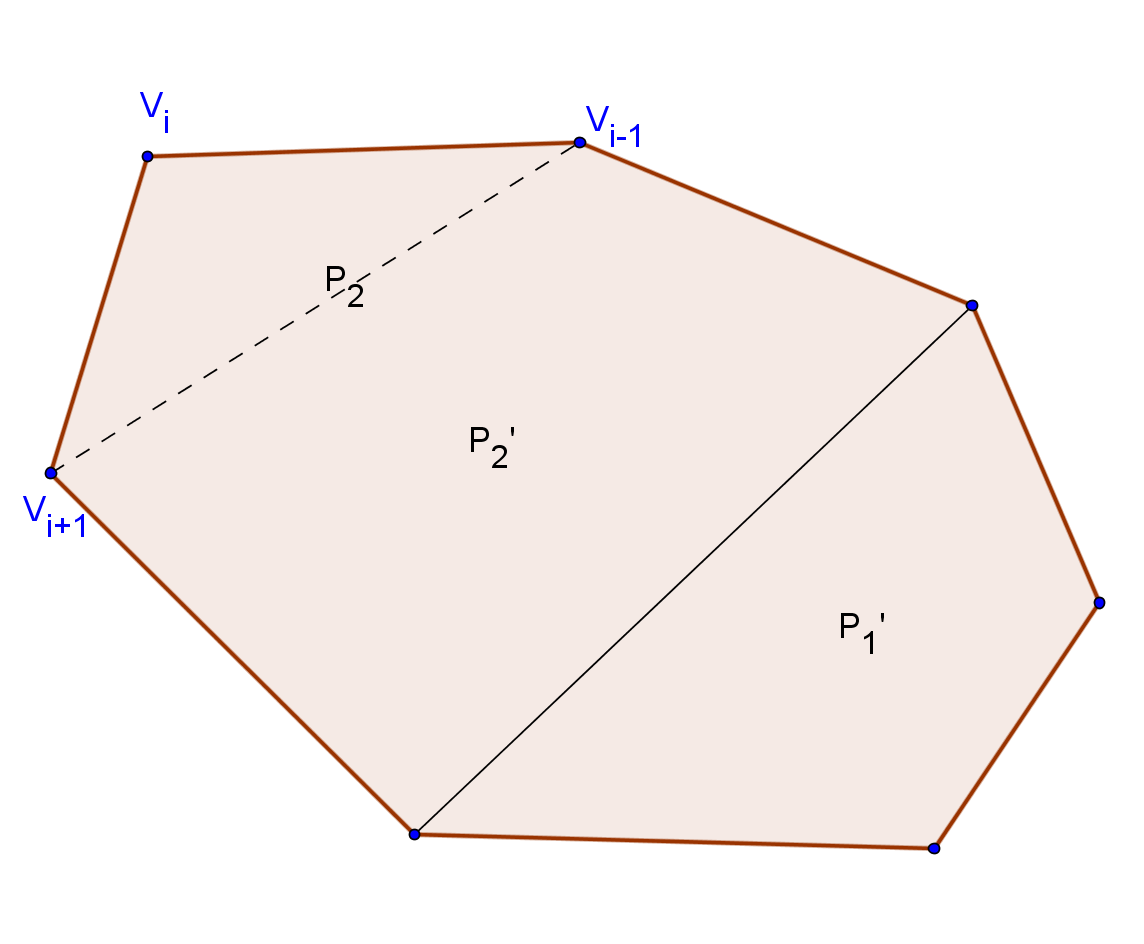}}
\caption[Proving Theorem \ref{thm:globaldecomp} (1)]{}
\end{figure}

Before we apply our inductive hypothesis to $P'$, we make a quick observation. We note that if $V_i$ is maximal-extremal for $P,$ then it is maximal-extremal for $P_2.$ This is because the corresponding empty circle $C_i$ will stay empty after any cutting not involving a diagonal passing through $V_i.$
\vskip1mm
We now apply our induction hypothesis to $P'$ and obtain that $$s_{-}(P')\geq s_{-}(P_1')+s_{-}(P_2')-2.$$ Combining our observation above that $V_i$ is extremal for $P_2$ and (1), we obtain that either
\vskip2mm
$s_{-}(P_2)=s_{-}(P_2')$ or $s_{-}(P_2)=s_{-}(P_2')+1.\hskip20mm (2)$
\vskip2mm
Applying the same idea to $P$ and $P',$ we obtain that either
\vskip2mm
$s_{-}(P)=s_{-}(P')$ or $s_{-}(P)=s_{-}(P')+1.\hskip22mm (3)$
\vskip2mm
We now see that combining (2) and (3), we obtain several cases. After examining these cases, we see that there is only one that gives us a problem, the case where
\vskip2mm
$s_{-}(P_2)=s_{-}(P_2')+1$ and $s_{-}(P)=s_{-}(P').$
\vskip2mm
If this holds, then we cannot achieve that $s_{-}(P)\geq s_{-}(P_1')+s_{-}(P_2)-2,$ so we must show that such a situation is not possible. Indeed, since our diagonal does not pass through $V_{i-1}$ or $V_{i+1}$, it follows that if either $V_{i-1}$ or $V_{i+1}$ is maximal for $P'$, then it is maximal for $P_2'.$ Then contrapositive of this statement shows that such a situation was not possible, and \textit{Case 1} is proved.
\vskip2mm
\noindent\textit{Case 2:} Our diagonal passes through either $V_{i-1}$ or $V_{i+1}$.
\vskip1mm
It follows since our vertex $V_i$ is maximal-extremal in $P$, then it is maximal-extremal in $P_2$, and necessarily the vertices $V_{i-1}$ and $V_{i+1}$ are not maximal-extremal in either $P$ or $P_2.$ Without loss of generality, we assume that our diagonal passes through $V_{i+1}$. The following figure illustrates our situation:

\begin{figure}[H]
\centerline{\includegraphics{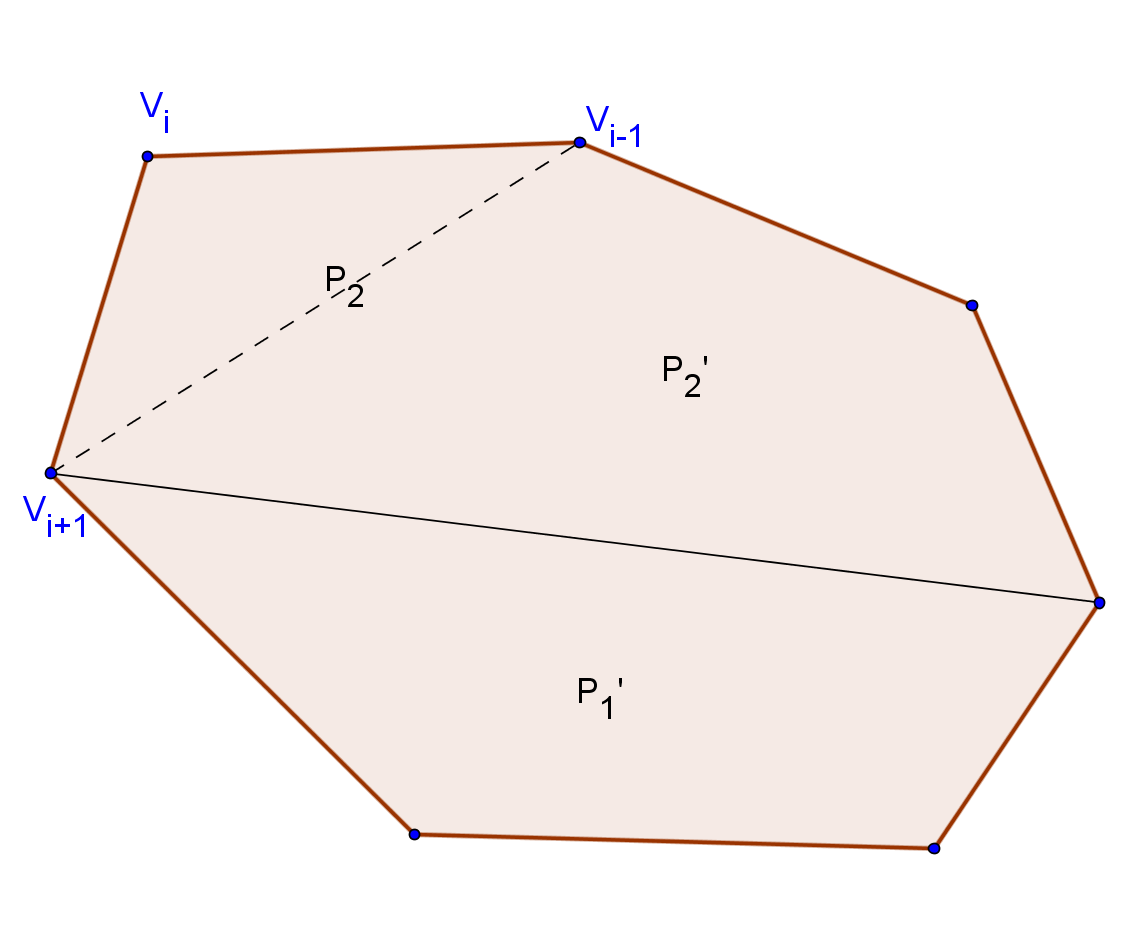}}
\caption[Proving Theorem \ref{thm:globaldecomp} (2)]{}
\end{figure}

By induction, we have that $s_{-}(P')\geq s_{-}(P_1')+s_{-}(P_2')-2$. We must now run through all the cases of maximal-extremality and non-maximal-extremality of the vertices $V_{i-1}$ and $V_{i+1}$ in the smaller polygons, $P'$ and $P_2'$.

\vskip2mm
\noindent\textit{Subcase 1:} $V_{i+1}$ is not maximal in $P'$.
\vskip1mm

We first assume that $V_{i+1}$ is not maximal in $P'$. Then we have two cases, either $V_{i-1}$ can be either maximal or not maximal in $P'$.
\vskip2mm
Assume that $V_{i-1}$ is maximal in $P'$. It follows that $s_{-}(P)=s_{-}(P')$ and that $V_{i-1}$ is also maximal in $P_2'$. Now we closely observe our situation on the polygon $P_2'$. Since $V_{i-1}$ is maximal in $P_2'$, then $V_{i+1}$ cannot be maximal in $P_2'$.  But, it then follows that $s_{-}(P_2)=s_{-}(P_2')$, and applying this to our induction hypothesis, we are done.
\vskip2mm

We now assume that $V_{i-1}$ is not maximal in $P'$. Since $V_{i-1}$ is not maximal in $P'$ while $V_{i+1}$ is maximal in $P'$, it follows that $s_{-}(P)=s_{-}(P')+1.$ We now have two possibilities, either $V_{i+1}$ is maximal or not maximal in $P_2'.$ If it is maximal in $P_2'$, then $V_{i-1}$ is not maximal in $P_2'$, and necessarily we have that $s_{-}(P_2)=s_{-}(P_2')$. If $V_{i+1}$ is not maximal in $P_2'$, then it is possibly for $V_{i-1}$ to be either maximal or not maximal in $P_2'$. If it is maximal, then $s_{-}(P_2)=s_{-}(P_2')$, and if it is not maximal, then  $s_{-}(P_2)=s_{-}(P_2')+1$. In all situations discussed, our inequality follows.
\vskip2mm
\noindent\textit{Subcase 2:} $V_{i+1}$ is maximal in $P'$.
\vskip1mm
From our assumption, it automatically follows that $V_{i-1}$ cannot be maximal in $P'$ and that $s_{-}(P)=s_{-}(P')$.  We now have two possibilities, either $V_{i+1}$ is maximal or not maximal in $P_2'.$ If it is maximal in $P_2'$, then $V_{i-1}$ is not maximal in $P_2'$, and we have that $s_{-}(P_2)=s_{-}(P_2')$. If $V_{i+1}$, then we have that either $V_{i-1}$ is maximal or not maximal in $P_2'.$ If it is maximal then $s_{-}(P_2)=s_{-}(P_2')$, and our assertion follows. If $V_{i-1}$, then by our reasoning it follows that $s_{-}(P_2)=s_{-}(P_2')+1$, which cannot happen if we wish to achieve the stronger inequality $s_{-}(P')\geq s_{-}(P_1')+s_{-}(P_2')-2$. Unfortunately, this situation is geometrically feasible, forcing us to weaken our inequality to $s_{-}(P')\geq s_{-}(P_1')+s_{-}(P_2')-3$.
\end{proof}

The following figure illustrates the situation where $s_{-}(P)\geq s_{-}(P_1)+s_{-}(P_2)-3$:

\begin{figure}[H]
\centerline{\includegraphics[scale=0.9]{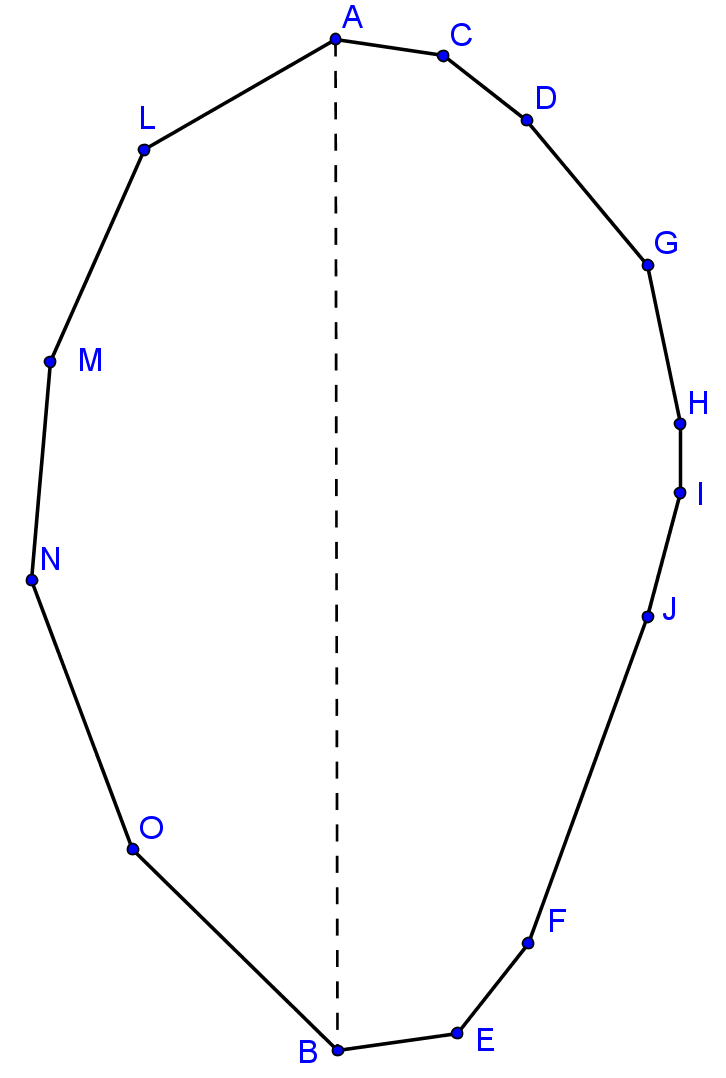}}
\caption[An Example Where $s_{-}(P)\geq s_{-}(P_1)+s_{-}(P_2)-3$]{}
\end{figure}

In our above figure, set $V_i=C$, $V_{i-1}=D$ and $V_{i+1}=A$. Also, call the polygon on the left $P_1$ and the one on the right $P_2$. We obtain $P'$ and $P_2'$ by cutting the vertex $C.$ Now, we are in the situation of our above proof. It is not too difficult to check that $A$ is maximal in $P'$ and not maximal in $P_2'$. Also, $D$ is not maximal in $P'$ and $P_2'$. So, we are actually in the ``bad" case of our proof. After rigorously checking, we can verify that the globally extremal vertices of $P$ are $C$, $E$, and $G$, of $P_1$ are $L$ and $B$, and of $P_2$ are $C$, $G$, $I$ and $E$. So, indeed $s_{-}(P)=3\geq s_{-}(P_1)+s_{-}(P_2)-3=2+4-3=3$.
\vskip2mm
So why did the stronger inequality not hold? The problem lies in the fact that when decomposing a polygon into two smaller polygons, not only are the vertices of the cutting diagonal as well as the neighboring vertices affected, but also the remaining vertices can be affected. While the maximal-extremal vertices of the original polygon which are not on the cutting diagonal stay maximal-extremal in the smaller polygons, it is possible that non-extremal vertices in the larger polygon can become maximal-extremal in the smaller polygons after a decomposition. Hence the choice of an inductive approach in Theorem \ref{thm:globaldecomp}. It turns out that with an extra assumption, we can actually obtain our stronger inequality.

\begin{theorem}
\label{theorem:strongglobal}
Let $P$ be a generic convex polygon with six or more vertices. Assume that the cutting diagonal of a decomposition is Delaunay. Then
$$s_{-}(P)\geq s_{-}(P_1)+s_{-}(P_2)-2.$$
\end{theorem}
\begin{proof}
The worst case scenario would be if we had no globally maximal-extremal vertices on on our edge. By applying a Delaunay triangulation to $P_1$ and $P_2$ we see that our result follows immediately. If we would like to avoid Delaunay triangulation, it turns out that we can prove this result using the same technique as in the proof of Lemma \ref{lemma:globalinductionmax}.
\end{proof}

Now, recall from Section 2.2 that we do not necessarily have the same number of globally maximal-extremal vertices as globally minimal-extremal vertices. Observe that the essence of proving the above result was Lemma \ref{lemma:globalinductionmax}. We proved an analogous lemma for globally minimal-extremal vertices, Lemma \ref{lemma:globalinductionmin}. So, we have a similar inequality for globally minimal-extremal vertices.

\begin{theorem}
\label{theorem:strongglobal2}
Let $P$ be a generic convex polygon with six or more vertices. Assume that the cutting diagonal of a decomposition is Anti-Delaunay. Then
$$s_{+}(P)\geq s_{+}(P_1)+s_{+}(P_2)-2.$$
\end{theorem}

\subsection{Decomposition and Locally Extremal Vertices}

Now that we have investigated decomposition and globally extremal vertices, a natural question is, what happens with locally extremal vertices? Luckily, we will not have the situation described above. In fact, it is easy to see that the only vertices that will be affected by a decomposition of a polygon will be the vertices on the cutting diagonal and the neighboring vertices. Which means that we have a total six vertices that are impacted by a decomposition, which leads us to a feasible case by case analysis. In this section, we will heavily be using Proposition \ref{prop:circleprop}.

\begin{lemma}
\label{lemma:localdec1}
Let $P$ be a generic convex polygon and $B$ and $D$ the vertices of a cutting diagonal. Let $A$ and $C$ be the neighbors of $B$ in $P$ and let $P_1$ and $P_2$ be the polygons obtained after a decomposition, with $P_1$ possessing vertex $A$ and $P_2$ possessing vertex $C$. Assume that $A$ is locally maximal-extremal in $P_1$ but not in $P$, and that $C$ is locally maximal-extremal for $P_2$ but not in $P$. Then, $B$ is a locally maximal-extremal vertex for $P$.
\end{lemma}
\begin{proof}
Let $X$ be the neighbor of $A$ in $P_1$ and $Y$ be the neighbor of $C$ in $P_2$. Denote the circle passing through vertices $A$, $B$ and $C$ by $C_B$, the circle passing through vertices $X$, $A$ and $B$ by $C_A$, and the circle passing through vertices $B$, $C$ and $Y$ by $C_C$.
Since $A$ is not maximal-extremal in $P$, it follows that $A$ lies inside the circle $C_C$. By Proposition \ref{prop:circleprop}, it follows that $Y$ lies outside of the circle $C_B$. Since $C$ is not maximal-extremal in $P$, it follows that $C$ lies inside the circle $C_A$. By Proposition \ref{prop:circleprop}, it follows that $X$ lies outside of the circle $C_B$. The following figure illustrates our situation:

\begin{figure}[H]
\centerline{\includegraphics[scale=0.65]{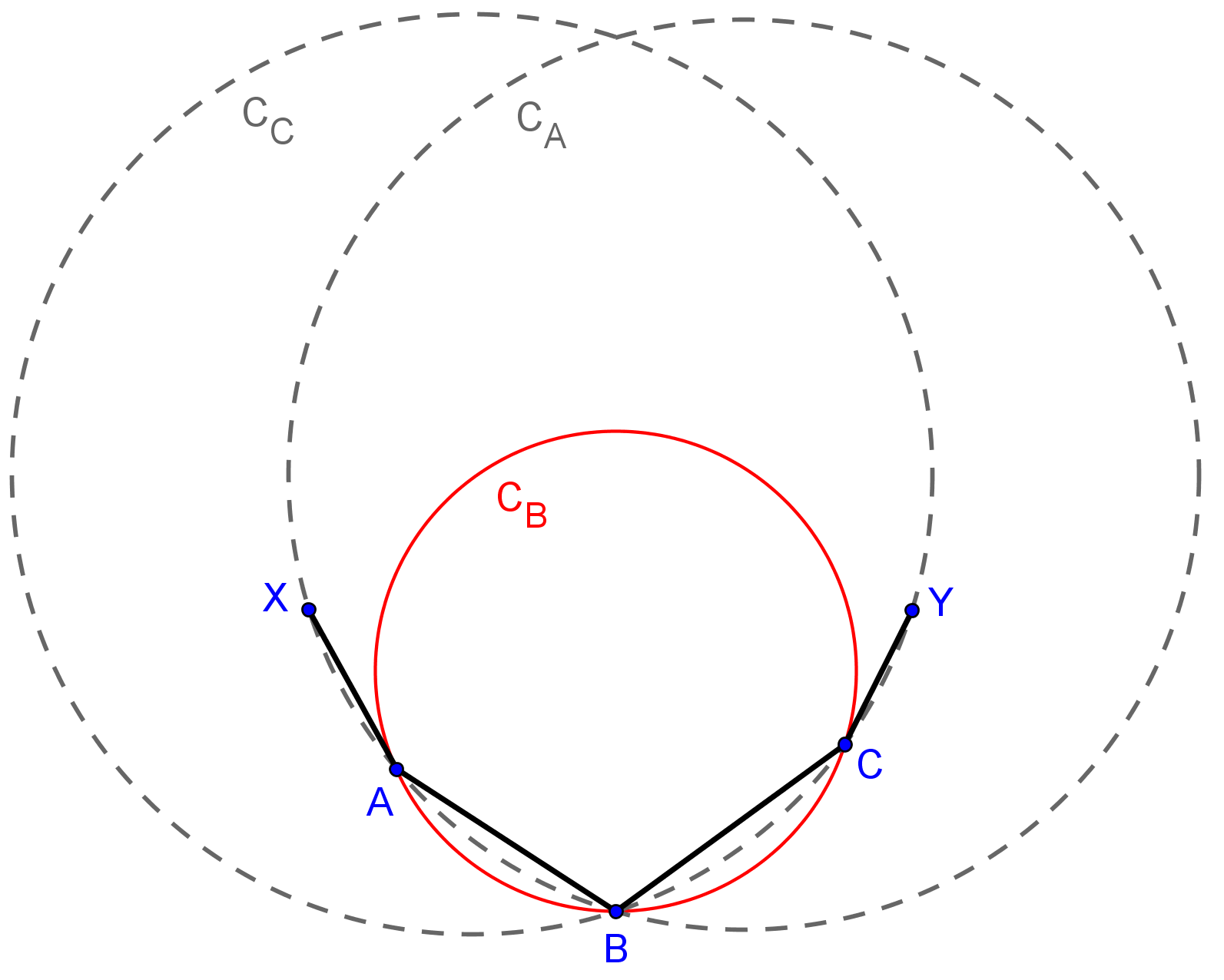}}
\caption[Proving Lemma \ref{lemma:localdec1}]{}
\end{figure}

So, it follows that, since both $X$ and $Y$ lie outside of the circle $C_B$, $B$ is maximal-extremal in $P$.
\end{proof}

\begin{lemma}
\label{lemma:localdec2}
Let $P$ be a generic convex polygon and $B$ and $D$ the vertices of a cutting diagonal. Let $A$ and $C$ be the neighbors of $B$ in $P$ and let $P_1$ and $P_2$ be the polygons obtained after a decomposition, with $P_1$ possessing vertex $A$ and $P_2$ possessing vertex $C$. Assume that $A$ is locally maximal-extremal in $P_1$ but not in $P$, and that $B$ is locally maximal-extremal in $P_2$. Then, $B$ is locally maximal-extremal in $P$.
\end{lemma}
\begin{proof}
For simplicity, consider the following figure, which will illustrate our configuration of points and circles:

\begin{figure}[H]
\centerline{\includegraphics[scale=0.75]{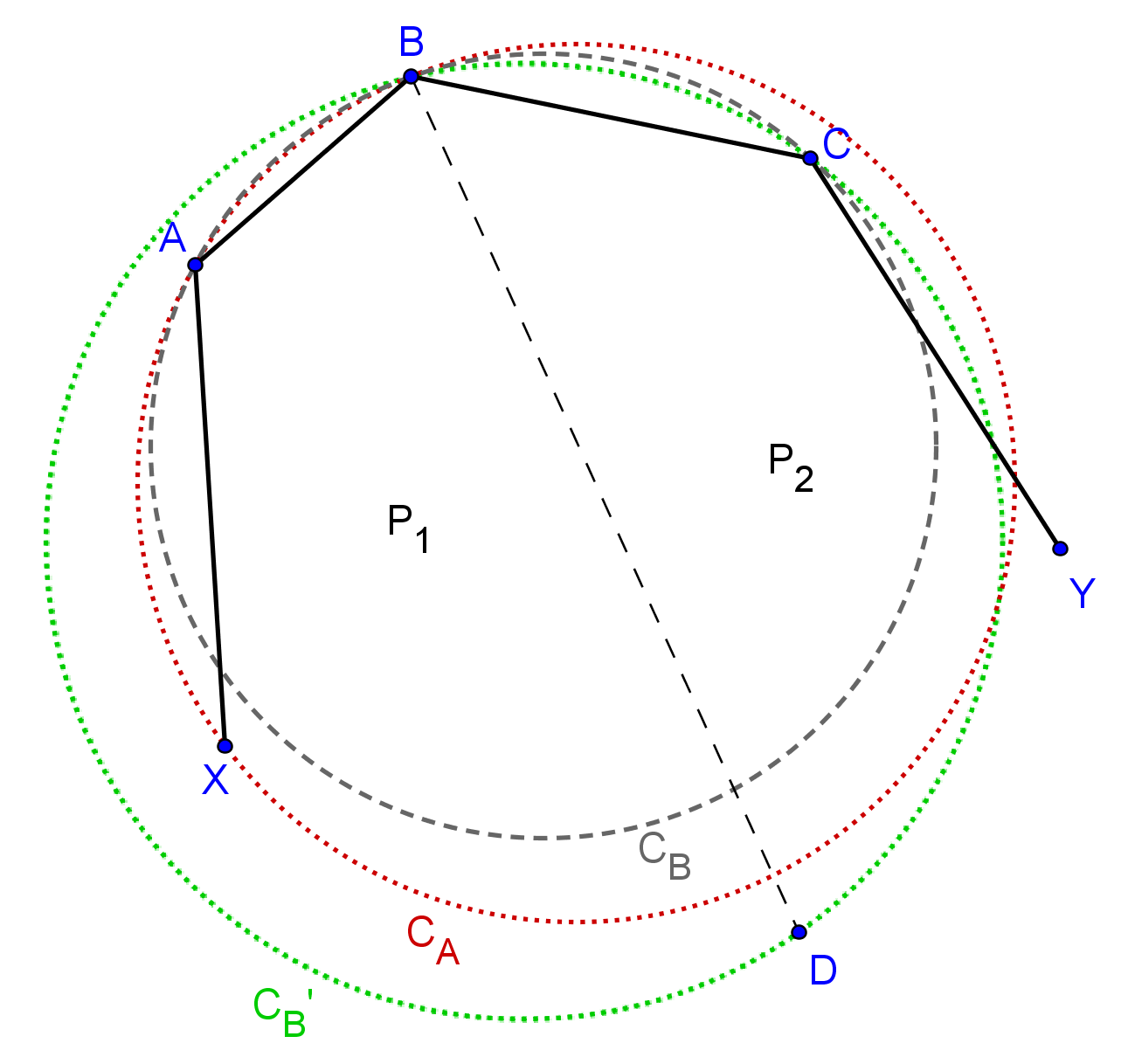}}
\caption[Proving Lemma \ref{lemma:localdec2}]{}
\end{figure}

Let $X$ be the neighbor of $A$ and $Y$ the neighbor of $C$. Denote by $C_A$ the circle passing through vertices $X$, $A$, and $B$. Since $A$ is maximal-extremal in $P_1$, it follows that $D$ lies outside of the circle $C_A$. Since $A$ is not maximal-extremal in $P$, it follows that $C$ must lie inside the circle $C_A$. Now, denote the circle passing through vertices $A$, $B$, and $C$ by $C_B$. Our goal is to show that vertices $X$ and $Y$ lie outside of the circle $C_B$.

We first will show that $X$ lies outside of $C_B$. Denote by $H^{+}_{AB}$ the half plane formed by the infinite line $AB$ and containing the point $X$. By convexity, it follows that $C$ lies in $H^{+}_{AB}$. Since $C$ lies inside circle $C_A$, Proposition \ref{prop:circleprop} yields that $X$ lies outside of circle $C_B$.

Now, we will show that $Y$ lies outside of the circle $C_B$. Denote by $C_B'$ the circle passing through the points $C$, $B$ and $D$. We will show that if $Y$ lies outside of $C_B'$, then it lies outside of $C_B$. To do this, we first must show that $A$ lies inside the circle $C_B'$.

Consider the circles $C_A$ and $C_B'$. These circles intersect at two points, point $B$ and some other point, say $Z$. The infinite line $BZ$ divides the plane into two half-planes. By our convex arrangement of points, it follows that points $A$ and $D$ lie in the same half-plane. Since $D$ lies outside of the circle $C_A$, it follows by Proposition \ref{prop:circleprop} that $A$ lies inside the circle $C_B'$.

Lastly, consider the circles $C_B$ and $C_B'$. These two circles intersect at the points $B$ and $C$. The infinite line $BC$ divides the plane into two half-planes, and our convexity assumption yields that points $A$ and $D$ lie in the same half-plane. Since $A$ lies inside the circle $C_B'$, it follows from Proposition \ref{prop:circleprop} that $D$ lies outside of the circle $C_B$.

Now, since $B$ is maximal-extremal in $P_2$, it follows that $Y$ lies outside of $C_B'$. By our above observation, it follows immediately from Proposition \ref{prop:circleprop} that $Y$ lies outside of $C_B$, since $Y$ lies in the same half-plane as $A$ and $D$ with respect to the infinite line $BC$. Since points $X$ and $Y$ both lie outside of the circle $C_B$, it follows that $B$ is maximal-extremal in $P$.
\end{proof}

\begin{lemma}
\label{lemma:localdec3}
Let $P$ be a generic convex polygon and $B$ and $D$ the vertices of a cutting diagonal. Let $A$ and $C$ be the neighbors of $B$ in $P$ and let $P_1$ and $P_2$ be the polygons obtained after a decomposition, with $P_1$ possessing vertex $A$ and $P_2$ possessing vertex $C$. Assume that $A$ is locally maximal-extremal for $P_1$ and $D$ is locally maximal-extremal for both $P_1$ and $P_2$, but not for $P$. Then $A$ is locally maximal-extremal for $P$.
\end{lemma}
\begin{proof}
Let $X$ be the neighbor of $A$ in $P_1$, $E$ be the neighbor of $D$ in $P_1$, and $F$ be the neighbor of $D$ in $P_2$. Denote by $C_{D1}$ the circle passing through vertices $B$, $D$ and $E$, by $C_{D2}$ the circle passing through vertices $B$, $E$ and $F$, and by $C_A$ the circle passing through vertices $X$, $A$ and $B$. The following figure illustrates our configuration:

\begin{figure}[H]
\centerline{\includegraphics[scale=1.1]{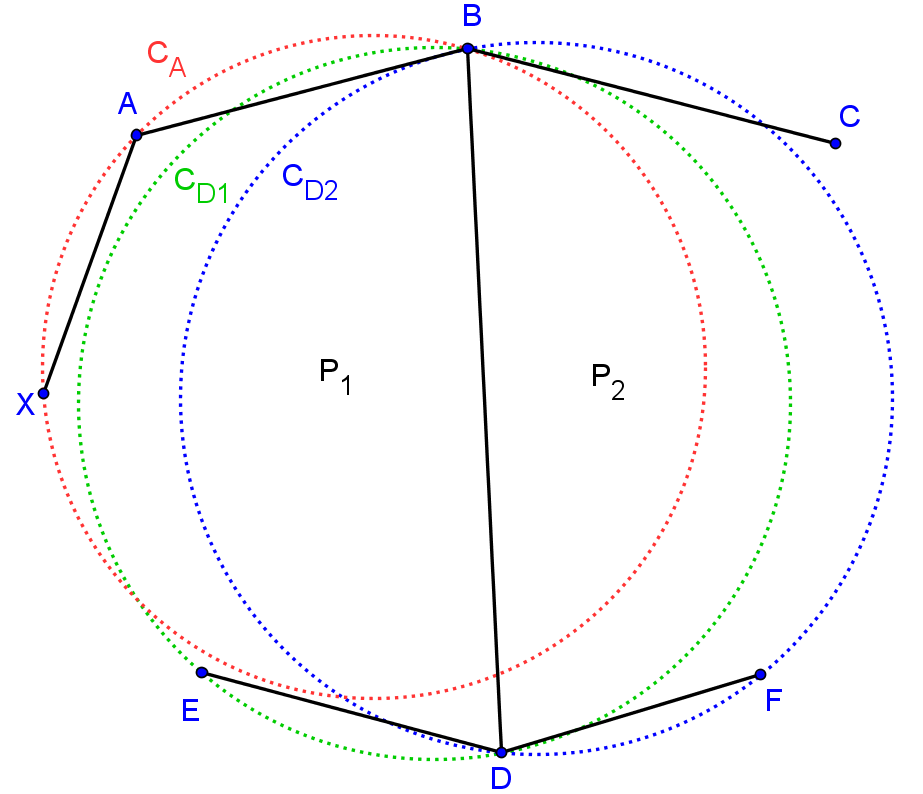}}
\caption[Proving Lemma \ref{lemma:localdec3} (1)]{}
\end{figure}

Our goal is to show that vertex $C$ lies outside of the circle $C_A$. We will do this by showing that if $C$ lies outside the circle $C_{D2}$, then it also lies outside of circle $C_A$. Since $A$ is maximal-extremal in $P_1$, it follows that $D$ lies outside of $C_A$. Since $D$ is maximal-extremal in $P_1$, it follows that $A$ lies outside of circle $C_{D1}$.

We will first show that if $C$ would lie outside of $C_{D1}$, then $C$ lies outside of $C_A$. This requires a clever application of Proposition \ref{prop:circleprop}. Consider the points of intersection of the circle $C_A$ with the circle $C_{D1}$. We know that both circles intersect at $B$ and some second point, say $Y$, and that the infinite line $BY$ divides the plane into two half-planes, one containing the point $D$ ($H^{+}_{BY}$) and the other one not ($H^{-}_{BY}$). Since $A$ lies outside of $C_{D1}$ and $D$ lies outside of $C_A$, it follows that $A$ must lie in $H^{-}_{BY}$.

By convexity, we have that $C$ lies in $H^{+}_{BY}$. It now follows from Proposition \ref{prop:circleprop} that if $C$ lies outside of the circle passing through points $Y$, $B$ and $D$, then it lies outside of $C_A$. But, the circle passing through $Y$, $B$ and $D$ is the circle $C_{D1}$. The following figure illustrates this situation:

\begin{figure}[H]
\centerline{\includegraphics[scale=1.1]{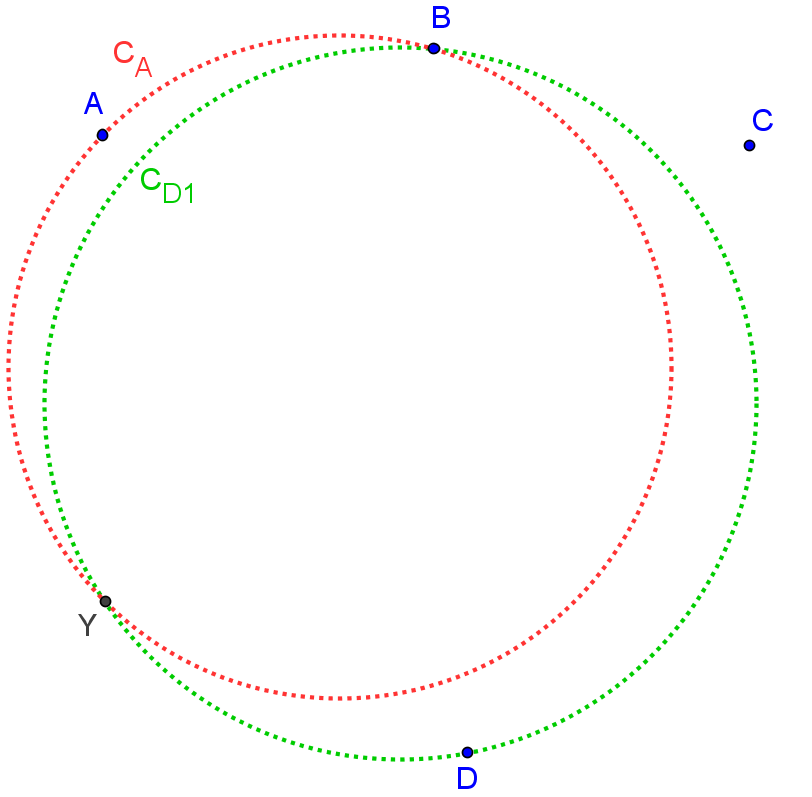}}
\caption[Proving Lemma \ref{lemma:localdec3} (2)]{}
\end{figure}

Now we consider circles $C_{D1}$ and $C_{D2}$. Since $D$ is maximal-extremal in $P_2$, it follows that $C$ lies outside of the circle $C_{D2}$. If we show that $C$ also lies outside of $C_{D1}$, then we are done. To do this, we will heavily use the fact that $D$ is not maximal-extremal in $P$. We will show that if $E$ lies inside the circle $C_{D2}$ or if $F$ lies in $C_{D1}$, then $D$ is maximal-extremal in $P$, contradicting our assumption.

It is enough just to check this for $E$. Denote the circle passing through vertices $E$, $D$ and $F$ by $C_D$. If $E$ lies inside the circle $C_{D2}$, then by Proposition \ref{prop:circleprop}, it follows that $B$ lies outside of the circle $C_D$, since both $E$ and $B$ lie in the same half-plane with respect to the infinite line $DF$. It also follows by Proposition \ref{prop:circleprop} that $F$ lies inside the circle $C_{D1}$, since $F$ lies in a different half-plane than $E$ with respect to the infinite line $BD$. Now denote by $E'$ the neighbor of $E$ and by $F'$ the neighbor of $F$. The following figure illustrates the situation:

\begin{figure}[H]
\centerline{\includegraphics[scale=1.2]{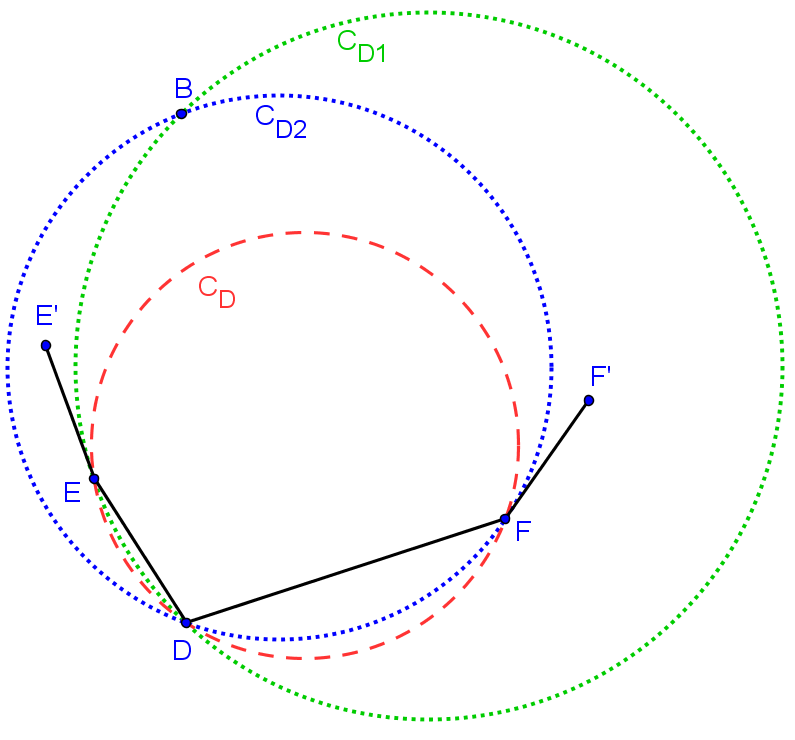}}
\caption[Proving Lemma \ref{lemma:localdec3} (3)]{}
\end{figure}

Since $D$ is maximal-extremal in $P_1$, it follows that $E'$ lies outside of the circle $C_{D1}$. Similarly, since $D$ is maximal-extremal in $P_2$, it follows that $F'$ lies outside the circle $C_{D2}$. Now, recall that $B$ lies outside of the circle $C_D$. Since points $E'$ and $B$ lie on the same half-plane with respect to the infinite line $ED$, it follows by Proposition \ref{prop:circleprop} that $E'$ lies outside of $C_D$. Since points $B$ and $F'$ lie in the same half-plane with respect to the infinite line $DF$, it follows that $F'$ also lies outside of $C_D$. So, we obtain that $D$ is maximal-extremal in $P$, a contradiction.

So now we know that $E$ must lie outside of the circle $C_{D2}$. By Proposition \ref{prop:circleprop}, we see that $F$ lies outside of the circle $C_{D1}$, since $F$ lies in a different half-plane than $E$ with respect to the infinite line $BD$. So, if $C$ were to lie outside of circle $C_{D2}$, then it would also lie outside of the circle $C_{D1}$. But earlier we proved that if $C$ would lie outside of circle $C_{D1}$, then $C$ would lie outside of the circle $C_A$. Indeed, by assumption, $C$ lies outside of $C_{D2}$ and hence outside of $C_A$. Since $A$ was maximal-extremal in $P_1$, it also follows that $X$ lies outside of the circle $C_A$. Therefore $A$ is maximal-extremal in $P$.
\end{proof}

\begin{theorem}
\label{thm:localineq}
Let $P$ be a generic convex polygon with at least 6 vertices and let $P_1$ and $P_2$ be the resulting polygons of a decomposition. Then
$$l_{-}(P)\geq l_{-}(P_1)+l_{-}(P_2)-2.$$
\end{theorem}
\begin{proof}
We note that only six vertices are affected by a decomposition from the local point of view: the vertices of the cutting diagonal and the neighbors of those vertices. So, we will want to eliminate the cases which violate our inequality. Denote the 6 vertices effected by $A$, $B$, $C$, $D$, $E$ and $F$ so that after decomposition, $P_1=ABED$ and $P_2=CBEF$, as shown in the following figure:

\begin{figure}[H]
\centerline{\includegraphics{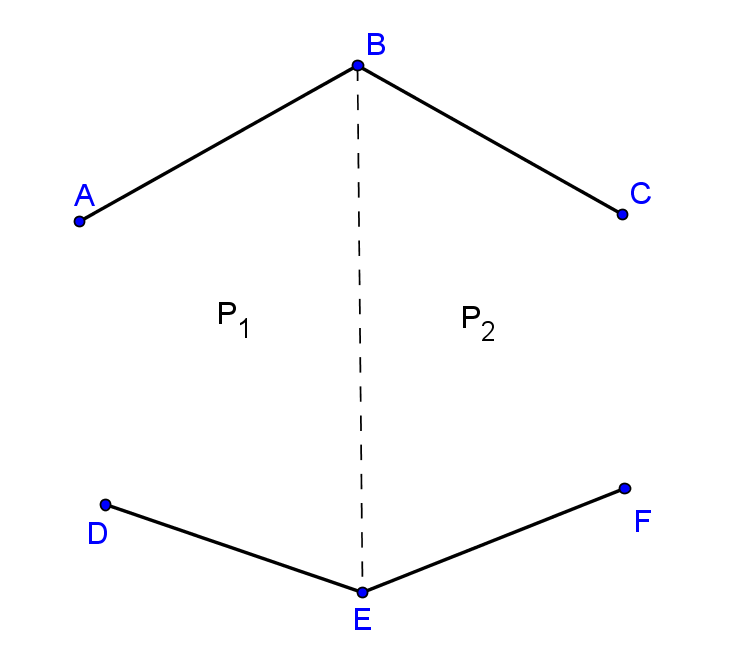}}
\caption[Proving Theorem \ref{thm:localineq}]{}
\end{figure}

\noindent\textit{Case 1:} We gain two maximal-extremal vertices in $P_1$, as well as $P_2$, but none of the six vertices are maximal-extremal in $P$.
\vskip1mm
We now will run through all possible configurations in $P_1$, and see that we cannot have two maximal vertices in $P_2$, since this will violate at least one of our three lemmas from this section. By the symmetry of our cases, without loss of generality, it is sufficient for us to check two configurations in $P_1$.

Firstly, we assume that $A$ and $D$ are maximal-extremal in $P_1$. It then follows by Lemma \ref{lemma:localdec1} and Lemma \ref{lemma:localdec2} that neither $C$, $B$, $E$ or $F$ can be maximal-extremal for $P_2$, because if any of those four vertices were maximal-extremal in $P_2$, then we risk having $B$ or $E$ as maximal vertices for $P$, or even both. So, such a situation is not possible.

Secondly, we assume that $A$ and $E$ are maximal-extremal in $P_1$. It follows from \ref{lemma:localdec1} and Lemma \ref{lemma:localdec2} that $F$ and $C$ cannot be maximal-extremal in $P_2$. It then follows that only $B$ or $E$ can be maximal-extremal in $P_2$, but not both since $P$ was assumed to be generic. So, such a case is not possible.

\vskip2mm
\noindent\textit{Case 2:} We gain two maximal-extremal vertices in $P_1$ and gain two maximal-extremal vertex in $P_2$, and one of the six vertices is maximal-extremal in $P$.
\vskip1mm
Without loss of generality and by the symmetry of our cases, it will be sufficient just to check the cases where $B$ or $C$ is maximal-extremal in $P$.
\vskip2mm
\noindent\textit{Subcase 1:} $B$ is maximal in $P$.
\vskip1mm
It is sufficient to check all the possible configurations in $P_1$ and see that we cannot have two maximal-extremal vertices in $P_2$. By the symmetry of our cases, it is sufficient to check only three configurations.

If $A$ and $D$ are maximal in $P_1$, then by Lemma \ref{lemma:localdec1} and Lemma \ref{lemma:localdec2}, it follows that $F$ or $E$ cannot be maximal in $P_2$, since then $E$ would become maximal-extremal in $P$. It follows that either $B$ or $C$ is maximal-extremal in $P_2$, but not both since we cannot have to maximal-extremal vertices next to each other by our generic assumption. So, we cannot have such a case.

If $A$ and $E$ are maximal-extremal in $P_1$, it follows by Lemma \ref{lemma:localdec2} that $F$ cannot be maximal-extremal in $P_2$, since it would follow that $E$ would be maximal in $P$. It follows that only $C$ and $E$ can be maximal-extremal in $P_2$. By Lemma \ref{lemma:localdec3}, it follows that $E$ cannot be maximal-extremal in $P_2$. So, such a situation cannot happen.

If $B$ and $D$ are maximal-extremal in $P_1$, it then follows by Lemma \ref{lemma:localdec1} and Lemma \ref{lemma:localdec2} that $E$ or $F$ cannot be maximal-extremal in $P_2$, since it would follow that $E$ would be maximal-extremal in $P$. This simply leaves us simply with $B$ and $C$, which both cannot be extremal in $P_2$ by our generic assumption. So, such a situation is not plausible.

\vskip2mm
\noindent\textit{Subcase 2:} $C$ is maximal in $P$.
\vskip1mm
As before, it is sufficient to check the same three configurations in $P_1$.

If $A$ and $D$ are maximal-extremal in $P_1$, then by Lemma \ref{lemma:localdec1} and Lemma \ref{lemma:localdec2} it follows that $F$ or $E$ cannot be maximal-extremal in $P_2$, because if they were, it would follow that $E$ is maximal-extremal in $P$. So, by our generic assumption it follows that we cannot have two maximal-extremal vertices in $P_2$.

If $A$ and $E$ are maximal-extremal in $P_1$, by Lemma \ref{lemma:localdec1} it follows that $F$ cannot be maximal-extremal in $P_2$ since, if this were the case, $E$ would become maximal-extremal in $P$. By Lemma \ref{lemma:localdec3}, it follows that $E$ cannot be maximal-extremal in $P_2$ since $A$ would be maximal-extremal in $P$. So, it follows by our generic assumption that we cannot have two maximal-extremal vertices in $P_2$.

Lastly, if $B$ and $D$ are maximal-extremal in $P_1$, it then follows by Lemma \ref{lemma:localdec1} and Lemma \ref{lemma:localdec2} that $E$ or $F$ cannot be maximal-extremal in $P_2$, since $E$ would be maximal-extremal in $P$. Again, we cannot have two maximal-extremal vertices in $P_2$.
\end{proof}

\vskip2mm
\noindent\textit{Case 3:} We gain two maximal-extremal vertices in $P_1$ and gain one maximal-extremal vertex in $P_2$, and none of the six vertices is maximal-extremal in $P$.
\vskip1mm

By the symmetry of our cases, it is just sufficient to check two configurations in $P_1$ and see that we cannot gain any new vertices in $P_2$.

For the first case, we assume that $A$ and $D$ are maximal-extremal in $P_1$. It then follows by Lemma \ref{lemma:localdec1} and Lemma \ref{lemma:localdec2} that neither $C$, $B$, $E$ or $F$ can be maximal-extremal for $P_2$, because if any of those four vertices were maximal-extremal in $P_2$, then we risk gaining a maximal-extremal vertex for $P$.

Secondly, assume that $A$ and $E$ are maximal-extremal in $P_1$. It follows from Lemma \ref{lemma:localdec1} and Lemma \ref{lemma:localdec2} that $F$, $C$ and $B$ cannot be maximal-extremal for $P_2$. Lemma \ref{lemma:localdec3} eliminates the possibility of $E$ being maximal-extremal in $P_2$, showing that this case is not feasible.
\vskip2mm
We now see that we have eliminated all possible cases that violate our inequality. Due to the nice symmetry of our cases, we only had three major situations to check. So, our assertion is proved.

\subsection{New Proofs of The Four-Vertex Theorems}

Now that we have our new results in the previous two sections, we will show that we can actually can derive our Four-Vertex Theorems from Section 2.

\subsubsection{Deriving the Global Four-Vertex Theorem}

Before we can prove our Global Four-Vertex Theorem, we need to prove an important fact about the triangulation of polygons.

\begin{lemma}
\label{lemma:triangulation}
Let $P$ be a convex polygon with seven or more vertices and let $T(P)$ be a triangulation of $P$. Then, there exists a diagonal of our triangulation such that if we apply a decomposition of $P$ using this diagonal, then both $P_1$ and $P_2$ have four or more vertices.
\end{lemma}
\begin{proof}
We will perform induction on the number of vertices. For the base case $n=7$, it is just a simple routine checking of all possible ways to triangulate $P$. We see by checking all possibilities, there always exists a diagonal which divides the polygon into two polygons with four or more vertices.

Now we consider the case when $n>7$. We pick a vertex and call it $V_i$, and let $V_{i-1}$ and $V_{i+1}$ be the neighboring vertices of $V_i$. Now, we remove our vertex $V_i$ and connect vertices $V_{i-1}$ and $V_{i+1}$ by an edge to obtain a polygon $P'$ with $n-1$ vertices. By induction, there exists a diagonal $d$ such that if we decompose $P'$ into $P_1$ and $P_2$ by this diagonal, then $P_1$ and $P_2$ each have four or more vertices. We now have two cases.
\vskip2mm
\noindent\textit{Case 1:} Our diagonal $d$ does not pass through $V_{i-1}$ or $V_{i+1}$.
\vskip1mm
Since our diagonal $d$ does not pass through $V_{i-1}$ or $V_{i+1}$, it follows that either $P_1$ or $P_2$ posseses both $V_{i-1}$ and $V_{i+1}$. Without loss of generality, assume that $P_1$ possesses the vertices $V_{i-1}$ and $V_{i+1}$. Since $P_1$ has four or more vertices, it follows that if we add $V_i$ back to our polygon, then we add $V_i$ to $P_1$. So, our assertion is proved.
\vskip2mm
\noindent\textit{Case 2:} Our diagonal $d$ passes through $V_{i-1}$ or $V_{i+1}$.
\vskip1mm
Without loss of generality, we assume that it passes through $V_{i-1}$. It follows that either $P_1$ or $P_2$ possesses the vertex $V_{i+1}$, so without loss of generality we assume that $P_2$ possesses $V_{i+1}$. It follows that if we re-attach our vertex $V_i$ to $P'$, then we are then attaching it to $P_2$. Since $P_2$ already had four or more vertices by the inductive assumption, our assertion follows.
\end{proof}

\begin{defn}
Let $P$ be a polygon and let $T(P)$ be a triangulation of $P$ such that every edge and diagonal is Delaunay. We call $T(P)$ a Delaunay triangulation. Similarly, if every edge and diagonal is Anti-Delaunay, then we call $T(P)$ an Anti-Delaunay triangulation.
\end{defn}

\begin{remark}
We note that these definitions are equivalent to the definitions of Delaunay and Anti-Delaunay triangulations given in Section 2.1.1.
\end{remark}

It is interesting to note that Lemma \ref{lemma:triangulation} holds for any triangulation of $P$, in particular for both Delaunay triangulation and Anti-Delaunay triangulation. The question is, why did Lemma \ref{lemma:triangulation} not hold for the case when $n=6$? The following figure illustrates a triangulation for which the case $n=6$ fails:

\begin{figure}[H]
\centerline{\includegraphics[scale=0.47]{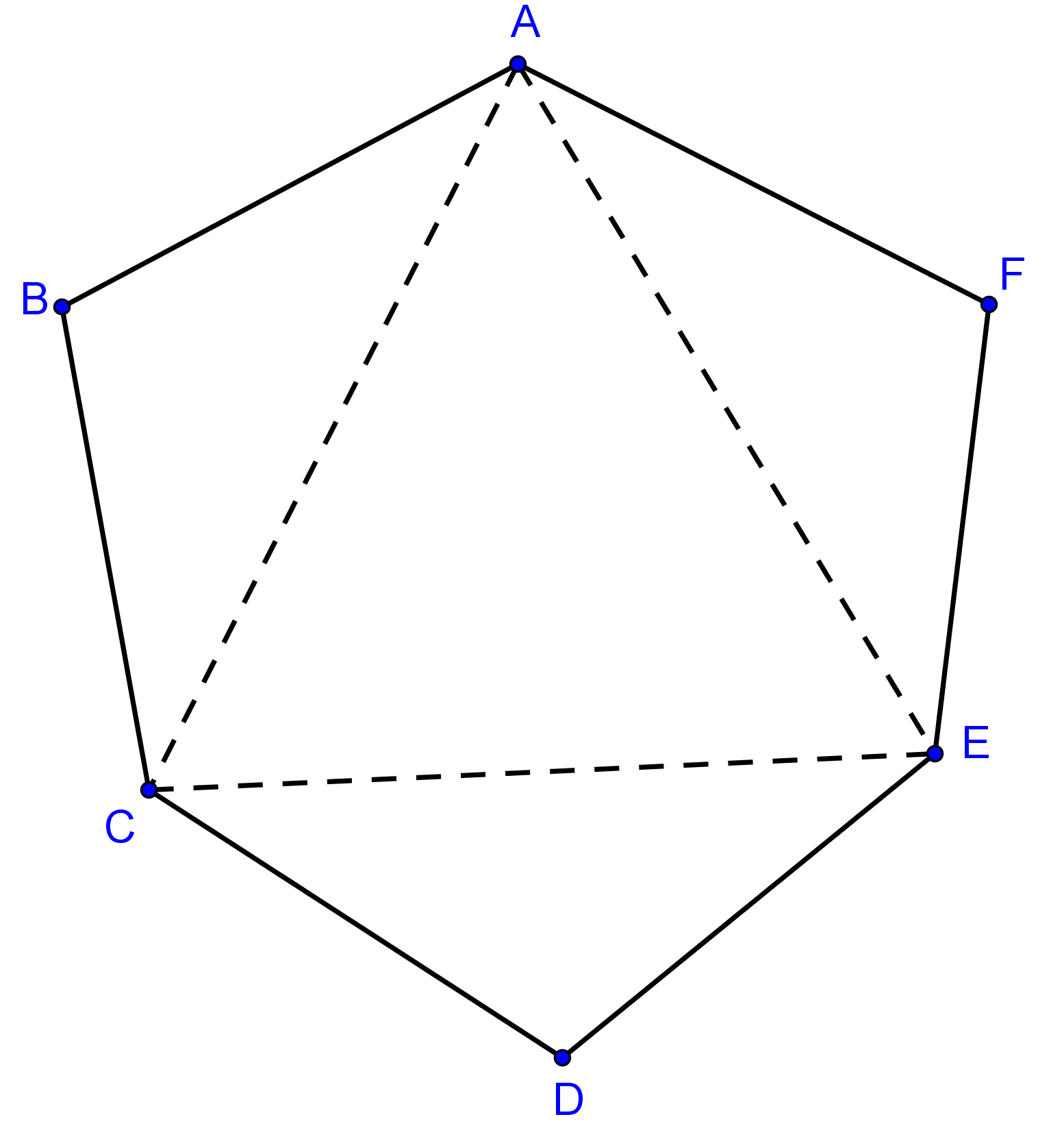}}
\caption[An Example Where Lemma \ref{lemma:triangulation} Fails]{}
\end{figure}

\begin{theorem}[The Global Four-Vertex Theorem]
\label{thm:newglobal}
Let $P$ be a generic convex polygon with six or more vertices. Then $$s_{+}(P)+s_{-}(P)\geq 4.$$
\end{theorem}
\begin{proof}
We now will prove our theorem by induction on the number of vertices. For the base case $n=6$, if we apply a decomposition to $P$, it follows that $P_1$ and $P_2$ are both quadrilaterals. By Proposition \ref{prop:quad}, we obtain that $P_1$ and $P_2$ each have four globally extremal vertices. It follows from Lemma \ref{lemma:decompglobal} that $P$ has four globally extremal vertices.

We now consider the case where $n\geq7$. We begin by applying a Delaunay triangulation to $P$. By Lemma \ref{lemma:triangulation}, it follows that there exists a diagonal $d$ such that when we decompose $P$ by this diagonal, $P_1$ and $P_2$ each have four or more vertices. Since we applied a Delaunay triangulation, it follows that $d$ is Delaunay. By Theorem \ref{theorem:strongglobal}, we have that $s_{-}(P)\geq s_{-}(P_1)+s_{-}(P_2)-2$. Since $P_1$ and $P_2$ have less vertices than $P$, by induction it follows that $s_{-}(P_1)\geq 2$ and $s_{-}(P_2)\geq 2$. Applying this to Theorem \ref{theorem:strongglobal}, we get that $s_{-}(P)\geq 2$.

By applying an Anti-Delaunay triangulation to our polygon, we can apply Lemma \ref{lemma:triangulation} and analogously, by Theorem \ref{theorem:strongglobal2}, we obtain that $s_{+}(P)\geq s_{+}(P_1)+s_{+}(P_2)-2$.  Applying our induction hypothesis to $P_1$ and $P_2$, we obtain that $s_{+}(P)\geq 2$. So we have that $s_{+}(P)+s_{-}(P)\geq 4$, proving the theorem.
\end{proof}

\subsubsection{Deriving the Local Four-Vertex Theorem}

In this section, we will derive the Local Four-Vertex Theorem from our results.

\begin{theorem}[The Local Four-Vertex Theorem]
\label{thm:newlocal}
Let $P$ be a generic convex polygon with at least six vertices. Then
$$l_{+}(P)+l_{-}(P)\geq 4.$$
\end{theorem}
\begin{proof}
We apply induction on the number of vertices of $P$. For the case where $n=6$, we know that if we apply a decomposition to $P$, then both $P_1$ and $P_2$ will be quadrilaterals. Proposition \ref{prop:quad} yields that $l_{-}(P_1)=l_{-}(P_2)=2$. Applying this to Theorem \ref{thm:localineq} completes the proof for this case.

Now, for the inductive step, by Theorem \ref{thm:localineq}, we have that $l_{-}(P)\geq l_{-}(P_1)+l_{-}(P_2)-2$. We now apply induction to the smaller polygons $P_1$ and $P_2$ and obtain that $l_{-}(P_1)\geq 2$ and $l_{-}(P_2)\geq 2$. Applying this to our inequality and using Proposition \ref{prop:maxmin}, we prove our assertion.
\end{proof}

\section{Appendix}

\subsection{Our Counterexamples}
In this section we will provide the coordinates of our important counterexamples. In the following matrices, the first column always will have the coordinates of the vertex $A$, and the rest of the coordinates are given by traveling counterclockwise along the polygon.
\vskip2mm
Figure 2.7:
$$\left(
  \begin{array}{cccccc}
    18.38 & 17.59 & 13.58 & 26.21 & 23.68 & 21.88 \\
    -2.05 & -2.41 & -6.13 & -5.82 & -3.54 & -2.9 \\
  \end{array}
\right)$$
\vskip1mm
Figure 4.1:
$$\left(
  \begin{array}{cccccccccccc}
    1.46 & -2.19 & -2.79 & -2.74 & -1.48 & 1.54 & 4.72 & 6.57 & 7.78 & 8.34 & 6.53 & 4.44 \\
    5.59 & 5.17 & 2.55 & -0.49 & -2.08 & -2.72 & -2.04 & -0.62 & 0.84 & 2.39 & 4.01 & 5.22 \\
  \end{array}
\right)$$
\vskip1mm
Figure 4.2:
$$\left(
  \begin{array}{cccccccccccc}
    1.78 & 1.24 & 0.37 & 1 & 1.32 & 1.82 & 2.48 & 3 & 3.36 & 3.45 & 3.32 & 2.44 \\
    4.76 & 4.58 & 3.77 & 2.23 & 1.86 & 1.7 & 1.7 & 2 & 2.41 & 3.08 & 4.3 & 4.68 \\
  \end{array}
\right)$$
\vskip1mm
Figure 4.6:

{\tiny
$$\left(
  \begin{array}{ccccccccccccccc}
    0.6 & -0.98 & -1.82 & -1.85 & -1.12 & 0.62 & 1.63 & 2.23 & 2.68 & 3.24 & 3.52 & 3.52 & 3.24 & 2.15 & 1.51 \\
    5.12 & 4.08 & 2.39 & 0.52 & -1.74 & -3.44 & -3.29 & -2.53 & -1.35 & 0.23 & 1.28 & 1.86 & 3.21 & 4.32 & 4.98 \\
  \end{array}
\right)$$
}

\subsection{MATLAB Code}

\subsubsection{Finding the Evolute}

The following code creates a function which returns the center of the circumcircle of three points in the plane given in matrix form:

\begin{verbatim}
function x=circumc(p);
A=2*[p(1,1)-p(1,3), p(2,1)-p(2,3);p(1,2)-p(1,3), p(2,2)-p(2,3)];
b=[p(1,1)^2+p(2,1)^2-p(1,3)^2-p(2,3)^2;p(1,2)^2+p(2,2)^2-p(1,3)^2-p(2,3)^2];
x=A\b;
x=x';
\end{verbatim}
\vskip2mm

The following code creates a function which returns the coordinates of the evolute in matrix form, given a polygonal curve in matrix form (Note that this uses the previous function):

\begin{verbatim}
function x=evolute(p);
n=size(p,2);
q=[p(:,n),p,p(:,1)];
x=[1:0];
for i=1:n
    r=q(:,i:i+2);
    x=[x,circumc(r)']
end
\end{verbatim}

\subsubsection{Finding The Winding Number of a Polygon}

The following code creates a function which, after inputting coordinates of the polygon in matrix form, returns the measures of the angles at each of the vertices. It always measures the left angle with respect to orientation to be consistent with Definition 2.1. If the input matrix is $2\times n$, then it returns a $1\times n$ matrix.

\begin{verbatim}
function q1=vertexangle(p);
n=size(p,2);
q=[p(:,n),p,p(:,1)];
pangle=[1:0];
for i=1:n
    s1=[q(:,i)-q(:,i+1)];
    s2=[q(:,i+2)-q(:,i+1)];
    m1=sqrt(s1(1,1)^2+s1(2,1)^2);
    m2=sqrt(s2(1,1)^2+s2(2,1)^2);
    if det([s1,s2]) < 0
       ang=acos(dot(s1',s2')/(m1*m2));
       pangle=[pangle,ang];
    elseif  det([s1,s2])>0
       ang=2*pi-acos(dot(s1',s2')/(m1*m2));
       pangle=[pangle,ang];
    end
end
q1=pangle
\end{verbatim}
\vskip2mm
The following code defines a function which returns the winding number of a polygon. As input, it accepts the coordinates of the vertices of the polygon in matrix form. You must define the previous function to use this function.

\begin{verbatim}
function w=windingnumber(p);
q=vertexangle(p);
n=size(q,2);
w=0;
for i=1:n
    w=w+(pi-q(:,i));
end
w=w/(2*pi)
\end{verbatim}
\listoffigures

\section{Acknowledgements}

The author would like to thank his advisor Dr. Oleg R. Musin for his guidance and support.

Wiktor Mogilski, Department of Mathematics, University of Texas at Brownsville, 80 Fort Brown, Brownsville, TX, 78520.

{\it E-mail address:} wiktormd@gmail.com

\end{document}